\documentclass[a4paper]{amsart}
\usepackage{amssymb}
\usepackage{amsmath}
\usepackage{amsthm}
\usepackage{float}
\usepackage{latexsym}
\usepackage{graphicx}
\usepackage{color}

\newtheorem{theorem}{Theorem}

\newtheorem{corollary}[theorem]{Corollary}

\newtheorem{definition}[theorem]{Definition}

\newtheorem{lemma}[theorem]{Lemma}

\newtheorem{remark}[theorem]{Remark}

\numberwithin{equation}{section}
\newcommand{\kom}[1]{}
\renewcommand{\kom}[1]{{\bf [#1]}}
\definecolor{blau}{rgb}{0.1,0.0,0.9}

\newcounter{komcounter}
\numberwithin{komcounter}{section}

\begin{document}
\title[SDEs and PDEs on non-smooth time-dependent domains]{Stochastic and
partial differential equations on non-smooth time-dependent domains}
\address{Niklas L.P. Lundstr\"{o}m\\
Department of Mathematics and Mathematical Statistics, Ume\aa\ University\\
SE-901 87 Ume\aa , Sweden}
\email{niklas.lundstrom@umu.se}
\address{Thomas \"{O}nskog\\
Department of Mathematics, Royal Institute of Technology (KTH)\\
SE-100 44 Stockholm, Sweden}
\email{onskog@kth.se}
\author{Niklas L.P. Lundstr\"{o}m, Thomas \"{O}nskog}

\begin{abstract}
In this article, we consider non-smooth time-dependent domains whose boundary is $\mathcal{W}^{1,p}$ in time and single-valued, smoothly varying directions of reflection at the boundary. In this setting, we first prove existence and uniqueness of strong solutions to stochastic differential equations with oblique reflection. Secondly, we prove, using the theory of viscosity solutions, a comparison principle for fully nonlinear second-order parabolic partial differential equations with oblique derivative boundary conditions. As a consequence, we obtain uniqueness, and, by barrier construction and Perron's method, we also conclude existence of viscosity
solutions. Our results generalize two articles by Dupuis and Ishii to time-dependent domains.

\vspace{0.2cm}%

\noindent 2000\textit{\ Mathematics Subject Classification. }35D05, 49L25,
60J50, 60J60.
\vspace{0.2cm}%

\end{abstract}

\maketitle

\setcounter{equation}{0} \setcounter{theorem}{0}

\section{Introduction\label{intro}}

In this article we establish existence and uniqueness of strong solutions to
stochastic differential equations (SDE) with single-valued, smoothly varying
oblique reflection at the boundary of a bounded, non-smooth time-dependent
domain whose boundary is $\mathcal{W}^{1,p}$ in time. In the same geometric
setting, we also prove a comparison principle, uniqueness and existence of
viscosity solutions to partial differential equations (PDE) with oblique
derivative boundary conditions.

In the SDE case, our approach is based on the Skorohod problem, which, in
the form studied in this article, was first described by Tanaka \cite%
{Tanaka1979}. Tanaka established existence and uniqueness of solutions to
the Skorohod problem in convex domains with normal reflection. These results
were subsequently substantially generalized by, in particular, Lions and
Sznitman \cite{LionsSznitman1984} and Saisho \cite{Saisho1987}. To the
authors' knowledge, the most general results on strong solutions to
reflected SDEs in time-independent domains based on the Skorohod problem are
those established by Dupuis and Ishii \cite{DupuisIshii1993}. The aim here
is to generalize the SDE results mentioned above, in particular those of
Case 1 in \cite{DupuisIshii1993}, to the setting of time-dependent domains.

There is, by now, a number of articles on reflected SDEs in time-dependent
domains. Early results on this topic include the exhaustive study of the
heat equation and reflected Brownian motion in smooth time-dependent domains
by Burdzy, Chen, and Sylvester \cite{BurdzyChenSylvester2004AP} and the
study of reflected SDEs in smooth time-dependent domains with reflection in
the normal direction by Costantini, Gobet, and El Karoui \cite%
{CostantiniGobetKaroui2006}. We also mention that Burdzy, Kang, and Ramanan
\cite{BurdzyKangRamanan2009} investigated the Skorohod problem in a
one-dimensional, time-dependent domain and, in particular, found conditions
for when there exists a solution to the Skorohod problem in the event that
the two boundaries meet. Existence of weak solutions to SDEs with oblique
reflection in non-smooth time-dependent domains was established by Nystr\"{o}%
m and \"{O}nskog \cite{NystromOnskog2010a} under fairly general conditions
using the approach of \cite{Costantini1992}. In the article at hand, we use
the approach of \cite{DupuisIshii1993} and derive regularity conditions,
under which we can obtain existence and also uniqueness of strong solutions
to SDEs with oblique reflection in time-dependent domains.

Turning to the PDE case, we recall that the approach of \cite%
{DupuisIshii1993} relies on the construction of test functions used earlier
in Dupuis and Ishii \cite{DupuisIshii1990} to prove the comparison
principle, existence and uniqueness for fully nonlinear second-order
elliptic PDEs in non-smooth time-independent domains. Here we generalize
these test functions to our time-dependent setting, and obtain the
corresponding results for both SDEs and PDEs in time-dependent domains. In
particular, our PDE results generalize the main part of \cite%
{DupuisIshii1990} to hold in the setting of fully nonlinear second-order
parabolic PDEs in non-smooth time-dependent domains. Our proofs are based on
the theory of viscosity solutions. The first step is to observe that the
maximum principle for semicontinuous functions by Crandall and Ishii \cite%
{CrandallIshii1990} holds in time-dependent domains. Using the maximum
principle and the above-mentioned test functions, we prove the comparison
principle by following the nowadays standard method, see Crandall, Ishii,
and Lions \cite{CrandallIshiiLions1992} and \cite{DupuisIshii1990}. Next, we
prove existence of a unique solution to the PDE problem by means of Perron's
method, the comparison principle and by constructing several explicit sub-
and supersolutions (barriers) to the PDE.

To the authors' knowledge, there are no previous results on the oblique
derivative problem for parabolic PDEs in non-smooth time-dependent domains.
For time-independent domains, however, there are several articles in the
literature. Besides \cite{DupuisIshii1990}, Dupuis and Ishii studied oblique
derivative problems for fully nonlinear elliptic PDEs on domains with
corners in \cite{DupuisIshii1991}. Moreover, Barles \cite{Barles1993} proved
a comparison principle and existence of unique solutions to degenerate
elliptic and parabolic boundary value problems with nonlinear Neumann type
boundary conditions in bounded domains with $\mathcal{W}^{3,\infty }$%
-boundary. Ishii and Sato \cite{IshiiSato2004} proved similar theorems for
boundary value problems for some singular degenerate parabolic partial
differential equations with nonlinear oblique derivative boundary conditions
in bounded $\mathcal{C}^{1}$-domains. Further, in bounded domains with $%
\mathcal{W}^{3,\infty }$-boundary, Bourgoing \cite{Bourgoing2008} considered
singular degenerate parabolic equations and equations having $L^{1}$
dependence in time.

Concerning PDEs in the setting of time-dependent domains, we mention that Bj%
\"{o}rn\textit{\ et al}.~\cite{BjornBjornGianazzaParviainen2015} proved,
among other results, a comparison principle for solutions of degenerate and
singular parabolic equations with Dirichlet boundary conditions using a
different technique and that Avelin \cite{Avelin2016} proved boundary
estimates of solutions to the degenerate $p$-parabolic equation.

As a motivation for considering SDEs and PDEs in time-dependent domains, we
mention that such geometries arise naturally in a wide range of applications
in which the governing equation of interest is a differential equation, for
example in modelling of crack propagation \cite{NicaiseSandig2007},
modelling of fluids \cite{FiloZauskova2008}, \cite{HeHsiao2000} and
modelling of chemical, petrochemical and pharmaceutical processes \cite%
{IzadiAbdollahiDubljevic2014}.

The rest of the paper is organized as follows. In Section \ref{DNA} we give
preliminary definitions, notations, assumptions and also state our main
results. In Section \ref{test} we construct the test functions crucial for
the proofs of both the SDE and the PDE results. Using these test functions,
we prove existence of solutions to the Skorohod problem in Section \ref{SP}.
The results on the Skorohod problem are subsequently used, in Section \ref%
{RSDE}, to prove the main results for SDEs. Finally, in Section \ref{PDE},
we use the theory of viscosity solutions together with the test functions
derived in Section \ref{test} to establish the PDE results.

\setcounter{equation}{0} \setcounter{theorem}{0}

\section{Preliminaries and statement of main results\label{DNA}}

Throughout this article we will use the following definitions and
assumptions. Given $n\geq 1$, $T>0$ and a bounded, open, connected set $%
\Omega ^{\prime }\subset
\mathbb{R}
^{n+1}$ we will refer to
\begin{equation}
\Omega =\Omega ^{\prime }\cap ([0,T]\times
\mathbb{R}
^{n}),  \label{timedep}
\end{equation}%
as a time-dependent domain. Given $\Omega $ and $t\in \left[ 0,T\right] $,
we define the time sections of $\Omega $ as $\Omega _{t}=\left\{ x:\left(
t,x\right) \in \Omega \right\} $, and we assume that%
\begin{equation}
\Omega _{t}\neq \emptyset \text{ and that }\Omega _{t}\text{ is bounded and
connected for every }t\in \left[ 0,T\right] .  \label{timesect}
\end{equation}%
Let $\partial \Omega _{t}$, for $t\in \left[ 0,T\right] $, denote the
boundary of $\Omega _{t}$. Let $\left\langle \cdot ,\cdot \right\rangle $
and $\left\vert \cdot \right\vert =\left\langle \cdot ,\cdot \right\rangle
^{1/2}$ define the Euclidean inner product and norm, respectively, on $%
\mathbb{R}
^{n}$ and define, whenever $a\in
\mathbb{R}
^{n}$ and $\,b>0$, the sets $B\left( a,b\right) =\left\{ x\in
\mathbb{R}
^{n}:\left\vert x-a\right\vert \leq b\right\} $ and $S\left( a,b\right)
=\left\{ x\in
\mathbb{R}
^{n}:\left\vert x-a\right\vert =b\right\} $. For any Euclidean spaces $E$
and $F$, we define the following spaces of functions mapping $E$ into $F$. $%
\mathcal{C}\left( E,F\right) $ denotes the set of continuous functions, $%
\mathcal{C}^{k}\left( E,F\right) $ denotes the set of $k$ times continuously
differentiable functions and $\mathcal{W}^{1,p}\left( E,F\right) $ denotes
the Sobolev space of functions whose first order weak derivatives belong to $%
L^{p}\left( E\right) $. If we can distinguish the time variable from the
spatial variables, we let $\mathcal{C}^{1,2}\left( E,F\right) $ denote the
set of functions, whose elements are continuously differentiable once with
respect to the time variable and twice with respect to any space variable,
and by $\mathcal{C}_{b}^{1,2}\left( E,F\right) $ we denote the space of
bounded functions in $\mathcal{C}^{1,2}\left( E,F\right) $ having bounded
derivatives. Moreover, $\mathcal{BV}\left( E,F\right) $ denotes the set of
functions with bounded variation. In particular, for $\eta \in \mathcal{BV}%
\left( \left[ 0,T\right] ,%
\mathbb{R}
^{n}\right) $, we let $\left\vert \eta \right\vert \left( t\right) $ denote
the total variation of $\eta $ over the interval $\left[ 0,t\right] $.

\subsection{Assumptions on the domain and directions of reflection\label%
{geoassume}}

Throughout this article we consider non-smooth time-dependent domains of the
following type. Let $\Omega \subset
\mathbb{R}
^{n+1}$ be a time-dependent domain satisfying \eqref{timesect}. The
direction of reflection at $x\in \partial \Omega _{t}$, $t\in \left[ 0,T%
\right] $, is given by $\gamma \left( t,x\right) $ satisfying
\begin{equation}
\gamma \in \mathcal{C}_{b}^{1,2}\left( \mathbb{R}^{n+1},B\left( 0,1\right)
\right) ,  \label{smooth_gamma}
\end{equation}%
such that $\gamma \left( t,x\right) \in S\left( 0,1\right) $ for all $\left(
t,x\right) \in V$, where $V$ is an open set satisfying $\Omega
_{t}^{c}\subset V$ for all $t\in \lbrack 0,T]$. Moreover, there is a
constant $\rho \in \left( 0,1\right) $ such that the exterior cone condition%
\begin{equation}
\bigcup_{0\leq \zeta \leq \rho }B\left( x-\zeta \gamma \left( t,x\right)
,\zeta \rho \right) \subset \Omega _{t}^{c},  \label{boundarylip}
\end{equation}%
holds, for all $x\in \partial \Omega _{t}$, $t\in \left[ 0,T\right] $. Note
that it follows from \eqref{boundarylip} that $\gamma $ points into the
domain and this is indeed the standard convention for SDEs. For PDEs,
however, the standard convention is to let $\gamma $ point out of the
domain. To facilitate for readers accustomed with either of these
conventions we, in the following, let $\gamma $ point inward whenever SDEs
are treated, whereas when we treat PDEs we assume the existence of a
function
\begin{equation}
\widetilde{\gamma }\in \mathcal{C}_{b}^{1,2}\left( \mathbb{R}^{n+1},B\left(
0,1\right) \right) ,  \label{smooth_gamma2}
\end{equation}%
defined as $\widetilde{\gamma }\left( t,x\right) =-\gamma \left( t,x\right) $%
, with $\gamma $ as in \eqref{smooth_gamma}. In particular, we have
\begin{equation}
\bigcup_{0\leq \zeta \leq \rho }B\left( x+\zeta \widetilde{\gamma }\left(
t,x\right) ,\zeta \rho \right) \subset \Omega _{t}^{c},  \label{boundarylip2}
\end{equation}%
for all $x\in \partial \Omega _{t}$, $t\in \left[ 0,T\right] $. Finally,
regarding the temporal variation of the domain, we define $d\left(
t,x\right) =\inf_{y\in \Omega _{t}}\left\vert x-y\right\vert $, for all $%
t\in \left[ 0,T\right] $, $x\in
\mathbb{R}
^{n}$, and assume that for some fixed $p\in \left( 1,\infty \right) $ and
all $x\in
\mathbb{R}
^{n}$,
\begin{equation}
d\left( \cdot ,x\right) \in \mathcal{W}^{1,p}\left( \left[ 0,T\right] ,\left[
0,\infty \right) \right) ,  \label{templip}
\end{equation}%
with Sobolev norm uniformly bounded in space. We also assume that $%
D_{t}d(t,x)$ is jointly measurable in $(t,x)$.

\begin{remark}
\label{spaceremark}A simple contradiction argument based on the exterior
cone condition \eqref{boundarylip} for the time sections and the regularity
of $\gamma $ and $\Omega _{t}$, shows that the time sections satisfy the
interior cone condition
\begin{equation*}
\bigcup_{0\leq \zeta \leq \rho }B\left( x+\zeta \gamma \left( t,x\right)
,\zeta \rho \right) \subset \overline{\Omega }_{t},
\end{equation*}%
for all $x\in \partial \Omega _{t}$, $t\in \left[ 0,T\right] $. The exterior
and interior cone conditions together imply that the boundary of $\Omega
_{t} $ is Lipschitz continuous (in space) with a Lipschitz constant $K_{t}$
satisfying $\sup_{t\in \left[ 0,T\right] }K_{t}<\infty $. Moreover, these
conditions imply that for a suitable constant $\theta \in \left( 0,1\right) $%
, $\theta ^{2}>1-\rho ^{2}$, there exists $\delta >0$ such that%
\begin{equation*}
\left\langle y-x,\gamma \left( t,x\right) \right\rangle \geq -\theta
\left\vert y-x\right\vert ,
\end{equation*}%
for all $x\in \partial \Omega _{t}$, $y\in \overline{\Omega }_{t}$, $t\in %
\left[ 0,T\right] $ satisfying $\left\vert x-y\right\vert \leq \delta $.
\end{remark}

\begin{remark}
\label{timeholder}By Morrey's inequality, condition \eqref{templip} implies
the existence of a H\"{o}lder exponent $\widehat{\alpha }=1-1/p\in \left(
0,1\right) $ and a H\"{o}lder constant $K\in \left( 0,\infty \right) $ such
that, for all $s,t\in \left[ 0,T\right] $, $x\in
\mathbb{R}
^{n}$,
\begin{equation}
\left\vert d\left( s,x\right) -d\left( t,x\right) \right\vert \leq
K\left\vert s-t\right\vert ^{\widehat{\alpha }}.  \label{tempholder}
\end{equation}
\end{remark}

\begin{remark}
The assumptions imposed on the time sections of the time-dependent domain in %
\eqref{smooth_gamma}, \eqref{boundarylip} coincide with those imposed on the
time-independent domains in \cite{DupuisIshii1990} and in Case 1 of \cite%
{DupuisIshii1993}. For time-independent domains, existence and uniqueness
results for SDE and PDE have also been obtained under the conditions given
in \cite{DupuisIshii1991a} and in Case 2 of \cite{DupuisIshii1993}. It is
likely that also these results can be extended to time-dependent domains
using a procedure similar to that of the article at hand, but we leave this
as a topic of future research.
\end{remark}

\begin{remark}
Consider the function%
\begin{equation*}
l\left( r\right) =\sup_{s,t\in \left[ 0,T\right] ,\text{ }\left\vert
s-t\right\vert \leq r}\,\,\sup_{x\in \overline{\Omega }_{s}}\,\,\inf_{y\in
\overline{\Omega }_{t}}\left\vert x-y\right\vert ,
\end{equation*}%
introduced in \cite{CostantiniGobetKaroui2006} and frequently used in \cite%
{NystromOnskog2010a}. Condition \eqref{tempholder} is equivalent to,
\begin{equation*}
l\left( r\right) \leq Kr^{\widehat{\alpha }},
\end{equation*}%
which is considerably stronger than the condition $\lim_{r\rightarrow
0^{+}}l\left( r\right) =0$ assumed in \cite{NystromOnskog2010a}. On the
other hand, it was assumed in \cite{NystromOnskog2010a} that $\Omega _{t}$
satisfies a uniform exterior sphere condition, and this does not hold in
general for domains satisfying \eqref{boundarylip}.
\end{remark}

\subsection{Statement of main result for SDEs}

We consider the Skorohod problem in the following form.

\begin{definition}
\label{skorohodprob}Given $\psi \in \mathcal{C}\left( \left[ 0,T\right] ,%
\mathbb{R}
^{n}\right) $, with $\psi \left( 0\right) \in \overline{\Omega }_{0}$, we
say that the pair $\left( \phi ,\lambda \right) \in \mathcal{C}\left( \left[
0,T\right] ,%
\mathbb{R}
^{n}\right) \times \mathcal{C}\left( \left[ 0,T\right] ,%
\mathbb{R}
^{n}\right) $ is a solution to the Skorohod problem for $\left( \Omega
,\gamma ,\psi \right) $ if $\left( \psi ,\phi ,\lambda \right) $ satisfies,
for all $t\in \left[ 0,T\right] $,
\begin{eqnarray}
\phi \left( t\right) &=&\psi \left( t\right) +\lambda \left( t\right) ,\quad
\phi \left( 0\right) =\psi \left( 0\right) ,  \label{SP1} \\
\phi \left( t\right) &\in &\overline{\Omega }_{t},  \label{SP2} \\
\left\vert \lambda \right\vert \left( T\right) &<&\infty ,  \label{SP3} \\
\left\vert \lambda \right\vert \left( t\right) &=&\int_{\left( 0,t\right]
}I_{\left\{ \phi \left( s\right) \in \partial \Omega _{s}\right\}
}d\left\vert \lambda \right\vert \left( s\right) ,  \label{SP4} \\
\lambda \left( t\right) &=&\int_{\left( 0,t\right] }\widehat{\gamma }\left(
s\right) d\left\vert \lambda \right\vert \left( s\right) ,  \label{SP5}
\end{eqnarray}%
for some measurable function $\widehat{\gamma }:\left[ 0,T\right]
\rightarrow
\mathbb{R}
^{n}$ satisfying $\widehat{\gamma }\left( s\right) =\gamma \left( s,\phi
\left( s\right) \right) $ $d\left\vert \lambda \right\vert $-a.s.
\end{definition}

We use the Skorohod problem to construct solutions to SDEs confined to the
given time-dependent domain $\overline{\Omega }$ and with direction of
reflection given by $\gamma $. We shall consider the following notion of
SDEs. Let $\left( \Omega ,\mathcal{F},\mathbb{P}\right) $ be a complete
probability space and let $\left\{ \mathcal{F}_{t}\right\} _{t\geq 0}$ be a
filtration satisfying the usual conditions. Let $m$ be a positive integer,
let $W=\left( W_{i}\right) $ be an $m$-dimensional Wiener process and let $b:%
\left[ 0,T\right] \times
\mathbb{R}
^{n}\rightarrow
\mathbb{R}
^{n}$ and $\sigma :\left[ 0,T\right] \times
\mathbb{R}
^{n}\rightarrow
\mathbb{R}
^{n\times m}$ be continuous functions.

\begin{definition}
\label{strong}A strong solution to the SDE in $\overline{\Omega }$ driven by
the Wiener process $W$ and with coefficients $b$ and $\sigma $, direction of
reflection along $\gamma $ and initial condition $x\in \overline{\Omega }%
_{0} $ is an $\left\{ \mathcal{F}_{t}\right\} $-adapted continuous
stochastic process $X\left( t\right) $ which satisfies, $\mathbb{P}$-almost
surely, whenever $t\in \left[ 0,T\right] $,
\begin{equation}
X\left( t\right) =x+\int_{0}^{t}b\left( s,X\left( s\right) \right)
ds+\int_{0}^{t}\left\langle \sigma \left( s,X\left( s\right) \right)
,dW\left( s\right) \right\rangle +\Lambda \left( t\right) ,  \label{RSDE1}
\end{equation}%
where%
\begin{equation}
X\left( t\right) \in \overline{\Omega }_{t},\quad \left\vert \Lambda
\right\vert \left( t\right) =\int_{\left( 0,t\right] }I_{\left\{ X\left(
s\right) \in \partial \Omega _{s}\right\} }d\left\vert \Lambda \right\vert
\left( s\right) <\infty ,  \label{RSDE2}
\end{equation}%
and where
\begin{equation}
\Lambda \left( t\right) =\int_{\left( 0,t\right] }\widehat{\gamma }\left(
s\right) d|\Lambda |\left( s\right) ,  \label{RSDE3}
\end{equation}%
for some measurable stochastic process $\widehat{\gamma }:\left[ 0,T\right]
\rightarrow
\mathbb{R}
^{n}$ satisfying $\widehat{\gamma }\left( s\right) =\gamma \left( s,X\left(
s\right) \right) $ $d|\Lambda |$-a.s.
\end{definition}

Comparing Definition \ref{skorohodprob} with Definition \ref{strong}, it is
clear that $\left( X\left( \cdot \right) ,\Lambda \left( \cdot \right)
\right) $ should solve the Skorohod problem for $\psi \left( \cdot \right)
=x+\int_{0}^{\cdot }b\left( s,X\left( s\right) \right) ds+\int_{0}^{\cdot
}\left\langle \sigma \left( s,X\left( s\right) \right) ,dW\left( s\right)
\right\rangle $ on an a.s.~pathwise basis. We assume that the coefficient
functions $b\left( t,x\right) $ and $\sigma \left( t,x\right) $ satisfy the
Lipschitz continuity condition
\begin{equation}
\left\vert b_i\left( t,x\right) -b_i\left( t,y\right) \right\vert \leq
K\left\vert x-y\right\vert \quad \text{and\quad }\left\vert \sigma_{i,j}
\left( t,x\right) -\sigma_{i,j} \left( t,y\right) \right\vert \leq
K\left\vert x-y\right\vert ,  \label{lipcoeff}
\end{equation}
for all $(i,j) \in \{1,\dots n\} \times \{1,\dots, m\}$, $x,y \in \mathbb{R}%
^n$ and for some positive constant $K\in \left( 0,\infty \right)$.
Our main result for SDEs is the following theorem.

\begin{theorem}
\label{main}Let $\Omega \subset
\mathbb{R}
^{n+1}$ be a time-dependent domain satisfying \eqref{timesect} and assume
that \eqref{smooth_gamma}, \eqref{boundarylip}, \eqref{templip} and %
\eqref{lipcoeff} hold. Then there exists a unique strong solution to the SDE
in $\overline{\Omega }$ driven by the Wiener process $W$ and with
coefficients $b$ and $\sigma$, direction of reflection along $\gamma $ and
initial condition $x\in \overline{\Omega}_{0}$.
\end{theorem}

We prove Theorem \ref{main} by completing the following steps. First, in
Lemma \ref{smoothexist}, we use a penalty method to prove existence of
solutions to the Skorohod problem for smooth functions. In Lemma \ref%
{compactest}, we then derive a compactness estimate for solutions to the
Skorohod problem. Based on the compactness estimate, we are, in Lemma \ref%
{contexist}, able to generalize the existence result for the Skorohod
problem to all continuous functions. Finally, in Section \ref{RSDE}, we use
two classes of test functions and the existence result for the Skorohod
problem to obtain existence and uniqueness of strong solutions to SDEs with
oblique reflection at the boundary of a bounded, time-dependent domain. Note
that we are able to obtain uniqueness of the reflected SDE although the
solution to the corresponding Skorohod problem need not be unique.

\subsection{Statement of main results for PDEs}

To state and prove our results for PDEs we introduce some more notation. Let
$\Omega ^{\prime }$ be as in \eqref{timedep} and put
\begin{equation*}
\Omega ^{\circ }=\Omega ^{\prime }\cap \left( \left(0,T\right)\times \mathbb{R}%
^{n}\right), \quad \widetilde{\Omega }=\overline{\Omega }^{\prime }\cap
\left( [0,T)\times \mathbb{R}^{n}\right) , \quad \partial \Omega = \left(%
\overline{\Omega}^{\prime} \setminus \Omega ^{\prime }\right) \cap \left(
\left(0,T\right)\times \mathbb{R}^{n}\right).  
\end{equation*}%
%
%
We consider fully nonlinear parabolic PDEs of the form
\begin{equation}
u_{t}+F\left(t,x,u,Du,D^{2}u\right)=0\quad \text{in}\;\Omega ^{\circ }.
\label{huvudekvationen}
\end{equation}%
Here $F$ is a given real function on $\overline{\Omega }\times \mathbb{R}%
\times \mathbb{R}^{n}\times \mathbb{S}^{n}$, where $\mathbb{S}^{n}$ denotes
the space of $n\times n$ real symmetric matrices equipped with the positive
semi-definite ordering; that is, for $X,Y\in \mathbb{S}^{n}$, we write $%
X\leq Y$ if $\langle \left(X-Y\right)\xi ,\xi \rangle \leq 0$ for all $\xi \in \mathbb{R%
}^{n}$. We also adopt the matrix norm notation
\begin{equation*}
\left\Vert A\right\Vert =\sup \{|\lambda |:\lambda \text{ is an eigenvalue
of }A\}=\sup \{|\langle A\xi ,\xi \rangle |:|\xi |\leq 1\}.
\end{equation*}
Moreover, $u$ represents a real function in $\Omega ^{\circ }$ and $Du$ and $%
D^{2}u$ denote the gradient and Hessian matrix, respectively, of $u$ with
respect to the spatial variables. On the boundary we impose the oblique
derivative condition to the unknown $u$
\begin{equation}
\frac{\partial u}{\partial \widetilde{\gamma }}+f\left(t,x,u\left(t,x\right)\right)=0\quad \text{on%
}\;\partial \Omega ,  \label{randvillkor}
\end{equation}%
where $f$ is a real valued function on $\overline{\partial \Omega }\times
\mathbb{R}$ and $\widetilde{\gamma }\left( t,\cdot \right) $ is the vector
field on $\mathbb{R}^{n}$, oblique to $\partial \Omega _{t}$, introduced in %
\eqref{smooth_gamma2} and \eqref{boundarylip2}.

Regarding the function $F$, we make the following assumptions.
\begin{equation}
F\in C\left(\overline{\Omega }\times \mathbb{R}\times \mathbb{R}^{n}\times
\mathbb{S}^{n}\right).  \label{ass_F_cont}
\end{equation}%
For some $\lambda \in \mathbb{R}$ and each $\left(t,x,p,A\right)\in \overline{\Omega }%
\times \mathbb{R}^{n}\times \mathbb{S}^{n}$ the function
\begin{equation}
\text{$r\rightarrow F\left(t,x,r,p,A\right)-\lambda r$ is nondecreasing on $\mathbb{R}$.%
}  \label{ass_F_nondecreasing}
\end{equation}%
There is a function $m_{1}\in C\left([0,\infty)\right)$ satisfying $m_{1}\left(0\right)=0$ for
which
\begin{align}
& F\left(t,y,r,p,-Y\right)-F\left(t,x,r,p,X\right)\leq m_{1}\left(|x-y|\left(|p|+1\right)+\alpha |x-y|^{2}\right)
\label{ass_F_XY} \\
& \text{if}\qquad -\alpha \left(
\begin{array}{cc}
I & 0 \\
0 & I%
\end{array}%
\right) \leq \left(
\begin{array}{cc}
X & 0 \\
0 & Y%
\end{array}%
\right) \leq \alpha \left(
\begin{array}{cc}
I & -I \\
-I & I%
\end{array}%
\right) ,  \notag
\end{align}%
for all $\alpha \geq 1$, $\left(t,x\right),\left(t,y\right)\in \overline{\Omega }$, $r\in \mathbb{R%
}$, $p\in \mathbb{R}^{n}$ and $X,Y\in \mathbb{S}^{n}$, where $I$ denotes the
unit matrix of size $n\times n$. There is a neighborhood $U$ of $\partial
\Omega $ in $\overline{\Omega }$ and a function $m_{2}\in C\left([0,\infty )\right)$
satisfying $m_{2}\left(0\right)=0$ for which
\begin{equation}
|F\left(t,x,r,p,X\right)-F\left(t,x,r,q,Y\right)|\leq m_{2}\left(|p-q|+||X-Y||\right),  \label{ass_F_boundary}
\end{equation}%
for $\left(t,x\right)\in U$, $r\in \mathbb{R}$, $p,q\in \mathbb{R}^{n}$ and $X,Y\in
\mathbb{S}^{n}$. Regarding the function $f$ we assume that
\begin{equation}
f\left(t,x,r\right)\in C\left(\overline{\partial \Omega }\times \mathbb{R}\right),
\label{f_kontinuerlig}
\end{equation}%
and that for each $\left(t,x\right)\in \overline{\partial \Omega }$ the function
\begin{equation}
\text{ $r\rightarrow f\left(t,x,r\right)$ is nondecreasing on $\mathbb{R}$}.
\label{ass_f_nondecreasing}
\end{equation}

We remark that assumptions \eqref{ass_F_cont} and \eqref{ass_F_XY} imply the
degenerate ellipticity
\begin{equation}
F\left(t,x,r,p,A+B\right)\leq F\left(t,x,r,p,A\right)\quad \text{if}\;B\geq 0,
\label{F_fundamental}
\end{equation}%
for $\left(t,x\right)\in \overline{\Omega }$, $r\in \mathbb{R}$, $p\in \mathbb{R}^{n}$
and $A,B\in \mathbb{S}^{n}$, see Remark 3.4 in \cite{CrandallIshiiLions1992}
for a proof. To handle the strong degeneracy allowed, we will adapt the
notion of viscosity solutions \cite{CrandallIshiiLions1992}, which we recall
for problem \eqref{huvudekvationen}-\eqref{randvillkor} in Section \ref{PDE}.
Let $USC(E)$ ($LSC(E)$) denote the set of upper (lower) semi-continuous functions on $E\subset\mathbb{R}^{n+1}$.
Our main results for PDEs are given in the following theorems.

\begin{theorem}
\label{comparison}Let $\Omega ^{\circ }$ be a time-dependent domain
satisfying \eqref{timesect} and assume that \eqref{smooth_gamma2}-%
\eqref{templip} and \eqref{ass_F_cont}-\eqref{ass_f_nondecreasing} hold. Let
$u\in USC(\widetilde{\Omega})$ be a viscosity subsolution, and $v\in LSC(%
\widetilde{\Omega})$ be a viscosity supersolution of problem %
\eqref{huvudekvationen}-\eqref{randvillkor} in $\Omega ^{\circ }$. If $%
u\left(0,x\right)\leq v\left(0,x\right)$ for all $x\in \overline{\Omega }_{0}$, then
$u\leq v\;\text{in}\;\widetilde{\Omega }$. 
\end{theorem}

\begin{theorem}
\label{existence}Let $\Omega ^{\circ }$ be a time-dependent domain
satisfying \eqref{timesect} and assume that \eqref{smooth_gamma2}-%
\eqref{templip} and \eqref{ass_F_cont}-\eqref{ass_f_nondecreasing} hold.
Then there exists a unique viscosity solution, continuous on $\widetilde{%
\Omega }$, to the initial value problem
\begin{align}
u_{t}+F\left(t,x,u,Du,D^{2}u\right)& =0\qquad \quad \,\text{in}\quad \Omega ^{\circ },
\notag  \label{initial_value_problem} \\
\frac{\partial u}{\partial \widetilde{\gamma }}+f\left(t,x,u\left(t,x\right)\right)& =0\qquad
\quad \,\text{on}\quad \partial \Omega ,  \notag \\
u\left(0,x\right)& =g\left(x\right)\qquad \text{for}\quad x\in \overline{\Omega }_{0},
\end{align}%
where $g\in C\left(\overline{\Omega }_{0}\right)$.
\end{theorem}

Theorems \ref{comparison} and \ref{existence} are proved in Section \ref{PDE}%
. The comparison principle in Theorem \ref{comparison} is obtained using two
of the test functions constructed in Section \ref{test} together with
nowadays standard techniques from the theory of viscosity solutions for
fully nonlinear PDEs as described in \cite{CrandallIshiiLions1992}. Our
proof uses ideas from the corresponding elliptic result given in \cite%
{DupuisIshii1990}. The uniqueness part of Theorem \ref{existence} is
immediate from the formulation of Theorem \ref{comparison}, which also,
together with the maximum principle in Lemma \ref{maxrand}, allows comparison
in the setting of mixed boundary conditions, as follows.

\begin{corollary}
\label{maxrand_partial} Let $\Omega ^{\circ }$ be a time-dependent domain
satisfying \eqref{timesect} and assume that \eqref{smooth_gamma2}-%
\eqref{templip} and \eqref{ass_F_cont}-\eqref{ass_f_nondecreasing} hold. Let
$u\in USC(\widetilde{\Omega })$ be a viscosity subsolution, and $v\in LSC(%
\widetilde{\Omega })$ be a viscosity supersolution of \eqref{huvudekvationen}
in $\Omega ^{\circ }$. Suppose also that $u$ and $v$ satisfy the oblique
derivative boundary condition \eqref{randvillkor} on a subset $G\subset
\partial \Omega$.
Then $\sup_{\widetilde{\Omega}}u-v\leq \sup_{\left(\partial \Omega \setminus
G\right)\cup \,\overline{\Omega }_{0}}\left(u-v\right)^{+}$.
\end{corollary}

The existence part of Theorem \ref{existence} is proved using Perron's
method and Corollary \ref{maxrand_partial}, together with constructions of
several explicit viscosity sub- and supersolutions to the problem %
\eqref{huvudekvationen}-\eqref{randvillkor}.

\setcounter{equation}{0} \setcounter{theorem}{0}

\section{Construction of test functions\label{test}}

In this section we show how the classes of test functions constructed in
\cite{DupuisIshii1990} for time-independent domains can be generalized to
similar classes of test functions valid for time-dependent domains. Lemma %
\ref{testlemma3} and Lemma \ref{testlemma4} provide test functions that are
modifications of the square function, but which interact with the direction
of $\gamma $ in a suitable way. The derivations of these functions follow
the lines of the derivations of the corresponding test functions in \cite%
{DupuisIshii1990} with the addition that it has to be verified that the time
derivative of the test functions has a certain order. Lemma \ref{testlemma5}
provides a non-negative test function in $\mathcal{C}^{1,2}\left( \overline{%
\Omega },%
\mathbb{R}
\right) $, whose gradient is aligned with $\gamma $ at the boundary. To
verify the existence of this function, the proof for the corresponding
function in \cite{DupuisIshii1990} has to be extended considerably due to
the time-dependence of the domain. In particular, new methods have to be
used to obtain differentiability with respect to the time variable.

The constructions of the test functions below are given with sufficient
detail and for those parts of the constructions that are identical in
time-dependent and time-independent domains, we refer the reader to \cite%
{DupuisIshii1990}. We start by stating a straightforward extension of Lemma
4.4 in \cite{DupuisIshii1990} from $\xi \in S\left( 0,1\right) $ to $\xi \in
B\left( 0,1\right) $. The proof follows directly from the construction in
Lemma 4.4 in \cite{DupuisIshii1990} and is omitted. For any $\theta \in
\left( 0,1\right) $, there exists a function $g\in \mathcal{C}\left(
\mathbb{R}
^{n}\times
\mathbb{R}
^{n},%
\mathbb{R}
\right) $ and positive constants $\chi ,C$ such that%
\begin{equation}
g\in \mathcal{C}^{1}\left(
\mathbb{R}
^{n}\times
\mathbb{R}
^{n},%
\mathbb{R}
\right) \cap \mathcal{C}^{2}\left(
\mathbb{R}
^{n}\times \left(
\mathbb{R}
^{n}\setminus \left\{ 0\right\} \right) ,%
\mathbb{R}
\right) ,  \label{testlemma21}
\end{equation}%
\begin{equation}
g\left( \xi ,p\right) \geq \chi \left\vert p\right\vert ^{2},\quad \text{for
}\xi \in B\left( 0,1\right) \text{, }p\in
\mathbb{R}
^{n},  \label{testlemma22}
\end{equation}%
\begin{equation}
g\left( \xi ,0\right) =0,\quad \text{for }\xi \in
\mathbb{R}
^{n},  \label{testlemma23}
\end{equation}%
\begin{equation}
\left\langle D_{p}g\left( \xi ,p\right) ,\xi \right\rangle \geq 0,\quad
\text{for }\xi \in S\left( 0,1\right) \text{, }p\in
\mathbb{R}
^{n}\text{ and }\left\langle p,\xi \right\rangle \geq -\theta \left\vert
p\right\vert \text{,}  \label{testlemma24}
\end{equation}%
\begin{equation}
\left\langle D_{p}g\left( \xi ,p\right) ,\xi \right\rangle \leq 0,\quad
\text{for }\xi \in S\left( 0,1\right) \text{, }p\in
\mathbb{R}
^{n}\text{ and }\left\langle p,\xi \right\rangle \leq \theta \left\vert
p\right\vert ,  \label{testlemma25}
\end{equation}%
\begin{equation}
\left\vert D_{\xi }g\left( \xi ,p\right) \right\vert \leq C\left\vert
p\right\vert ^{2},\quad \left\vert D_{p}g\left( \xi ,p\right) \right\vert
\leq C\left\vert p\right\vert ,\quad \text{for }\xi \in B\left( 0,1\right)
\text{, }p\in
\mathbb{R}
^{n},  \label{testlemma26}
\end{equation}%
and
\begin{equation}
\left\Vert D_{\xi }^{2}g\left( \xi ,p\right) \right\Vert \leq C\left\vert
p\right\vert ^{2},\quad \left\Vert D_{\xi }D_{p}g\left( \xi ,p\right)
\right\Vert \leq C\left\vert p\right\vert ,\quad \left\Vert D_{p}^{2}g\left(
\xi ,p\right) \right\Vert \leq C,  \label{testlemma28}
\end{equation}%
for $\xi \in B\left( 0,1\right) $, $p\in
\mathbb{R}
^{n}\setminus \left\{ 0\right\} $. The test function provided by the
following lemma will be used to assert relative compactness of solutions to
the Skorohod problem in Lemma \ref{compactest} below.

\begin{lemma}
\label{testlemma3}For any $\theta \in \left( 0,1\right) $, there exists a
function $h\in \mathcal{C}^{1,2}\left( \left[ 0,T\right] \times
\mathbb{R}
^{n}\times
\mathbb{R}
^{n},%
\mathbb{R}
\right) $ and positive constants $\chi ,C$ such that, for all $\left(
t,x,p\right) \in \lbrack 0,T]\times
\mathbb{R}
^{n}\times
\mathbb{R}
^{n}$,
\begin{equation}
h\left( t,x,p\right) \geq \chi \left\vert p\right\vert ^{2},
\label{testlemma31}
\end{equation}%
\begin{equation}
h\left( t,x,0\right) =1,  \label{testlemma32}
\end{equation}%
\begin{equation}
\left\langle D_{p}h\left( t,x,p\right) ,\gamma \left( t,x\right)
\right\rangle \geq 0,\quad \text{for }x\in \partial \Omega _{t}\text{ and }%
\left\langle p,\gamma \left( t,x\right) \right\rangle \geq -\theta
\left\vert p\right\vert ,  \label{testlemma33}
\end{equation}%
\begin{equation}
\left\langle D_{p}h\left( t,x,p\right) ,\gamma \left( t,x\right)
\right\rangle \leq 0,\quad \text{for }x\in \partial \Omega _{t}\text{ and }%
\left\langle p,\gamma \left( t,x\right) \right\rangle \leq \theta \left\vert
p\right\vert ,  \label{testlemma34}
\end{equation}%
\begin{equation}
\left\vert D_{t}h\left( t,x,p\right) \right\vert \leq C\left\vert
p\right\vert ^{2},\quad \left\vert D_{x}h\left( t,x,p\right) \right\vert
\leq C\left\vert p\right\vert ^{2},\left\vert D_{p}h\left( t,x,p\right)
\right\vert \leq C\left\vert p\right\vert ,  \label{testlemma35}
\end{equation}%
and%
\begin{equation}
\left\Vert D_{x}^{2}h\left( t,x,p\right) \right\Vert \leq C\left\vert
p\right\vert ^{2},\quad \left\Vert D_{x}D_{p}h\left( t,x,p\right)
\right\Vert \leq C\left\vert p\right\vert ,\quad \left\Vert D_{p}^{2}h\left(
t,x,p\right) \right\Vert \leq C.  \label{testlemma37}
\end{equation}
\end{lemma}

\noindent \textbf{Proof.} Let $\nu \in \mathcal{C}^{2}\left(
\mathbb{R}
,%
\mathbb{R}
\right) $ be such that $\nu \left( t\right) =t$ for $t\geq 2$, $\nu \left(
t\right) =1$ for $t\leq 1/2$, $\nu ^{\prime }\left( t\right) \geq 0$ and $%
\nu \left( t\right) \geq t$ for all $t\in
\mathbb{R}
$. Let $\theta \in \left( 0,1\right) $ be given, choose $g\in \mathcal{C}%
\left(
\mathbb{R}
^{n}\times
\mathbb{R}
^{n},%
\mathbb{R}
\right) $ satisfying \eqref{testlemma21}-\eqref{testlemma28} and define%
\begin{equation*}
h\left( t,x,p\right) =\nu \left( g\left( \gamma \left( t,x\right) ,p\right)
\right) .
\end{equation*}%
The regularity of $h$ follows easily from the regularity of $g$ and $\nu $
and \eqref{testlemma23}. It is straightforward to deduce properties %
\eqref{testlemma31}-\eqref{testlemma37} from \eqref{testlemma21}-%
\eqref{testlemma28} and we limit the proof to two examples, which are not
fully covered in \cite{DupuisIshii1990}. We have%
\begin{equation*}
\left\vert D_{t}h\left( t,x,p\right) \right\vert =\left\vert \nu ^{\prime
}\left( g\left( \gamma \left( t,x\right) ,p\right) \right) \right\vert
\left\vert D_{\xi }g\left( \gamma \left( t,x\right) ,p\right) \right\vert
\left\vert \frac{\partial \gamma }{\partial t}\right\vert \leq C\left\vert
p\right\vert ^{2},
\end{equation*}%
by \eqref{testlemma26} and the regularity of $\nu $ and $\gamma $. Moreover,
%
\begin{eqnarray*}
\left\Vert D_{x}^{2}h\left( t,x,p\right) \right\Vert &\leq &C(n)\bigg(%
\left\vert \nu ^{\prime \prime }\left( g\left( \gamma \left( t,x\right)
,p\right) \right) \right\vert \left\vert D_{\xi }g\left( \gamma \left(
t,x\right) ,p\right) \right\vert ^{2}\left\Vert \frac{\partial \gamma }{%
\partial x}\right\Vert ^{2} \\
&&+\left\vert \nu ^{\prime }\left( g\left( \gamma \left( t,x\right)
,p\right) \right) \right\vert \left\Vert D_{\xi }^{2}g\left( \gamma \left(
t,x\right) ,p\right) \right\Vert \left\Vert \frac{\partial \gamma }{\partial
x}\right\Vert ^{2} \\
&&+\left\vert \nu ^{\prime }\left( g\left( \gamma \left( t,x\right)
,p\right) \right) \right\vert \left\vert D_{\xi }g\left( \gamma \left(
t,x\right) ,p\right) \right\vert \max_{1\leq k\leq n}\left\Vert \frac{%
\partial ^{2}\gamma _{k}}{\partial x^{2}}\right\Vert \bigg).
\end{eqnarray*}%
Since $\nu ^{\prime \prime }$ is zero unless $2\geq g\left( \gamma \left(
t,x\right) ,p\right) \geq \chi \left\vert p\right\vert ^{2}$, the first
term, which is of order $C\left\vert p\right\vert ^{4}$, only contributes
for small $\left\vert p\right\vert ^{2}$ and can thus be bounded from above
by $C\left\vert p\right\vert ^{2}$. By \eqref{testlemma26}-%
\eqref{testlemma28}, the two latter terms are also bounded from above by $%
C\left\vert p\right\vert ^{2}$.\hfill $\Box $

\vspace{0.2cm}

The test function in Lemma \ref{testlemma3} is also used to verify the
existence of the following test function, which will be useful in the proofs
of Theorem \ref{comparison} and Lemma \ref{rsdetheorem}.

\begin{lemma}
\label{testlemma4}For any $\theta \in \left( 0,1\right) $, there exists a
family $\left\{ w_{\varepsilon }\right\} _{\varepsilon >0}$ of functions $%
w_{\varepsilon }\in \mathcal{C}^{1,2}\left( \left[ 0,T\right] \times
\mathbb{R}
^{n}\times
\mathbb{R}
^{n},%
\mathbb{R}
\right) $ and positive constants $\chi ,C$ (independent of $\varepsilon $)
such that, for all $\left( t,x,y\right) \in \lbrack 0,T]\times
\mathbb{R}
^{n}\times
\mathbb{R}
^{n}$,%
\begin{equation}
w_{\varepsilon }\left( t,x,y\right) \geq \chi \frac{\left\vert
x-y\right\vert ^{2}}{\varepsilon },  \label{testlemma41}
\end{equation}%
\begin{equation}
w_{\varepsilon }\left( t,x,y\right) \leq C\left( \varepsilon +\frac{%
\left\vert x-y\right\vert ^{2}}{\varepsilon }\right) ,  \label{testlemma42}
\end{equation}%
\begin{equation}
\left\langle D_{x}w_{\varepsilon }\left( t,x,y\right) ,\gamma \left(
t,x\right) \right\rangle \leq C\frac{\left\vert x-y\right\vert ^{2}}{%
\varepsilon },\quad \text{for }x\in \partial \Omega _{t}\text{, }%
\left\langle y-x,\gamma \left( t,x\right) \right\rangle \geq -\theta
\left\vert x-y\right\vert ,  \label{testlemma43}
\end{equation}%
\begin{equation}
\left\langle D_{y}w_{\varepsilon }\left( t,x,y\right) ,\gamma \left(
t,x\right) \right\rangle \leq 0,\quad \text{for }x\in \partial \Omega _{t}%
\text{, }\left\langle x-y,\gamma \left( t,x\right) \right\rangle \geq
-\theta \left\vert x-y\right\vert ,  \label{testlemma49}
\end{equation}%
\begin{equation}
\left\langle D_{y}w_{\varepsilon }\left( t,x,y\right) ,\gamma \left(
t,y\right) \right\rangle \leq C\frac{\left\vert x-y\right\vert ^{2}}{%
\varepsilon },\quad \text{for }y\in \partial \Omega _{t}\text{, }%
\left\langle x-y,\gamma \left( t,y\right) \right\rangle \geq -\theta
\left\vert x-y\right\vert ,  \label{testlemma44}
\end{equation}%
\begin{equation}
\left\vert D_{t}w_{\varepsilon }\left( t,x,y\right) \right\vert \leq C\frac{%
\left\vert x-y\right\vert ^{2}}{\varepsilon },  \label{testlemma45}
\end{equation}%
\begin{equation}
\left\vert D_{y}w_{\varepsilon }\left( t,x,y\right) \right\vert \leq C\frac{%
\left\vert x-y\right\vert }{\varepsilon },\quad \left\vert
D_{x}w_{\varepsilon }\left( t,x,y\right) +D_{y}w_{\varepsilon }\left(
t,x,y\right) \right\vert \leq C\frac{\left\vert x-y\right\vert ^{2}}{%
\varepsilon },  \label{testlemma46}
\end{equation}%
and%
\begin{equation}
D^{2}w_{\varepsilon }\left( t,x,y\right) \leq \frac{C}{\varepsilon }\left(
\begin{array}{cc}
I & -I \\
-I & I%
\end{array}%
\right) +\frac{C\left\vert x-y\right\vert ^{2}}{\varepsilon }\left(
\begin{array}{cc}
I & 0 \\
0 & I%
\end{array}%
\right) .  \label{testlemma47}
\end{equation}
\end{lemma}

\noindent \textbf{Proof.} Let $\theta \in \left( 0,1\right) $ be given and
choose $h\in \mathcal{C}^{1,2}\left( \left[ 0,T\right] \times
\mathbb{R}
^{n}\times
\mathbb{R}
^{n},%
\mathbb{R}
\right) $ as in Lemma \ref{testlemma3}. For all $\varepsilon >0$, we define
the function $w_{\varepsilon }$ as%
\begin{equation*}
w_{\varepsilon }\left( t,x,y\right) =\varepsilon h\left( t,x,\frac{x-y}{%
\varepsilon }\right) .
\end{equation*}%
Property \eqref{testlemma41} follows easily from \eqref{testlemma31} and
property \eqref{testlemma42} was verified in Remark 3.3 in \cite%
{DupuisIshii1993}. Moreover, properties \eqref{testlemma43}, %
\eqref{testlemma49}, \eqref{testlemma46} and \eqref{testlemma47} were
verified in the proof of Theorem 4.1 in \cite{DupuisIshii1990} and %
\eqref{testlemma45} is a simple consequence of \eqref{testlemma35}. To prove %
\eqref{testlemma44}, we note that
\begin{eqnarray*}
\left\langle D_{y}w_{\varepsilon }\left( t,x,y\right) ,\gamma \left(
t,y\right) \right\rangle &=&-\left\langle D_{p}h\left( t,x,\frac{x-y}{%
\varepsilon }\right) ,\gamma \left( t,y\right) \right\rangle \\
&=&-\left\langle D_{p}h\left( t,y,\frac{x-y}{\varepsilon }\right) ,\gamma
\left( t,y\right) \right\rangle \\
&&+\left\langle D_{p}h\left( t,y,\frac{x-y}{\varepsilon }\right)
-D_{p}h\left( t,x,\frac{x-y}{\varepsilon }\right) ,\gamma \left( t,y\right)
\right\rangle .
\end{eqnarray*}%
Moreover, if $\left\langle x-y,\gamma \left( t,y\right) \right\rangle \geq
-\theta \left\vert x-y\right\vert $, then by \eqref{testlemma33}, $%
\left\langle D_{p}h\left( t,y,p\right) ,\gamma \left( t,y\right)
\right\rangle \geq 0$ with $p=\left( x-y\right) /\varepsilon $. Hence, for
some $\xi $ in the segment joining $x$ and $y$, we obtain, with the aid of
the mean value theorem and \eqref{testlemma37},
\begin{eqnarray*}
\left\langle D_{y}w_{\varepsilon }\left( t,x,y\right) ,\gamma \left(
t,y\right) \right\rangle &\leq &\left\Vert D_{x}D_{p}h\left( t,\xi ,\frac{x-y%
}{\varepsilon }\right) \right\Vert \left\vert x-y\right\vert \\
&\leq &C\left\vert \frac{x-y}{\varepsilon }\right\vert \left\vert
x-y\right\vert =C\frac{\left\vert x-y\right\vert ^{2}}{\varepsilon }.
\end{eqnarray*}%
\hfill \hfill $\Box $

We conclude this section by proving Lemma \ref{testlemma5} using an
appropriate Cauchy problem. The test function $\alpha $ in Lemma \ref%
{testlemma5} will be crucial for the proofs of Theorems \ref{comparison}, %
\ref{existence} and Lemma \ref{rsdetheorem}.

\begin{lemma}
\label{testlemma5} There exists a nonnegative function $\alpha \in \mathcal{C%
}^{1,2}\left( \overline{\Omega },%
\mathbb{R}
\right) $, which satisfies
\begin{equation}
\left\langle D_{x}\alpha \left( t,x\right) ,\gamma \left( t,x\right)
\right\rangle \geq 1,  \label{alfaprop}
\end{equation}%
for $x\in \partial \Omega _{t}$, $t\in \left[ 0,T\right] $. Moreover, the
support of $\alpha$ can be assumed to lie in the neighbourhood $U$ defined
in \eqref{ass_F_boundary}.
\end{lemma}

\noindent \textbf{Proof.} Fix $s\in \left[ 0,T\right] $ and $z\in \partial
\Omega _{s}$ and define $H_{s,z}$ as the hyperplane%
\begin{equation*}
H_{s,z}=\left\{ x\in
\mathbb{R}
^{n}:\left\langle x-z,\gamma \left( s,z\right) \right\rangle =0\right\} .
\end{equation*}%
Given a function $u_{0}\in \mathcal{C}^{2}\left( H_{s,z},%
\mathbb{R}
\right) $, such that $u_{0}\left( z\right) =1$, $u_{0}\geq 0$ and supp $%
u_{0}\subset B\left( z,\delta ^{2}/4\right) \cap H_{s,z}$, we can use the
method of characteristics to solve the Cauchy problem%
\begin{eqnarray*}
\left\langle D_{x}u_{\left( t\right) }\left( x\right) ,\gamma \left(
t,x\right) \right\rangle &=&0, \\
\left. u_{\left( t\right) }\right\vert _{H_{s,z}} &=&u_{0}.
\end{eqnarray*}%
Choosing the positive constants $\delta $ and $\eta $ sufficiently small,
the Cauchy problem above has, for all $t\in \left[ s-\eta ,s+\eta \right] $,
a solution $u_{\left( t\right) }\in \mathcal{C}^{2}\left( B\left( z,\delta
\right) ,%
\mathbb{R}
\right) $ satisfying $u_{\left( t\right) }\geq 0$. Based on the continuity
of $\gamma $ and the restriction on the support of $u_{0}$, we may also
assume that%
\begin{equation*}
\text{supp }u_{\left( t\right) }\subset \bigcup_{\zeta \in
\mathbb{R}
}B(z-\zeta \gamma (s,z),\delta ^{2}/3)\cap B\left( z,\delta \right) .
\end{equation*}%
Next, we define the combined function
\begin{equation*}
u\left( t,x\right) =u_{\left( t\right) }\left( x\right) ,
\end{equation*}%
and we claim for now that $u\in \mathcal{C}^{1,2}\left( \left[ s-\eta
,s+\eta \right] \times B\left( z,\delta \right) ,%
\mathbb{R}
\right) $ and postpone the proof of this claim to the end of the proof of
the lemma. By the exterior and interior cone conditions, we can, for
sufficiently small $\delta $, find $\varepsilon >0$ such that%
\begin{eqnarray*}
&&\bigcup_{\zeta >0}B(z-\zeta \gamma (s,z),\delta ^{2}/3)\cap \left(B\left(
z,\delta \right) \setminus B\left( z,\delta -2\varepsilon \right)\right) \\
&\subset &\bigcup_{\zeta >0}B(z-\zeta \gamma (s,z),\zeta \delta )\cap
B\left( z,\delta \right) \subset \Omega _{s}^{c},
\end{eqnarray*}%
and such that the similar union over $\zeta <0$ belongs to $\Omega _{s}$.
Hence%
\begin{equation*}
\partial \Omega _{s}\cap \left(\text{supp }u_{\left( t\right) }\setminus
B\left( z,\delta -2\varepsilon \right)\right) =\emptyset ,
\end{equation*}%
and, by \eqref{tempholder}, it follows that if $\eta $ also satisfies the
constraint $\eta <\left( \varepsilon /K\right) ^{1/\widehat{\alpha }}$, then%
\begin{equation}
\partial \Omega _{t}\cap \left(\text{supp }u_{\left( t\right) }\setminus
B\left( z,\delta -\varepsilon \right)\right) =\emptyset \text{,}
\label{suppbound}
\end{equation}%
for all $t\in \left[ s-\eta ,s+\eta \right] $.

Now, choose a function $\xi \in \mathcal{C}_{0}^{1,2}\left( \left[ s-\eta
,s+\eta \right] \times B\left( z,\delta \right) ,%
\mathbb{R}
\right) $ so that $\xi \left( t,x\right) =1$ for $t\in \left[ s-\eta +
\varepsilon ,s+\eta - \varepsilon \right] $, $x\in B\left( z,\delta
-\varepsilon \right) $ and $\xi \geq 0$, and set%
\begin{equation*}
v_{s,z}\left( t,x\right) =u\left( t,x\right) \xi \left( t,x\right) .
\end{equation*}%
Then $v_{s,z}\in \mathcal{C}_{0}^{1,2}\left( \left[ s-\eta ,s+\eta \right]
\times B\left( z,\delta \right) ,%
\mathbb{R}
\right) $ satisfies $v_{s,z}\geq 0$. By \eqref{suppbound} and the
construction of $u$ and $\xi $, we obtain%
\begin{equation*}
\left\langle D_{x}v_{s,z}\left( t,x\right) ,\gamma \left( t,x\right)
\right\rangle =0\text{ for }x\in B\left( z,\delta \right) \cap \partial
\Omega _{t}\text{, }t\in \left[ s-\eta ,s+\eta \right] .
\end{equation*}%
Define $w_{s,z}\in \mathcal{C}^{2}\left( B\left( z,\delta \right) ,%
\mathbb{R}
\right) $ by%
\begin{equation*}
w_{s,z}\left( x\right) =\left\langle x-z,\gamma \left( s,z\right)
\right\rangle +M,
\end{equation*}%
where $M$ is large enough so that $w_{s,z}\geq 0$. Using the continuity of $%
\gamma $, we can find $\delta $ and $\eta $ such that $\left\langle \gamma
\left( s,z\right) ,\gamma \left( t,x\right) \right\rangle \geq 0$ for all $%
\left( t,x\right) \in \left[ s-\eta ,s+\eta \right] \times B\left( z,\delta
\right) $. Setting%
\begin{equation*}
g_{s,z}\left( t,x\right) =v_{s,z}\left( t,x\right) w_{s,z}\left( x\right) ,
\end{equation*}%
we find that $g_{s,z}\in \mathcal{C}_{0}^{1,2}\left( \left[ s-\eta ,s+\eta %
\right] \times B\left( z,\delta \right) ,%
\mathbb{R}
\right) $ satisfies $g_{s,z}\geq 0$. Moreover, using $\left\vert \gamma
\left( t,x\right) \right\vert =1$, we have%
\begin{eqnarray*}
\left\langle D_{x}g_{s,z}\left( s,z\right) ,\gamma \left( s,z\right)
\right\rangle &=&v_{s,z}\left( s,z\right) \left\langle D_{x}w_{s,z}\left(
z\right) ,\gamma \left( s,z\right) \right\rangle \\
&&+w_{s,z}\left( z\right) \left\langle D_{x}v_{s,z}\left( s,z\right) ,\gamma
\left( s,z\right) \right\rangle \\
&=&u\left( s,z\right) \xi \left( s,z\right) \left\vert \gamma \left(
s,z\right) \right\vert ^{2}=1,
\end{eqnarray*}%
and a similar calculation shows that%
\begin{eqnarray*}
\left\langle D_{x}g_{s,z}\left( t,x\right) ,\gamma \left( t,x\right)
\right\rangle &=&v_{s,z}\left( t,x\right) \left\langle D_{x}w_{s,z}\left(
x\right) ,\gamma \left( t,x\right) \right\rangle \\
&&+w_{s,z}\left( x\right) \left\langle D_{x}v_{s,z}\left( t,x\right) ,\gamma
\left( t,x\right) \right\rangle \\
&=&v_{s,z}\left( t,x\right) \left\langle \gamma \left( s,z\right) ,\gamma
\left( t,x\right) \right\rangle \geq 0,
\end{eqnarray*}%
for $x\in B\left( z,\delta \right) \cap \partial \Omega _{t}$, $t\in \left[
s-\eta ,s+\eta \right] $. Now, 
using a standard compactness argument we conclude the existence of a
nonnegative function $\alpha \in \mathcal{C}^{1,2}(\overline{\Omega },%
\mathbb{R})$, which satisfies $\left\langle D_{x}\alpha \left( t,x\right)
,\gamma \left( t,x\right) \right\rangle \geq 1$ for $x\in \partial \Omega
_{t}$, $t\in \left[ 0,T\right] $. Moreover, by the above construction, we
can assume that the support of $\alpha $ lies within the neighbourhood $U$
defined in \eqref{ass_F_boundary}.

It remains to prove the proposed regularity $u\in \mathcal{C}^{1,2}\left( %
\left[ s-\eta ,s+\eta \right] \times B\left( z,\delta \right) ,%
\mathbb{R}
\right) $. The regularity in the spatial variables follows directly by
construction, so it remains to show that $u$ is continuously differentiable
in the time variable. Let $x\in B\left( z,\delta \right) $ and let $t$ and $%
t+h$ belong to $\left[ s-\eta ,s+\eta \right] $. Denote by $y\left( t,\cdot
\right) $ and $y\left( t+h,\cdot \right) $ the characteristic curves through
$x$ for the vector fields $\gamma \left( t,\cdot \right) $ and $\gamma
\left( t+h,\cdot \right) $, respectively, so that
\begin{eqnarray*}
\frac{\partial y}{\partial r}\left( t,r\right) &=&\pm \gamma \left(
t,y\left( t,r\right) \right) , \\
y\left( t,0\right) &=&x,
\end{eqnarray*}%
and analogously for $y\left( t+h,\cdot \right) $. Choose the sign in the
parametrization of $y\left( t,\cdot \right) $ so that there exists some $%
r\left( t\right) >0$ such that $y\left( t,r\left( t\right) \right) =z\left(
t\right) \in H_{s,z}$. Choosing the same sign in the parametrization of $%
y\left( t+h,\cdot \right) $ asserts the existence of some $r\left(
t+h\right) >0$ such that $y\left( t+h,r\left( t+h\right) \right) =z\left(
t+h\right) \in H_{s,z}$. Without lack of generality, we assume the sign
above to be positive. Since $u\left( t,x\right) =u_{0}\left( z\left(
t\right) \right) $, where $u_{0}$ is continuously differentiable, it remains
to show that the function $z$ is continuously differentiable.

We will first show that $y\left( \cdot ,r\right) $ is continuously
differentiable by following an argument that can be found in e.g. \cite%
{Monti2010}. Differentiating the Cauchy problem formally with respect to the
time variable and introducing the function $\psi \left( t,r\right) =\dfrac{%
\partial y}{\partial t}\left( t,r\right) $, we obtain
\begin{eqnarray*}
\frac{\partial }{\partial r}\psi \left( t,r\right) &=&\frac{\partial \gamma
}{\partial t}\left( t,y\left( t,r\right) \right) +\frac{\partial \gamma }{%
\partial y}\left( t,y\left( t,r\right) \right) \psi \left( t,r\right) , \\
\psi \left( t,0\right) &=&0.
\end{eqnarray*}%
This Cauchy problem has a unique solution, which we will next show satisfies
\begin{equation}
\psi \left( t,r\right) =\lim_{h\rightarrow 0}\frac{y\left( t+h,r\right)
-y\left( t,r\right) }{h},  \label{eq:limit_true}
\end{equation}%
so that $\psi $ is in fact the time derivative of $y$ (not just formally).
Define
\begin{equation}
R\left( t,r,h\right) =\frac{y\left( t+h,r\right) -y\left( t,r\right) }{h}%
-\psi \left( t,r\right) .  \label{eq:R_def}
\end{equation}%
Now
\begin{eqnarray*}
R\left( t,r,h\right) &=&\int_{0}^{r}\left( \frac{\gamma \left( t+h,y\left(
t+h,u\right) \right) -\gamma \left( t,y\left( t,u\right) \right) }{h}\right)
du \\
&&-\int_{0}^{r}\left( \frac{\partial \gamma }{\partial t}\left( t,y\left(
t,u\right) \right) +\frac{\partial \gamma }{\partial y}\left( t,y\left(
t,u\right) \right) \psi \left( t,u\right) \right) du.
\end{eqnarray*}%
By the mean value theorem
\begin{eqnarray*}
&&\gamma _{i}\left( t+h,y\left( t+h,u\right) \right) -\gamma _{i}\left(
t,y\left( t,u\right) \right) \\
&=&\frac{\partial \gamma _{i}}{\partial t}\left( \overline{t}_{i},y\left(
t+h,u\right) \right) h+\frac{\partial \gamma _{i}}{\partial y}\left( t,%
\overline{y}_{i}\right) \left( y\left( t+h,u\right) -y\left( t,u\right)
\right) ,
\end{eqnarray*}%
for some $\overline{t}_{i}$ between $t$ and $t+h$, some $\overline{y}_{i}$
between $y\left( t,u\right) $ and $y\left( t+h,u\right) $ and all $i\in
\left\{ 1,...,n\right\} $. Hence, the $i^{\text{th}}$ component of $R\left(
t,r,h\right) $ is
\begin{eqnarray*}
R_{i}\left( t,r,h\right) &=&\int_{0}^{r}\left( \frac{\partial \gamma _{i}}{%
\partial t}\left( \overline{t}_{i},y\left( t+h,u\right) \right) -\frac{%
\partial \gamma _{i}}{\partial t}\left( t,y\left( t,u\right) \right) \right)
du \\
&&+\int_{0}^{r}\left( \frac{\partial \gamma _{i}}{\partial y}\left( t,%
\overline{y}_{i}\right) \frac{y\left( t+h,u\right) -y\left( t,u\right) }{h}-%
\frac{\partial \gamma _{i}}{\partial y}\left( t,y\left( t,u\right) \right)
\psi \left( t,u\right) \right) du,
\end{eqnarray*}%
where the second term on the right hand side can be rewritten as
\begin{equation*}
\int_{0}^{r}\left( \frac{\partial \gamma _{i}}{\partial y}\left( t,\overline{%
y}_{i}\right) R\left( t,u,h\right) +\left( \frac{\partial \gamma _{i}}{%
\partial y}\left( t,\overline{y}_{i}\right) -\frac{\partial \gamma _{i}}{%
\partial y}\left( t,y\left( t,u\right) \right) \right) \psi \left(
t,u\right) \right) du.
\end{equation*}%
Therefore we have
\begin{eqnarray*}
\left\vert R\left( t,r,h\right) \right\vert &\leq
&\int_{0}^{r}\sum_{i=1}^{n}\left\vert \frac{\partial \gamma _{i}}{\partial t}%
\left( \overline{t}_{i},y\left( t+h,u\right) \right) -\frac{\partial \gamma
_{i}}{\partial t}\left( t,y\left( t,u\right) \right) \right\vert du \\
&&+\int_{0}^{r}\left\vert R\left( t,u,h\right) \right\vert
\sum_{i=1}^{n}\left\vert \frac{\partial \gamma _{i}}{\partial y}\left( t,%
\overline{y}_{i}\right) \right\vert du \\
&&+\int_{0}^{r}\left\vert \psi \left( t,u\right) \right\vert
\sum_{i=1}^{n}\left\vert \frac{\partial \gamma _{i}}{\partial y}\left( t,%
\overline{y}_{i}\right) -\frac{\partial \gamma _{i}}{\partial y}\left(
t,y\left( t,u\right) \right) \right\vert du,
\end{eqnarray*}%
and by Gronwall's inequality we obtain
\begin{eqnarray}
\left\vert R\left( t,r,h\right) \right\vert &\leq
&C\int_{0}^{r}\sum_{i=1}^{n}\left\vert \frac{\partial \gamma _{i}}{\partial t%
}\left( \overline{t}_{i},y\left( t+h,u\right) \right) -\frac{\partial \gamma
_{i}}{\partial t}\left( t,y\left( t,u\right) \right) \right\vert du  \notag
\label{eq:gronvall_deriv} \\
&&+C\int_{0}^{r}\left\vert \psi \left( t,u\right) \right\vert
\sum_{i=1}^{n}\left\vert \frac{\partial \gamma _{i}}{\partial y}\left( t,%
\overline{y}_{i}\right) -\frac{\partial \gamma _{i}}{\partial y}\left(
t,y\left( t,u\right) \right) \right\vert du,
\end{eqnarray}%
for some positive constant $C$. Since $\left\vert \psi \left( t,u\right)
\right\vert $ exists and is bounded, and since the time and space
derivatives of $\gamma $ are continuous, \eqref{eq:gronvall_deriv} implies
boundedness of $\left\vert R\left( t,r,h\right) \right\vert $. Therefore, by %
\eqref{eq:R_def} we have $|y\left( t+h,r\right) -y\left( t,r\right) |\leq Ch$%
, for some constant $C$, and we can conclude that $\overline{y}%
_{i}\rightarrow y(t,u)$ and $\overline{t}_{i}\rightarrow t$ for all $i\in
\{1,\dots ,n\}$ as $h\rightarrow 0$. It follows that the differences in the
integrands in \eqref{eq:gronvall_deriv} vanish as $h\rightarrow 0$ and hence
$\lim_{h\rightarrow 0}R\left( t,r,h\right) =0$. This proves %
\eqref{eq:limit_true} and therefore that $y\left( \cdot ,r\right) $ is
continuously differentiable.

Now, by the mean value theorem,
\begin{eqnarray*}
z_{i}(t+h)-z_{i}(t) &=&y_{i}(t+h,r\left( t+h\right) )-y_{i}(t,r\left(
t\right) ) \\
&=&y_{i}(t+h,r\left( t+h\right) )-y_{i}(t+h,r\left( t\right) ) \\
&&+y_{i}(t+h,r\left( t\right) )-y_{i}(t,r\left( t\right) ) \\
&=&\frac{\partial y_{i}}{\partial r}\left( t+h,\overline{r}_{i}\right)
\left( r\left( t+h\right) -r\left( t\right) \right) +\frac{\partial y_{i}}{%
\partial t}\left( \overline{t}_{i},r\left( t\right) \right) h,
\end{eqnarray*}
for some $\overline{r}_{i}$ between $r\left( t\right) $ and $%
r\left(t+h\right)$, some $\overline{t}_{i}$ between $t$ and $t+h$ and all $%
i\in \left\{ 1,...,n\right\}$.
Since the function $r(t)$ is defined so that
\begin{align*}
\langle y(t, r(t)) - z, \gamma(s, z)\rangle = 0, \quad t \in (s - \eta, s +
\eta),
\end{align*}
it follows by the implicit function theorem and by the regularity of $y(t,r)$
that $r(t)$ is a continuously differentiable function. Hence, we conclude
that $\overline{r}_i \rightarrow r(t)$ and $\overline{t}_i \rightarrow t$,
all $i \in \{1,\dots,n\}$, as $h \rightarrow 0$ and therefore,
\begin{equation*}
\lim_{h\rightarrow 0}\frac{z(t+h)-z(t)}{h}=\frac{\partial y}{\partial r}%
\left( t,r\left( t\right) \right) r^{\prime }\left( t\right) +\frac{\partial
y}{\partial t}\left( t,r\left( t\right) \right) ,
\end{equation*}
where the right hand side is a continuous function. This proves that $z$ is
continuously differentiable and, hence, that $u\in \mathcal{C}^{1,2}\left( %
\left[ s-\eta ,s+\eta \right] \times B\left( z,\delta \right) ,%
\mathbb{R}
\right) $.\hfill $\Box $

\setcounter{equation}{0} \setcounter{theorem}{0}

\section{The Skorohod problem\label{SP}}

In this section we prove existence of solutions to the Skorohod problem
under the assumptions in Section \ref{geoassume}. This result could be
achieved using the methods in \cite{NystromOnskog2010a}, but as we here
assume more regularity on the direction of reflection and the temporal
variation of the domain compared to the setting in \cite{NystromOnskog2010a}
(and this is essential for the other sections of this article), we follow a
more direct approach using a penalty method. We first note that, mimicking
the proof of Lemma 4.1 in \cite{DupuisIshii1993}, we can prove the following
result.

\begin{lemma}
\label{dlemma} There is a constant $\mu >0$ such that, for every $t\in \left[
0,T\right] $, there exists a neighbourhood $U_{t}$ of $\partial \Omega _{t}$
such that%
\begin{equation}
\left\langle D_{x}d\left( t,x\right) ,\gamma \left( t,x\right) \right\rangle
\leq -\mu ,\quad \text{for a.e. }x\in U_{t}\setminus \overline{\Omega }_{t}%
\text{.}  \label{unmollineq}
\end{equation}
\end{lemma}

As \eqref{unmollineq} holds only for almost every point in a neighbourhood
of a non-smooth domain, we cannot apply \eqref{unmollineq} directly and will
use the following mollifier approach instead. Based on the construction of
the neighbourhoods $\left\{ U_{t}\right\} _{t\in \left[ 0,T\right] }$ in
Lemma \ref{dlemma} (see the proof of the corresponding lemma in \cite%
{DupuisIshii1993} for details), there exists a constant $\overline{\beta }>0$
such that $B\left( x,3\overline{\beta }\right) \subset U_{t}$ for all $x\in
\partial \Omega _{t}$, $t\in \left[ 0,T\right] $. For the value of $p$ given
in \eqref{templip}, let
\begin{equation*}
v\left( t,x\right) =\left( d\left( t,x\right) \right) ^{p}\quad \text{%
and\quad }\widetilde{v}\left( t,x\right) =\left( d\left( t,x\right) \right)
^{p-1}.
\end{equation*}%
Moreover, let $\varphi _{\beta }\in $ $\mathcal{C}^{\infty }\left(
\mathbb{R}
^{n},%
\mathbb{R}
\right) $ be a positive mollifier with support in $B\left( 0,\beta \right) $%
, for some $\beta >0$, and define the spatial convolutions
\begin{equation*}
v_{\beta }=v\ast \varphi _{\beta }\quad \text{and\quad }\widetilde{v}_{\beta
}=\widetilde{v}\ast \varphi _{\beta }.
\end{equation*}

\begin{lemma}
\label{vlemma} There is a constant $\kappa >0$ such that, for sufficiently
small $\beta > 0$ and every $t\in \left[ 0,T\right]$, there exists a
neighbourhood $\widetilde{U}_{t}$ of $\partial \Omega _{t}$, $\widetilde{U}%
_{t}\supset \left\{ x:d\left( x,\partial \Omega _{t}\right) <2\overline{%
\beta }\right\}$, such that
\begin{equation}
\left\langle D_{x}v_{\beta }\left( t,x\right) ,\gamma \left( t,x\right)
\right\rangle \leq -\kappa \widetilde{v}_{\beta }\left( t,x\right) ,\quad
\text{for }x\in \widetilde{U}_{t}\setminus \overline{\Omega }_{t}\text{.}
\label{mollineq}
\end{equation}
\end{lemma}

\noindent \textbf{Proof.} For all $x\in U_{t}\setminus \overline{\Omega }%
_{t} $ such that $B\left( x,\overline{\beta }\right) \subset U_{t}$ and for
all $\beta \leq \overline{\beta }$, we have%
\begin{eqnarray*}
&&\left\langle D_{x}v_{\beta }\left( t,x\right) ,\gamma \left( t,x\right)
\right\rangle =\int_{%
\mathbb{R}
^{n}}\left\langle \varphi _{\beta }\left( x-y\right) D_{y}v\left( t,y\right)
,\gamma \left( t,x\right) \right\rangle dy \\
&=&\int_{%
\mathbb{R}
^{n}}\left( \left\langle D_{y}v\left( t,y\right) ,\gamma \left( t,y\right)
\right\rangle +\left\langle D_{y}v\left( t,y\right) ,\gamma \left(
t,x\right) -\gamma \left( t,y\right) \right\rangle \right) \varphi _{\beta
}\left( x-y\right) dy.
\end{eqnarray*}%
The inner product in the second term is bounded from above by
\begin{equation*}
p\left( d\left( t,y\right) \right) ^{p-1}\left\vert D_{y}d\left( t,y\right)
\right\vert L\beta ,
\end{equation*}
where $L$ is the Lipschitz coefficient of $\gamma $ in spatial dimensions
over the compact set $\left[ 0,T\right] \times \bigcup_{t\in \left[ 0,T%
\right] }\overline{U}_{t}$. By Lemma \ref{dlemma}, we have, for almost every
$y\in U_{t}\setminus \overline{\Omega }_{t}$, $t\in \left[ 0,T\right] $,%
\begin{equation*}
\left\langle D_{y}v\left( t,y\right) ,\gamma \left( t,y\right) \right\rangle
=p\left( d\left( t,y\right) \right) ^{p-1}\left\langle D_{y}d\left(
t,y\right) ,\gamma \left( t,y\right) \right\rangle \leq -p\mu \left( d\left(
t,y\right) \right) ^{p-1},
\end{equation*}%
and, for sufficiently small $\beta >0$,
\begin{equation*}
p\left( d\left( t,y\right) \right) ^{p-1}L\beta -p\mu \left( d\left(
t,y\right) \right) ^{p-1}\leq -\kappa \left( d\left( t,y\right) \right)
^{p-1},
\end{equation*}%
for some constant $\kappa >0$. This proves \eqref{mollineq}.\hfill $\Box $

\vspace{0.2cm}

We next use a penalty method to verify the existence of a solution to the
Skorohod problem for continuously differentiable functions. The following
lemma generalizes Theorem 2.1 in \cite{LionsSznitman1984} and Lemma 4.5 in
\cite{DupuisIshii1993}.

\begin{lemma}
\label{smoothexist}Let $\psi \in \mathcal{C}^{1}\left( \left[ 0,T\right] ,%
\mathbb{R}
^{n}\right) $ with $\psi \left( 0\right) \in \overline{\Omega }_{0}$. Then
there exists a solution $\left( \phi ,\lambda \right) \in \mathcal{W}%
^{1,p}\left( \left[ 0,T\right] ,%
\mathbb{R}
^{n}\right) \times \mathcal{W}^{1,p}\left( \left[ 0,T\right] ,%
\mathbb{R}
^{n}\right) $ to the Skorohod problem for $\left( \Omega ,\gamma ,\psi
\right) $.
\end{lemma}

\noindent \textbf{Proof.} Choose $\varepsilon >0$ and consider the ordinary
differential equation%
\begin{equation}
\phi _{\varepsilon }^{\prime }\left( t\right) =\frac{1}{\varepsilon }d\left(
t,\phi _{\varepsilon }\left( t\right) \right) \gamma \left( t,\phi
_{\varepsilon }\left( t\right) \right) +\psi ^{\prime }\left( t\right)
,\quad \phi _{\varepsilon }\left( 0\right) =\psi \left( 0\right) ,
\label{ode}
\end{equation}%
for $\phi _{\varepsilon }\left( t\right) $, which has a unique solution on $%
\left[ 0,T\right] $. Let $\kappa >0$ and the family of neighbourhoods $\{%
\widetilde{U}_{t}\}_{t\in \left[ 0,T\right] }$ be as in Lemma \ref{vlemma}.
Choose a function $\zeta \in \mathcal{C}^{\infty }\left( \left[ 0,\infty
\right) ,\left[ 0,\infty \right) \right) $ such that%
\begin{equation*}
\zeta \left( r\right) =\left\{
\begin{array}{ll}
r, & \text{for }r\leq \overline{\beta }^{p}/2, \\
3\overline{\beta }^{p}/4, & \text{for }r\geq \overline{\beta }^{p},%
\end{array}%
\right.
\end{equation*}%
and $0\leq \zeta ^{\prime }\left( r\right) \leq 1$ for all $r\in \left[
0,\infty \right) $. Note that if $\phi _{\varepsilon }\left( t\right) \notin
\widetilde{U}_{t}\cup \overline{\Omega }_{t}$, then $d\left( t,\phi
_{\varepsilon }\left( t\right) \right) \geq 2\overline{\beta }$ and, as a
consequence, for all $\beta \leq \overline{\beta }$ it holds that $v_{\beta
}\left( t,\phi _{\varepsilon }\left( t\right) \right) \geq \overline{\beta }%
^{p}$ and $\zeta ^{\prime }\left( v_{\beta }\left( t,\phi _{\varepsilon
}\left( t\right) \right) \right) =0$. We next define the function $F\left(
t\right) =\zeta \left( v_{\beta }\left( t,\phi _{\varepsilon }\left(
t\right) \right) \right) $, for $t\in \left[ 0,T\right] $, and investigate
its time derivative. Let $D_{t}d$ denote the weak derivative guaranteed by %
\eqref{templip} and note that
\begin{eqnarray}
F^{\prime }\left( t\right) &=&\zeta ^{\prime }\left( v_{\beta }\left( t,\phi
_{\varepsilon }\left( t\right) \right) \right) \left( D_{t}v_{\beta }\left(
t,\phi _{\varepsilon }\left( t\right) \right) +\left\langle D_{x}v_{\beta
}\left( t,\phi _{\varepsilon }\left( t\right) \right) ,\phi _{\varepsilon
}^{\prime }\left( t\right) \right\rangle \right)  \notag \\
&=&\zeta ^{\prime }\left( v_{\beta }\left( t,\phi _{\varepsilon }\left(
t\right) \right) \right) \Bigg(D_{t}v_{\beta }\left( t,\phi _{\varepsilon
}\left( t\right) \right) \frac{{}}{{}}  \notag \\
&&+\left\langle D_{x}v_{\beta }\left( t,\phi _{\varepsilon }\left( t\right)
\right) ,\,\frac{1}{\varepsilon }d\left( t,\phi _{\varepsilon }\left(
t\right) \right) \gamma \left( t,\phi _{\varepsilon }\left( t\right) \right)
+\psi ^{\prime }\left( t\right) \right\rangle \Bigg),  \label{vprim}
\end{eqnarray}%
as $\phi _{\varepsilon }\left( t\right) $ solves \eqref{ode}. From Lemma \ref%
{vlemma}, we have
\begin{eqnarray*}
&&\zeta ^{\prime }\left( v_{\beta }\left( t,\phi _{\varepsilon }\left(
t\right) \right) \right) \left\langle D_{x}v_{\beta }\left( t,\phi
_{\varepsilon }\left( t\right) \right) ,\,\frac{1}{\varepsilon }d\left(
t,\phi _{\varepsilon }\left( t\right) \right) \gamma \left( t,\phi
_{\varepsilon }\left( t\right) \right) \right\rangle \\
&\leq &-\zeta ^{\prime }\left( v_{\beta }\left( t,\phi _{\varepsilon }\left(
t\right) \right) \right) \frac{\kappa }{\varepsilon }d\left( t,\phi
_{\varepsilon }\left( t\right) \right) \widetilde{v}_{\beta }\left( t,\phi
_{\varepsilon }\left( t\right) \right) ,
\end{eqnarray*}%
for $\phi _{\varepsilon }\left( t\right) \in \widetilde{U}_{t}\setminus
\overline{\Omega }_{t}$ and for all other $\phi _{\varepsilon }\left(
t\right) $ both sides vanish when $\beta \leq \overline{\beta }$.
Integrating the estimate for $F^{\prime }$, suppressing the $s$-dependence
in $\phi _{\varepsilon }$ and $\psi $ and denoting $\zeta ^{\prime }\left(
v_{\beta }\left( s,\phi _{\varepsilon }\right) \right) $ by $\zeta ^{\prime
}\left( v_{\beta }\right) $ for simplicity, we obtain, for all $t\in \left[
0,T\right] $,
\begin{eqnarray}
&&\zeta \left( v_{\beta }\left( t,\phi _{\varepsilon }\left( t\right)
\right) \right) -\zeta \left( v_{\beta }\left( 0,\phi _{\varepsilon }\left(
0\right) \right) \right) +\frac{\kappa }{\varepsilon }\int_{0}^{t}\zeta
^{\prime }\left( v_{\beta }\right) d\left( s,\phi _{\varepsilon }\right)
\widetilde{v}_{\beta }\left( s,\phi _{\varepsilon }\right) ds
\label{eq:nyref} \\
&\leq &\int_{0}^{t}\zeta ^{\prime }\left( v_{\beta }\right) \left\vert
D_{s}v_{\beta }\left( s,\phi _{\varepsilon }\right) \right\vert
ds+\int_{0}^{t}\zeta ^{\prime }\left( v_{\beta }\right) \left\vert
D_{x}v_{\beta }\left( s,\phi _{\varepsilon }\right) \right\vert \left\vert
\psi ^{\prime }\right\vert ds=I_{1}+I_{2}.  \notag
\end{eqnarray}%
Note that since $\left\vert D_{x}d\right\vert \leq 1$ a.e. we have $%
\left\vert D_{x}v_{\beta }\left( s,\phi _{\varepsilon }\right) \right\vert
\leq p\widetilde{v}_{\beta }\left( s,\phi _{\varepsilon }\right) $, and
hence, H\"{o}lder's inequality implies
\begin{eqnarray*}
I_{2} &=&\int_{0}^{t}\zeta ^{\prime }\left( v_{\beta }\right) \left\vert
D_{x}v_{\beta }\left( s,\phi _{\varepsilon }\right) \right\vert \left\vert
\psi ^{\prime }\right\vert ds\leq p\int_{0}^{t}\zeta ^{\prime }\left(
v_{\beta }\right) \widetilde{v}_{\beta }\left( s,\phi _{\varepsilon }\right)
\left\vert \psi ^{\prime }\right\vert ds \\
&\leq &p\left( \int_{0}^{t}\zeta ^{\prime }\left( v_{\beta }\right)
\left\vert \psi ^{\prime }\right\vert ^{p}ds\right) ^{1/p}\left(
\int_{0}^{t}\zeta ^{\prime }\left( v_{\beta }\right) \left( \widetilde{v}%
_{\beta }\left( s,\phi _{\varepsilon }\right) \right) ^{p/(p-1)}ds\right)
^{(p-1)/p}.
\end{eqnarray*}%
Moreover, since %
$|D_{s}v_{\beta }|\leq p\left( v_{\beta }\left( s,\phi _{\varepsilon
}\right) \right) ^{(p-1)/p}\left( |D_{s}d|^{p}\ast \varphi _{\beta }\right)
^{1/p}$, 
we also have
\begin{eqnarray*}
I_{1} &=&\int_{0}^{t}\zeta ^{\prime }\left( v_{\beta }\right) \left\vert
D_{s}v_{\beta }\left( s,\phi _{\varepsilon }\right) \right\vert ds \\
&\leq &p\int_{0}^{t}\zeta ^{\prime }\left( v_{\beta }\right) \left( v_{\beta
}\left( s,\phi _{\varepsilon }\right) \right) ^{(p-1)/p}\left(
|D_{s}d|^{p}\ast \varphi _{\beta }\right) ^{1/p}ds \\
&\leq &p\left( \int_{0}^{t}\zeta ^{\prime }\left( v_{\beta }\right) v_{\beta
}\left( s,\phi _{\varepsilon }\right) ds\right) ^{(p-1)/p}\left(
\int_{0}^{t}\zeta ^{\prime }\left( v_{\beta }\right) \left( |D_{s}d|^{p}\ast
\varphi _{\beta }\right) ds\right) ^{1/p}.
\end{eqnarray*}%
Inserting the bounds for $I_{1}$ and $I_{2}$ into \eqref{eq:nyref} yields
\begin{eqnarray}
&&\frac{1}{p}\zeta \left( v_{\beta }\left( t,\phi _{\varepsilon }\left(
t\right) \right) \right) +\frac{\kappa }{\varepsilon p}\int_{0}^{t}\zeta
^{\prime }\left( v_{\beta }\right) d\left( s,\phi _{\varepsilon }\right)
\widetilde{v}_{\beta }\left( s,\phi _{\varepsilon }\right) ds
\label{betabound} \\
&\leq &\left( \int_{0}^{t}\zeta ^{\prime }\left( v_{\beta }\right) v_{\beta
}\left( s,\phi _{\varepsilon }\right) ds\right) ^{(p-1)/p}\left(
\int_{0}^{t}\zeta ^{\prime }\left( v_{\beta }\right) \left( |D_{s}d|^{p}\ast
\varphi _{\beta }\right) ds\right) ^{1/p}  \notag \\
&+&\left( \int_{0}^{t}\zeta ^{\prime }\left( v_{\beta }\right) \left(
\widetilde{v}_{\beta }\left( s,\phi _{\varepsilon }\right) \right)
^{p/(p-1)}ds\right) ^{(p-1)/p}\left( \int_{0}^{t}\zeta ^{\prime }\left(
v_{\beta }\right) \left\vert \psi ^{\prime }\right\vert ^{p}ds\right)
^{1/p}+\rho \left( \beta \right) ,  \notag
\end{eqnarray}%
where $\rho \left( \beta \right) =p^{-1}\zeta \left( v_{\beta }\left( 0,\phi
_{\varepsilon }\left( 0\right) \right) \right) \rightarrow 0$ as $\beta
\rightarrow 0$. By spatial Lipschitz continuity of $d(t,x)$ we have $%
v_{\beta }\left( s,\phi _{\varepsilon }\right) \rightarrow v\left( s,\phi
_{\varepsilon }\right) $ and $\widetilde{v}_{\beta }\left( s,\phi
_{\varepsilon }\right) \rightarrow \widetilde{v}\left( s,\phi _{\varepsilon
}\right) $ as $\beta \rightarrow 0$. Moreover, since $d$ satisfies %
\eqref{templip}, uniformly in space, we also have
\begin{equation*}
\int_{0}^{t}\left\vert D_{s}d\right\vert ^{p}ds\leq C(T)^{p},
\end{equation*}%
for some constant $C(T)$ independent of $x$. Therefore, by the
Fubini-Tonelli theorem we can conclude, since $D_{s}d(t,x)$ is jointly
measurable in $(t,x)$, that
\begin{equation*}
\int_{0}^{t}\left( |D_{s}d|^{p}\ast \varphi _{\beta }\right) ds=\int_{%
\mathbb{R}^{n}}\left( \int_{0}^{t}\left\vert D_{s}d\right\vert ^{p}ds\right)
\varphi _{\beta }(x-y)dy\leq C(T)^{p},
\end{equation*}%
and so
\begin{equation*}
\left( \int_{0}^{t}\zeta ^{\prime }\left( v_{\beta }\right) \left(
|D_{s}d|^{p}\ast \varphi _{\beta }\right) ds\right) ^{1/p}+\left(
\int_{0}^{t}\zeta ^{\prime }\left( v_{\beta }\right) \left\vert \psi
^{\prime }\right\vert ^{p}ds\right) ^{1/p}\leq C\left( T\right) <\infty ,
\end{equation*}%
since by construction $\left\vert \zeta ^{\prime }\left( v_{\beta }\right)
\right\vert \leq 1$, and $\psi \in \mathcal{C}^{1}\left( \left[ 0,T\right] ,%
\mathbb{R}^{n}\right) $. Thus, letting $\beta $ tend to zero in %
\eqref{betabound}, we obtain
\begin{eqnarray*}
&&\frac{1}{p}\zeta \left( v\left( t,\phi _{\varepsilon }\left( t\right)
\right) \right) +\frac{\kappa }{\varepsilon p}\int_{0}^{t}\zeta ^{\prime
}\left( v\left( s,\phi _{\varepsilon }\right) \right) v\left( s,\phi
_{\varepsilon }\right) ds \\
&\leq &C\left( T\right) \left( \int_{0}^{t}\zeta ^{\prime }\left( v\left(
s,\phi _{\varepsilon }\right) \right) v\left( s,\phi _{\varepsilon }\right)
ds\right) ^{\left( p-1\right) /p}.
\end{eqnarray*}%
Both terms on the left hand side are positive and each of the terms are
therefore bounded from above by the right hand side. Hence
\begin{equation*}
\frac{\kappa }{\varepsilon p}\left( \int_{0}^{t}\zeta ^{\prime }\left(
v\left( s,\phi _{\varepsilon }\right) \right) v\left( s,\phi _{\varepsilon
}\right) ds\right) ^{1/p}\leq C\left( T\right) ,
\end{equation*}%
and, as a consequence,
\begin{equation*}
\zeta \left( v\left( t,\phi _{\varepsilon }\left( t\right) \right) \right) +%
\frac{\kappa }{\varepsilon }\int_{0}^{t}\zeta ^{\prime }\left( v\left(
s,\phi _{\varepsilon }\right) \right) v\left( s,\phi _{\varepsilon }\right)
ds\leq K\left( T\right) \varepsilon ^{p-1}.
\end{equation*}%
We may assume that $\varepsilon >0$ has been chosen small enough such that $%
v\left( t,\phi _{\varepsilon }\left( t\right) \right) \leq \overline{\beta }%
^{p}/2$, for all $t\in \left[ 0,T\right] $. Then, by the definition of $%
\zeta $,
\begin{equation}
\frac{1}{\varepsilon ^{p-1}}\left( d\left( t,\phi _{\varepsilon }\left(
t\right) \right) \right) ^{p}+\frac{\kappa }{\varepsilon ^{p}}%
\int_{0}^{t}\left( d\left( s,\phi _{\varepsilon }\left( s\right) \right)
\right) ^{p}ds\leq K\left( T\right) ,  \label{d2bound}
\end{equation}%
for $t\in \left[ 0,T\right] $.

The remainder of the proof follows along the lines of the proof of Lemma 4.5
in \cite{DupuisIshii1993}, but we give the details for completeness.
Relation \eqref{d2bound} asserts that the sequences $\left\{ l_{\varepsilon
}\right\} _{\varepsilon >0}$ and $\left\{ \lambda _{\varepsilon }\right\}
_{\varepsilon \geq 0}$, where%
\begin{equation*}
l_{\varepsilon }\left( t\right) =\frac{1}{\varepsilon }d\left( t,\phi
_{\varepsilon }\left( t\right) \right) ,\quad \lambda _{\varepsilon }\left(
t\right) =\frac{1}{\varepsilon }\int_{0}^{t}d\left( s,\phi _{\varepsilon
}\left( s\right) \right) \gamma \left( s,\phi _{\varepsilon }\left( s\right)
\right) ds,
\end{equation*}%
are bounded in $L^{p}\left( \left[ 0,T\right] ,%
\mathbb{R}
\right) $ and $\mathcal{W}^{1,p}\left( \left[ 0,T\right] ,%
\mathbb{R}
^{n}\right) $ respectively. Thus, we may assume that $l_{\varepsilon }$ and $%
\lambda _{\varepsilon }$ converge weakly to $l\in L^{p}\left( \left[ 0,T%
\right] ,%
\mathbb{R}
\right) $ and $\lambda \in \mathcal{W}^{1,p}\left( \left[ 0,T\right] ,%
\mathbb{R}
^{n}\right) \subset \mathcal{C}\left( [0,T],\mathbb{R}^{n}\right) $,
respectively, as $\varepsilon \rightarrow 0$. Moreover, from \eqref{ode} we
conclude that $\phi _{\varepsilon }$ converges weakly to $\phi \in \mathcal{W%
}^{1,p}\left( \left[ 0,T\right] ,%
\mathbb{R}
^{n}\right) $ and that $\phi \left( t\right) =\psi \left( t\right) +\lambda
\left( t\right) $, $\phi \left( 0\right) =\psi \left( 0\right) $. This
proves \eqref{SP1} and, moreover, \eqref{SP2} holds due to \eqref{d2bound}.
By construction, $\lambda _{\varepsilon }^{\prime }\left( t\right)
=l_{\varepsilon }\left( t\right) \gamma \left( t,\phi _{\varepsilon }\left(
t\right) \right) $ and this implies that $\lambda ^{\prime }\left( t\right)
=l\left( t\right) \gamma \left( t,\phi \left( t\right) \right) $. Moreover,
if we let $\tau =\left\{ t\in \left[ 0,T\right] :\phi \left( t\right) \in
\Omega _{t}\right\} $ and note that for each fixed $t\in \tau $, we have $%
l_{\varepsilon }\left( t\right) =0$ for all sufficiently small $\varepsilon $
and hence $l\left( t\right) =0$ on $\tau $. Therefore%
\begin{equation*}
\left\vert \lambda \right\vert \left( t\right) =\int_{0}^{t}\left\vert
\lambda ^{\prime }\left( s\right) \right\vert ds=\int_{0}^{t}l\left(
s\right) \left\vert \gamma \left( s,\phi \left( s\right) \right) \right\vert
ds=\int_{0}^{t}l\left( s\right) ds,\quad \text{for all }t\in \left[ 0,T%
\right] ,
\end{equation*}%
as $\left\vert \gamma \left( s,\phi \left( s\right) \right) \right\vert =1$
for all $s\in \left[ 0,T\right] \setminus \tau $. This proves \eqref{SP3}.
In addition,
\begin{equation*}
\lambda \left( t\right) =\int_{0}^{t}l\left( s\right) \gamma \left( s,\phi
\left( s\right) \right) ds=\int_{0}^{t}\gamma \left( s,\phi \left( s\right)
\right) d\left\vert \lambda \right\vert \left( s\right) ,\quad \text{for all
}t\in \left[ 0,T\right] ,
\end{equation*}%
which proves \eqref{SP5}. It remains to verify \eqref{SP4}, but this follows
readily from
\begin{equation*}
\left\vert \lambda \right\vert \left( t\right) =\int_{0}^{t}l\left( s\right)
\left\vert \gamma \left( s,\phi \left( s\right) \right) \right\vert
ds=\int_{0}^{t}I_{\left\{ \phi \left( s\right) \in \partial \Omega
_{s}\right\} }l\left( s\right) ds=\int_{0}^{t}I_{\left\{ \phi \left(
s\right) \in \partial \Omega _{s}\right\} }d\left\vert \lambda \right\vert
\left( s\right) .
\end{equation*}%
We have completed the proof that $\left( \phi ,\lambda \right) \in \mathcal{W%
}^{1,p}\left( \left[ 0,T\right] ,%
\mathbb{R}
^{n}\right) \times \mathcal{W}^{1,p}\left( \left[ 0,T\right] ,%
\mathbb{R}
^{n}\right) $ solves the Skorohod problem for $\left( \Omega ,\gamma ,\psi
\right) $.\hfill $\Box $

\vspace{0.2cm}

The next step is to prove relative compactness of solutions to the Skorohod
problem. The proof follows the proof of Lemma 4.7 in \cite{DupuisIshii1993},
but a number of changes must be made carefully to handle the time dependency
of the domain.

\begin{lemma}
\label{compactest}Let $A$ be a compact subset of $\mathcal{C}\left( \left[
0,T\right] ,%
\mathbb{R}
^{n}\right) $. Then

\begin{description}
\item[(i)] There exists a constant $L<\infty $ such that%
\begin{equation*}
\left\vert \lambda \right\vert \left( T\right) <L,
\end{equation*}%
for all solutions $\left( \psi +\lambda ,\lambda \right) $ to the Skorohod
problem for $\left( \Omega ,\gamma ,\psi \right) $ with $\psi \in A$.

\item[(ii)] The set of $\phi $, such that $\left( \phi ,\lambda \right) $
solves the Skorohod problem for $\left( \Omega ,\gamma ,\psi \right) $ with $%
\psi \in A$, is relatively compact.
\end{description}
\end{lemma}

\noindent \textbf{Proof.} By the compactness of $\overline{\Omega }$ and the
continuity of $\gamma $, there exists a constant $c>0$ such that for every $%
t\in \left[ 0,T\right] $ and $x\in \overline{\Omega }_{t}\cap V$, where $V$
is the set defined in connection with \eqref{smooth_gamma}, there exists a
vector $v\left( t,x\right) $ and a set $\left[ t,t+c\right] \times B\left(
x,c\right) $ such that $\left\langle \gamma \left( s,y\right) ,v\left(
t,x\right) \right\rangle >c$ for all $\left( s,y\right) \in \left[ t,t+c%
\right] \times B\left( x,c\right) $. Without lack of generality, we may
assume that $c<\delta $, for the $\delta $ introduced in Remark \ref%
{spaceremark}. Let $\psi \in A$ be given and let $\left( \phi ,\lambda
\right) $ be any solution to the Skorohod problem for $\left( \Omega ,\gamma
,\psi \right) $. Define $T_{1}$ to be smallest of $T$, $c$ and $\inf \left\{
t\in \left[ 0,T\right] :\phi \left( t\right) \notin B\left( \phi \left(
0\right) ,c\right) \right\} $. Next define $T_{2}$ to be the smallest of $T$%
, $T_{1}+c$ and $\inf \left\{ t\in \left[ T_{1},T\right] :\phi \left(
t\right) \notin B\left( \phi \left( T_{1}\right) ,c\right) \right\} $.
Continuing in this fashion, we obtain a sequence $\left\{ T_{m}\right\}
_{m=1,2,...}$ of time instants. By construction, for all $s\in \left[
T_{m-1},T_{m}\right) $ we have $s\in \left[ T_{m-1},T_{m-1}+c\right] $ and $%
\phi \left( s\right) \in B\left( \phi \left( T_{m-1}\right) ,c\right) $. For
all $m$ such that $\phi \left( T_{m-1}\right) \in $ $\overline{\Omega }%
_{T_{m-1}}\cap V$, we have $\left\langle \gamma \left( s,\phi \left(
s\right) \right) ,v\left( T_{m-1},\phi \left( T_{m-1}\right) \right)
\right\rangle >c$ and hence%
\begin{eqnarray*}
&&\left\langle \phi \left( T_{m}\right) -\phi \left( T_{m-1}\right) ,v\left(
T_{m-1},\phi \left( T_{m-1}\right) \right) \right\rangle \\
&&-\left\langle \psi \left( T_{m}\right) -\psi \left( T_{m-1}\right)
,v\left( T_{m-1},\phi \left( T_{m-1}\right) \right) \right\rangle \\
&=&\int_{T_{m-1}}^{T_{m}}\left\langle \gamma \left( s,\phi \left( s\right)
\right) ,v\left( T_{m-1},\phi \left( T_{m-1}\right) \right) \right\rangle
d\left\vert \lambda \right\vert \left( s\right) \geq c\left( \left\vert
\lambda \right\vert \left( T_{m}\right) -\left\vert \lambda \right\vert
\left( T_{m-1}\right) \right) .
\end{eqnarray*}%
Since $A$ is compact, the set $\left\{ \psi \left( t\right) :t\in \left[ 0,T%
\right] ,\psi \in A\right\} $ is bounded. Moreover, since $\overline{\Omega }
$ is compact and $\phi \left( t\right) \in \overline{\Omega }_{t}$ for all $%
t\in \left[ 0,T\right] $, there exists a constant $M<\infty $ such that%
\begin{equation*}
\left\vert \lambda \right\vert \left( T_{m}\right) -\left\vert \lambda
\right\vert \left( T_{m-1}\right) <M.
\end{equation*}%
Note also that, for all $m$ such that $\phi \left( T_{m-1}\right) \notin $ $%
\overline{\Omega }_{T_{m-1}}\cap V$, we have, for $c$ sufficiently small,
that $\left\vert \lambda \right\vert \left( T_{m}\right) -\left\vert \lambda
\right\vert \left( T_{m-1}\right) =0$.

Define the modulus of continuity of a function $f\in \mathcal{C}\left( \left[
0,T\right] ,%
\mathbb{R}
^{n}\right) $ as $\left\Vert f\right\Vert _{s,t}=\sup_{s\leq t_{1}\leq
t_{2}\leq t}\left\vert f\left( t_{2}\right) -f\left( t_{1}\right)
\right\vert $ for $0\leq s\leq t\leq T$. We next prove that there exists a
positive constant $R$ such that, for any $\psi \in A$ and $T_{m-1}\leq \tau
\leq T_{m}$, it holds that
\begin{equation}
\left\Vert \lambda \right\Vert _{T_{m-1},\tau }\leq R\left( \left\Vert \psi
\right\Vert _{T_{m-1},\tau }^{1/2}+\left\Vert \psi \right\Vert
_{T_{m-1},\tau }^{3/2}+\left( \tau -T_{m-1}\right) ^{\widehat{\alpha }%
/2}\right) ,  \label{apriori}
\end{equation}%
where $\widehat{\alpha }$ is the H\"{o}lder exponent in Remark \ref%
{timeholder}. As we are only interested in the behaviour during the time
interval $\left[ T_{m-1},T_{m}\right] $, we simplify the notation by
setting, without loss of generality, $T_{m-1}=0$, $\phi \left(
T_{m-1}\right) =x$, $\psi \left( T_{m-1}\right) =x$, $\lambda \left(
T_{m-1}\right) =0$ and $\left\vert \lambda \right\vert \left( T_{m-1}\right)
=0$. Let $h$ be the function in Lemma \ref{testlemma3} and let $\chi ,C$ be
the corresponding positive constants. Define $B_{\varepsilon }\left(
t\right) =\varepsilon h\left( t,x,-\lambda \left( t\right) /\varepsilon
\right) $ and $E\left( t\right) =e^{-\left( 2\left\vert \lambda \right\vert
\left( t\right) +t\right) C/\chi }$. Since $h\left( t,x,0\right) =1$, we get%
\begin{eqnarray*}
B_{\varepsilon }\left( \tau \right) E\left( \tau \right) &=&B_{\varepsilon
}\left( 0\right) E\left( 0\right) +\int_{0}^{\tau }\left( E\left( u\right)
dB_{\varepsilon }\left( u\right) +B_{\varepsilon }\left( u\right) dE\left(
u\right) \right) \\
&=&\varepsilon +\int_{0}^{\tau }E\left( u\right) dB_{\varepsilon }\left(
u\right) -\frac{2C}{\chi }\int_{0}^{\tau }B_{\varepsilon }\left( u\right)
E\left( u\right) d\left\vert \lambda \right\vert \left( u\right) \\
&&-\frac{C}{\chi }\int_{0}^{\tau }B_{\varepsilon }\left( u\right) E\left(
u\right) du,
\end{eqnarray*}%
where the first integral can be rewritten as%
\begin{eqnarray*}
\int_{0}^{\tau }E\left( u\right) dB_{\varepsilon }\left( u\right)
&=&\int_{0}^{\tau }E\left( u\right) \varepsilon D_{t}h\left( u,x,-\lambda
\left( u\right) /\varepsilon \right) du \\
&&-\int_{0}^{\tau }E\left( u\right) \left\langle D_{p}h\left( u,x,-\lambda
\left( u\right) /\varepsilon \right) ,d\lambda \left( u\right) \right\rangle
.
\end{eqnarray*}%
By \eqref{testlemma31} and \eqref{testlemma35}, the integral involving $%
D_{t}h$ has the upper bound
\begin{eqnarray*}
&&\int_{0}^{\tau }E\left( u\right) \varepsilon D_{t}h\left( u,x,-\lambda
\left( u\right) /\varepsilon \right) du \\
&\leq &C\varepsilon \int_{0}^{\tau }E\left( u\right) \left\vert \lambda
\left( u\right) /\varepsilon \right\vert ^{2}du\leq \frac{C}{\chi }%
\int_{0}^{\tau }E\left( u\right) B_{\varepsilon }\left( u\right) du.
\end{eqnarray*}%
Next, we would like to find an upper bound for the integral involving $%
D_{p}h $ using \eqref{testlemma33} in some appropriate way, but we have to
be somewhat careful due to the temporal variation of the domain. Assume that
$\phi \left( u\right) \in \partial \Omega _{u}$. If $x\notin \overline{%
\Omega }_{u}$, there exists at least one point $y_{u}\in \overline{\Omega }%
_{u}\cap B\left( x,c\right) $ such that $\left\vert x-y_{u}\right\vert
=d\left( u,x\right) $. We have chosen $c<\delta $, so $\left\langle
y_{u}-\phi \left( u\right) ,\gamma \left( u,\phi \left( u\right) \right)
\right\rangle \geq -\theta \left\vert y_{u}-\phi \left( u\right) \right\vert
$ holds by Remark \ref{spaceremark} and, due to \eqref{testlemma33}, we can
conclude%
\begin{equation*}
I_{1}:=-\int_{0}^{\tau }E\left( u\right) \left\langle D_{p}h\left( u,\phi
\left( u\right) ,\left( y_{u}-\phi \left( u\right) \right) /\varepsilon
\right) ,\gamma \left( u,\phi \left( u\right) \right) \right\rangle
d\left\vert \lambda \right\vert \left( u\right) \leq 0,
\end{equation*}%
since $d\left\vert \lambda \right\vert \left( u\right) =0$ if $\phi \left(
u\right) \notin \partial \Omega _{u}$. If $x\in \overline{\Omega }_{u}$, the
above estimate holds with $y_{u}$ replaced by $x$. The integral involving $%
D_{p}h$ can be decomposed into%
\begin{equation*}
-\int_{0}^{\tau }E\left( u\right) \left\langle D_{p}h\left( u,x,-\lambda
\left( u\right) /\varepsilon \right) ,d\lambda \left( u\right) \right\rangle
=I_{1}+I_{2}+I_{3},
\end{equation*}%
for $I_{1}$ as above and%
\begin{equation*}
I_{2}=\int_{0}^{\tau }E\left( u\right) \left\langle D_{p}h\left( u,\phi
\left( u\right) ,-\lambda \left( u\right) /\varepsilon \right) -D_{p}h\left(
u,x,-\lambda \left( u\right) /\varepsilon \right) ,d \lambda
\left( u\right) \right\rangle ,
\end{equation*}%
\begin{equation*}
I_{3}=\int_{0}^{\tau }E\left( u\right) \left\langle D_{p}h\left( u,\phi
\left( u\right) ,\left( y_{u}-\phi \left( u\right) \right) /\varepsilon
\right) -D_{p}h\left( u,\phi \left( u\right) ,-\lambda \left( u\right)
/\varepsilon \right) ,d \lambda \left( u\right)
\right\rangle .
\end{equation*}%
By \eqref{testlemma31} and \eqref{testlemma37}, these integrals can be
bounded from above by%
\begin{eqnarray*}
I_{2} &\leq &\frac{C}{\varepsilon }\int_{0}^{\tau }E\left( u\right)
\left\vert \lambda \left( u\right) \right\vert \left\vert x-\phi \left(
u\right) \right\vert d\left\vert \lambda \right\vert \left( u\right) \\
&\leq &\frac{C}{\varepsilon }\int_{0}^{\tau }E\left( u\right) \left(
\left\vert \lambda \left( u\right) \right\vert ^{2}+\left\vert x-\psi \left(
u\right) \right\vert \left\vert \lambda \left( u\right) \right\vert \right)
d\left\vert \lambda \right\vert \left( u\right) \\
&\leq &\frac{2C}{\varepsilon }\int_{0}^{\tau }E\left( u\right) \left(
\left\vert \lambda \left( u\right) \right\vert ^{2}+\left\vert x-\psi \left(
u\right) \right\vert ^{2}\right) d\left\vert \lambda \right\vert \left(
u\right) \\
&\leq &\frac{2C}{\chi }\int_{0}^{\tau }E\left( u\right) B_{\varepsilon
}\left( u\right) d\left\vert \lambda \right\vert \left( u\right) +\frac{2C}{%
\varepsilon }\int_{0}^{\tau }E\left( u\right) \left\vert x-\psi \left(
u\right) \right\vert ^{2}d\left\vert \lambda \right\vert \left( u\right) ,
\end{eqnarray*}%
and%
\begin{eqnarray*}
I_{3} &\leq &\frac{C}{\varepsilon }\int_{0}^{\tau }E\left( u\right)
\left\vert y_{u}-\phi \left( u\right) -\left( -\lambda \left( u\right)
\right) \right\vert d\left\vert \lambda \right\vert \left( u\right) \\
&=&\frac{C}{\varepsilon }\int_{0}^{\tau }E\left( u\right) \left\vert
y_{u}-\psi \left( u\right) \right\vert d\left\vert \lambda \right\vert
\left( u\right) \\
&\leq &\frac{C}{\varepsilon }\int_{0}^{\tau }E\left( u\right) \left(
\left\vert x-\psi \left( u\right) \right\vert +\left\vert y_{u}-x\right\vert
\right) d\left\vert \lambda \right\vert \left( u\right) \\
&\leq &\frac{C}{\varepsilon }\int_{0}^{\tau }E\left( u\right) \left(
\left\vert x-\psi \left( u\right) \right\vert +d\left( u,x\right) \right)
d\left\vert \lambda \right\vert \left( u\right) .
\end{eqnarray*}%
Collecting all the terms, we obtain%
\begin{equation*}
B_{\varepsilon }\left( \tau \right) E\left( \tau \right) \leq \varepsilon +%
\frac{C}{\varepsilon }\int_{0}^{\tau }E\left( u\right) \left( \left\vert
x-\psi \left( u\right) \right\vert +2\left\vert x-\psi \left( u\right)
\right\vert ^{2}+d\left( u,x\right) \right) d\left\vert \lambda \right\vert
\left( u\right) ,
\end{equation*}%
which implies%
\begin{equation*}
B_{\varepsilon }\left( \tau \right) \leq \left( \frac{2C}{\varepsilon }%
\int_{0}^{\tau }E\left( u\right) \left( \left\Vert \psi \right\Vert _{0,\tau
}+\left\Vert \psi \right\Vert _{0,\tau }^{2}+K\tau ^{\widehat{\alpha }%
}\right) d\left\vert \lambda \right\vert \left( u\right) +\varepsilon
\right) e^{\left( 2\left\vert \lambda \right\vert \left( \tau \right) +\tau
\right) C/\chi },
\end{equation*}%
where $K$ and $\widehat{\alpha }$ are the constants from Remark \ref%
{timeholder}. Now%
\begin{equation*}
\int_{0}^{\tau }E\left( u\right) d\left\vert \lambda \right\vert \left(
u\right) \leq \int_{0}^{\tau }e^{-2C\left\vert \lambda \right\vert \left(
u\right) /\chi }d\left\vert \lambda \right\vert \left( u\right) \leq \frac{%
\chi }{2C},
\end{equation*}%
so%
\begin{equation*}
B_{\varepsilon }\left( \tau \right) \leq \left( \frac{\chi }{\varepsilon }%
\left( \left\Vert \psi \right\Vert _{0,\tau }+\left\Vert \psi \right\Vert
_{0,\tau }^{2}+K\tau ^{\widehat{\alpha }}\right) +\varepsilon \right)
e^{\left( 2\left\vert \lambda \right\vert \left( \tau \right) +\tau \right)
C/\chi }.
\end{equation*}%
Another application of \eqref{testlemma31} gives%
\begin{eqnarray*}
\left\vert \lambda \left( \tau \right) \right\vert &\leq &\frac{1}{2}\left(
\varepsilon +\frac{1}{\varepsilon }\left\vert \lambda \left( \tau \right)
\right\vert ^{2}\right) \leq \frac{\varepsilon }{2}+\frac{B_{\varepsilon
}\left( \tau \right) }{2\chi } \\
&\leq &\frac{\varepsilon }{2}+\left( \frac{1}{2\varepsilon }\left(
\left\Vert \psi \right\Vert _{0,\tau }+\left\Vert \psi \right\Vert _{0,\tau
}^{2}+K\tau ^{\widehat{\alpha }}\right) +\frac{\varepsilon }{2\chi }\right)
e^{\left( 2M+T\right) C/\chi }.
\end{eqnarray*}%
Set $\varepsilon =\max \left\{ \left\Vert \psi \right\Vert _{0,\tau
}^{1/2},\tau ^{\widehat{\alpha }/2}\right\} $ so that $\varepsilon \leq
\left\Vert \psi \right\Vert _{0,\tau }^{1/2}+\tau ^{\widehat{\alpha }/2}$, $%
1/\varepsilon \leq \left\Vert \psi \right\Vert _{0,\tau }^{-1/2}$ and $%
1/\varepsilon \leq \tau ^{-\widehat{\alpha }/2}$. Then \eqref{apriori}
follows immediately from the above inequality. By \eqref{apriori} and the
compactness of $A$, there exists a $\hat{\tau}>0$ such that
\begin{equation*}
\max \left\{ \left\Vert \psi \right\Vert _{T_{m-1},T_{m-1}+\hat{\tau}%
},\left\Vert \lambda \right\Vert _{T_{m-1},T_{m-1}+\hat{\tau}}\right\} \leq
c/3,
\end{equation*}%
which implies $\left\Vert \phi \right\Vert _{T_{m-1},T_{m-1}+\hat{\tau}}\leq
2c/3$. The definition of $\left\{ T_{m}\right\} $ then implies that $%
T_{m}-T_{m-1}\geq \min \left\{ \hat{\tau},c\right\} $. This proves (i) with $%
L=M\left( T/\min \left\{ \hat{\tau},c\right\} +1\right) $. Part (ii) follows
from \eqref{apriori} and the bound $T_{m}-T_{m-1}\geq \min \left\{ \hat{\tau}%
,c\right\} $.\hfill $\Box $

\vspace{0.2cm}

Equipped with the results above, we are now ready to state and prove the
existence of solutions to the Skorohod problem. The proof is very similar to
the proof of Theorem 4.8 in \cite{DupuisIshii1993}, so we only sketch the
first half of the proof.

\begin{lemma}
\label{contexist}Let $\psi \in \mathcal{C}\left( \left[ 0,T\right] ,%
\mathbb{R}
^{n}\right) $ with $\psi \left( 0\right) \in \overline{\Omega }_{0}$. Then
there exists a solution $\left( \phi ,\lambda \right) $ to the Skorohod
problem for $\left( \Omega ,\gamma ,\psi \right) $.
\end{lemma}

\noindent \textbf{Proof.} Let $\psi _{n}\in \mathcal{C}^{1}\left( \left[ 0,T%
\right] ,%
\mathbb{R}
^{n}\right) $ form a sequence of functions converging uniformly to $\psi $.
According to Lemma \ref{smoothexist}, there exists a solution $\left( \phi
_{n},\lambda _{n}\right) $ to the Skorohod problem for $\left( \Omega
,\gamma ,\psi _{n}\right) $. By Lemma \ref{compactest}, we may assume that
the sequence $\left\{ \lambda _{n}\right\} _{n=1}^{\infty }$ is equibounded
and equicontinuous, that is%
\begin{eqnarray*}
\sup_{n}\left\vert \lambda _{n}\right\vert \left( T\right) &\leq &L<\infty ,
\\
\lim_{\left\vert s-t\right\vert \rightarrow 0}\sup_{n}\left\vert \lambda
_{n}\left( s\right) -\lambda _{n}\left( t\right) \right\vert &=&0.
\end{eqnarray*}%
The Arzela-Ascoli theorem asserts the existence of a function $\lambda \in
\mathcal{C}\left( \left[ 0,T\right] ,%
\mathbb{R}
^{n}\right) $ such that $\left\{ \lambda _{n}\right\} $ converges uniformly
to $\lambda $. Clearly $\left\vert \lambda \right\vert \left( T\right) \leq
L $. Defining the function $\phi $ by $\phi =\psi +\lambda $, we conclude
that \eqref{SP1}-\eqref{SP3} of Definition \ref{skorohodprob} hold. To show
properties \eqref{SP4} and \eqref{SP5} in the same definition, we define the
measure $\mu _{n}$ on $\overline{\Omega }\times S\left( 0,1\right) $ as%
\begin{equation*}
\mu _{n}\left( A\right) =\int_{\left[ 0,T\right] }I_{\left\{ \left( s,\phi
_{n}\left( s\right) ,\gamma \left( s,\phi _{n}\left( s\right) \right)
\right) \in A\right\} }d\left\vert \lambda _{n}\right\vert \left( s\right) ,
\end{equation*}%
for every Borel set $A\subset \overline{\Omega }\times S\left( 0,1\right) $.
Introducing the notation $\overline{\Omega }_{\left[ 0,t\right] }:=\overline{%
\Omega }\cap \left( \left[ 0,t\right] \times
\mathbb{R}
^{n}\right) $, we have, by definition and \eqref{SP5},
\begin{equation*}
\left\vert \lambda _{n}\right\vert \left( t\right) =\mu _{n}\left( \overline{%
\Omega }_{\left[ 0,t\right] }\times S\left( 0,1\right) \right) ,
\end{equation*}%
and%
\begin{equation*}
\lambda _{n}\left( t\right) =\int_{\overline{\Omega }_{\left[ 0,t\right]
}\times S\left( 0,1\right) }\gamma d\mu _{n}\left( s,x,\gamma \right) ,
\end{equation*}%
for all $t\in \left[ 0,T\right] $. Since $\left\vert \lambda _{n}\right\vert
\left( T\right) \leq L<\infty $ for all $n$, the Banach-Alaoglu theorem
asserts that a subsequence of $\mu _{n}$ converges to some measure $\mu $
satisfying $\mu \left( \overline{\Omega }\times S\left( 0,1\right) \right)
<\infty $. By weak convergence and the continuity of $\lambda $,%
\begin{equation*}
\lambda \left( t\right) =\int_{\overline{\Omega }_{\left[ 0,t\right] }\times
S\left( 0,1\right) }\gamma d\mu \left( s,x,\gamma \right) .
\end{equation*}%
Using the fact that $\left( \phi _{n},\lambda _{n}\right) $ solves the
Skorohod problem for $\left( \Omega ,\gamma ,\psi _{n}\right) $, we can draw
several conclusions regarding the properties of the measure $\mu _{n}$ and
then use weak convergence of $\mu _{n}$ to $\mu $ to deduce that $\lambda $
satisfies \eqref{SP4} and \eqref{SP5}. This procedure is carried out in the
proofs of Theorem 2.8 in \cite{Costantini1992}, Theorem 4.8 in \cite%
{DupuisIshii1993} and Theorem 5.1 in \cite{NystromOnskog2010a}, so we omit
further details.\hfill $\Box $

\setcounter{equation}{0} \setcounter{theorem}{0}

\section{SDEs with oblique reflection\label{RSDE}}

Using the existence of solutions $\left( \phi ,\lambda \right) $ to the
Skorohod problem for $\left( \Omega ,\gamma ,\psi \right) $, with $\psi \in
\mathcal{C}\left( \left[ 0,T\right] ,%
\mathbb{R}
^{n}\right) $ and $\psi \left( 0\right) \in \overline{\Omega }_{0}$, we can
now prove existence and uniqueness of solutions to SDEs with oblique
reflection at the boundary of a bounded, time-dependent domain. To this end,
assume that the triple $\left( X,Y,k\right) $ satisfies%
\begin{equation*}
Y\left( t\right) =x+\int_{0}^{t}b\left( s,X\left( s\right) \right)
ds+\int_{0}^{t}\sigma \left( s,X\left( s\right) \right) dM\left( s\right)
+k\left( t\right) ,
\end{equation*}%
\begin{equation*}
X\left( t\right) \in \overline{\Omega }_{t},\quad Y\left( t\right) \in
\overline{\Omega }_{t},
\end{equation*}%
\begin{equation*}
\left\vert k\right\vert \left( t\right) =\int_{\left( 0,t\right] }I_{\left\{
Y\left( s\right) \in \partial \Omega _{s}\right\} }d\left\vert k\right\vert
\left( s\right) <\infty ,\quad k\left( t\right) =\int_{\left( 0,t\right]
}\gamma \left( s\right) d|k|\left( s\right) ,
\end{equation*}%
where $x\in \overline{\Omega }_{0}$ is fixed, $\gamma \left( s\right)
=\gamma \left( s,Y\left( s\right) \right) $ $d\left\vert k\right\vert $%
-a.s.~and $M$ is a continuous $\mathcal{F}_{t}$-martingale satisfying%
\begin{equation}
d\left\langle M_{i},M_{j}\right\rangle \left( t\right) \leq Cdt,
\label{mart}
\end{equation}%
for some $C\in \left( 0,\infty \right) $. Let $\left( X^{\prime },Y^{\prime
},k^{\prime }\right) $ be a similar triple, 
but with $x$ replaced by $x^{\prime}\in \overline{\Omega }_{0}$, and $%
\gamma^{\prime} \left( s\right) =\gamma \left( s,Y^{\prime}\left( s\right)
\right)$ $d\left\vert k^{\prime}\right\vert $-a.s.

We shall prove uniqueness of solutions by a Picard iteration scheme and a
crucial ingredient is then the estimate provided by the following theorem.
Note that Lemma \ref{rsdetheorem} holds for a general continuous $\mathcal{F}%
_{t}$-martingale satisfying \eqref{mart}, whereas in Theorem \ref{main} we
restrict our interest to $M$ being a standard Wiener process.

\begin{lemma}
\label{rsdetheorem}There exists a positive constant $C<\infty $ such that%
\begin{equation*}
E\left[ \sup_{0\leq s\leq t}\left\vert Y\left( s\right) -Y^{\prime }\left(
s\right) \right\vert ^{%
{\acute{}}%
4}\right] \leq C\left( \left\vert x-x^{\prime }\right\vert ^{4}+\int_{0}^{t}E%
\left[ \sup_{0\leq u\leq s}\left\vert X\left( u\right) -X^{\prime }\left(
u\right) \right\vert ^{4}\right] ds\right) .
\end{equation*}
\end{lemma}

\noindent \textbf{Proof.} Fix $\varepsilon >0$, let $\lambda >0$ be a
constant to be specified later, and let $w_{\varepsilon }\in \mathcal{C}%
^{1,2}\left( \left[ 0,T\right] \times
\mathbb{R}
^{n}\times
\mathbb{R}
^{n},%
\mathbb{R}
\right) $ and $\alpha \in \mathcal{C}^{1,2}\left( \overline{\Omega },%
\mathbb{R}
\right) $ and be the functions defined in Lemma \ref{testlemma4}-\ref%
{testlemma5}. Define the stopping time%
\begin{equation*}
\tau =\inf \left\{ s\in \left[ 0,T\right] :\left\vert Y\left( s\right)
-Y^{\prime }\left( s\right) \right\vert \geq \delta \right\} ,
\end{equation*}%
where $\delta >0$ is the constant from Remark \ref{spaceremark}. Let $B$
denote the diameter of the smallest ball containing $\bigcup\nolimits_{t\in %
\left[ 0,T\right] }\overline{\Omega }_{t}$. Then, assuming without loss of
generality that $B/\delta \geq 1$, we have%
\begin{equation*}
E\left[ \sup_{0\leq s\leq t}\left\vert Y\left( s\right) -Y^{\prime }\left(
s\right) \right\vert ^{4}\right] \leq \left( \frac{B}{\delta }\right) ^{4}E%
\left[ \sup_{0\leq s\leq t\wedge \tau }\left\vert Y\left( s\right)
-Y^{\prime }\left( s\right) \right\vert ^{4}\right] ,
\end{equation*}%
so it is sufficient to prove the theorem for $t\wedge \tau $. To simplify
the notation, however, we write $t$ in place of $t\wedge \tau $ and assume
that $\left\vert Y\left( s\right) -Y^{\prime }\left( s\right) \right\vert
<\delta $ in the proof below.

Define, for all $\left( t,x,y\right) $ such that $\left( t,x\right) ,\left(
t,y\right) \in \overline{\Omega }$, the function $v$ as%
\begin{equation*}
v\left( t,x,y\right) =e^{-\lambda \left( \alpha \left( t,x\right) +\alpha
\left( t,y\right) \right) }w_{\varepsilon }\left( t,x,y\right) :=u\left(
t,x,y\right) w_{\varepsilon }\left( t,x,y\right) .
\end{equation*}%
The regularity of $v$ is inherited from that of $w_{\varepsilon }$ and $%
\alpha $. By It\={o}'s formula we have, suppressing the $s$-dependence for $%
X $, $X^{\prime }$, $Y$ and $Y^{\prime }$,%
\begin{eqnarray}
&&v\left( t,Y\left( t\right) ,Y^{\prime }\left( t\right) \right)  \label{Ito}
\\
&=&v\left( 0,x,x^{\prime }\right) +\int_{0}^{t}D_{s}v\left( s,Y,Y^{\prime
}\right) ds  \notag \\
&&+\int_{0}^{t}\left\langle D_{x}v\left( s,Y,Y^{\prime }\right) ,b\left(
s,X\right) \right\rangle ds+\int_{0}^{t}\left\langle D_{y}v\left(
s,Y,Y^{\prime }\right) ,b\left( s,X^{\prime }\right) \right\rangle ds  \notag
\\
&&+\int_{0}^{t}\left\langle D_{x}v\left( s,Y,Y^{\prime }\right) ,\sigma
\left( s,X\right) dM\left( s\right) \right\rangle +\int_{0}^{t}\left\langle
D_{y}v\left( s,Y,Y^{\prime }\right) ,\sigma \left( s,X^{\prime }\right)
dM\left( s\right) \right\rangle  \notag \\
&&+\int_{0}^{t}\left\langle D_{x}v\left( s,Y,Y^{\prime }\right) ,\gamma
\left( s\right) \right\rangle d\left\vert k\right\vert \left( s\right)
+\int_{0}^{t}\left\langle D_{y}v\left( s,Y,Y^{\prime }\right) ,\gamma
^{\prime }\left( s\right) \right\rangle d\left\vert k^{\prime }\right\vert
\left( s\right)  \notag \\
&&+\int_{0}^{t}\text{tr}\left( \left(
\begin{array}{c}
\sigma \left( s,X\right) \\
\sigma \left( s,X^{\prime }\right)%
\end{array}%
\right) ^{T}D^{2}v\left( s,Y,Y^{\prime }\right) \left(
\begin{array}{c}
\sigma \left( s,X\right) \\
\sigma \left( s,X^{\prime }\right)%
\end{array}%
\right) d\left\langle M\right\rangle \left( s\right) \right) .  \notag
\end{eqnarray}%
We define the martingale $N$ as
\begin{equation*}
N\left( t\right) =\int_{0}^{t}\left\langle D_{x}v\left( s,Y,Y^{\prime
}\right) ,\sigma \left( s,X\right) dM\left( s\right) \right\rangle
+\int_{0}^{t}\left\langle D_{y}v\left( s,Y,Y^{\prime }\right) ,\sigma \left(
s,X^{\prime }\right) dM\left( s\right) \right\rangle ,
\end{equation*}%
and simplify the remaining terms in \eqref{Ito}. From \eqref{testlemma42}, %
\eqref{testlemma45} and the regularity of $u$, we have
\begin{equation*}
\int_{0}^{t}D_{s}v\left( s,Y,Y^{\prime }\right) ds\leq C\left( \lambda
\right) \int_{0}^{t}\left( \varepsilon +\frac{\left\vert Y-Y^{\prime
}\right\vert ^{2}}{\varepsilon }\right) ds.
\end{equation*}%
Similarly, following the proof of Theorem 5.1 in \cite{DupuisIshii1993}, we
have%
\begin{eqnarray}
&&\int_{0}^{t}\left\langle D_{x}v\left( s,Y,Y^{\prime }\right) ,b\left(
s,X\right) \right\rangle ds+\int_{0}^{t}\left\langle D_{y}v\left(
s,Y,Y^{\prime }\right) ,b\left( s,X^{\prime }\right) \right\rangle ds
\label{driftterms} \\
&\leq &C\left( \lambda \right) \left( \varepsilon +\int_{0}^{t}\frac{%
\left\vert Y-Y^{\prime }\right\vert ^{2}}{\varepsilon }ds+\int_{0}^{t}\frac{%
\left\vert X-X^{\prime }\right\vert ^{2}}{\varepsilon }ds\right) .  \notag
\end{eqnarray}%
A simple extension of Lemma 5.7 in \cite{DupuisIshii1993} to the
time-dependent case shows that there exists a constant $K_{1}\left( \lambda
\right) <\infty $ such that for all $t\in \left[ 0,T\right] $, $x,y\in
\overline{\Omega }_{t}$, the second order derivatives of $v$ with respect to
the spatial variables satisfy%
\begin{equation*}
D^{2}v\left( t,x,y\right) \leq K_{1}\left( \lambda \right) \left( \frac{1}{%
\varepsilon }\left(
\begin{array}{cc}
I & -I \\
-I & I%
\end{array}%
\right) +\left( \varepsilon +\frac{\left\vert x-y\right\vert ^{2}}{%
\varepsilon }\right) \left(
\begin{array}{cc}
I & 0 \\
0 & I%
\end{array}%
\right) \right) .
\end{equation*}%
Moreover, it is an easy consequence of the Lipschitz continuity of $\sigma $
that there exists a constant $K_{2}\left( \lambda \right) <\infty $ such
that for all $t\in \left[ 0,T\right] $, $x,y,\xi ,\omega \in \overline{%
\Omega }_{t}$,%
\begin{equation*}
\left(
\begin{array}{c}
\sigma \left( t,\xi \right) \\
\sigma \left( t,\omega \right)%
\end{array}%
\right) ^{T}D^{2}v\left( t,x,y\right) \left(
\begin{array}{c}
\sigma \left( t,\xi \right) \\
\sigma \left( t,\omega \right)%
\end{array}%
\right) \leq K_{2}\left( \lambda \right) \left( \varepsilon +\frac{1}{%
\varepsilon }\left( \left\vert \xi -\omega \right\vert ^{2}+\left\vert
x-y\right\vert ^{2}\right) \right) I.
\end{equation*}%
Consequently, the last term in \eqref{Ito} may be simplified to%
\begin{eqnarray*}
&&\int_{0}^{t}\text{tr}\left( \left(
\begin{array}{c}
\sigma \left( s,X\right) \\
\sigma \left( s,X^{\prime }\right)%
\end{array}%
\right) ^{T}D^{2}v\left( s,Y,Y^{\prime }\right) \left(
\begin{array}{c}
\sigma \left( s,X\right) \\
\sigma \left( s,X^{\prime }\right)%
\end{array}%
\right) d\left\langle M\right\rangle \left( s\right) \right) \\
&\leq &C\left( \lambda \right) \left( \varepsilon +\int_{0}^{t}\frac{%
\left\vert X-X^{\prime }\right\vert ^{2}}{\varepsilon }ds+\int_{0}^{t}\frac{%
\left\vert Y-Y^{\prime }\right\vert ^{2}}{\varepsilon }ds\right) .
\end{eqnarray*}%
Considering now the terms containing $\left\vert k\right\vert $ and $%
\left\vert k^{\prime }\right\vert $, we see, following the proof of Theorem
5.1 in \cite{DupuisIshii1993}, that
\begin{eqnarray*}
&&\int_{0}^{t}\left\langle D_{x}v\left( s,Y,Y^{\prime }\right) ,\gamma
\left( s\right) \right\rangle d\left\vert k\right\vert \left( s\right)
+\int_{0}^{t}\left\langle D_{y}v\left( s,Y,Y^{\prime }\right) ,\gamma
^{\prime }\left( s\right) \right\rangle d\left\vert k^{\prime }\right\vert
\left( s\right) \\
&\leq &C\int_{0}^{t}u\left( s,Y,Y^{\prime }\right) \frac{\left\vert
Y-Y^{\prime }\right\vert ^{2}}{\varepsilon }d\left\vert k\right\vert \left(
s\right) +C\int_{0}^{t}u\left( s,Y,Y^{\prime }\right) \frac{\left\vert
Y-Y^{\prime }\right\vert ^{2}}{\varepsilon }d\left\vert k^{\prime
}\right\vert \left( s\right) \\
&&-\lambda \int_{0}^{t}v\left( s,Y,Y^{\prime }\right) \left\langle
D_{x}\alpha \left( s,Y\right) ,\gamma \left( s\right) \right\rangle
d\left\vert k\right\vert \left( s\right) \\
&&-\lambda \int_{0}^{t}v\left( s,Y,Y^{\prime }\right) \left\langle
D_{x}\alpha \left( s,Y^{\prime }\right) ,\gamma ^{\prime }\left( s\right)
\right\rangle d\left\vert k^{\prime }\right\vert \left( s\right) .
\end{eqnarray*}%
Moreover, \eqref{testlemma41} and \eqref{alfaprop} give, since $d\left\vert
k\right\vert \left( s\right) $ is zero unless $Y\left( s\right) \in \partial
\Omega _{s}$,
\begin{equation*}
-\lambda v\left( s,Y,Y^{\prime }\right) \left\langle D_{x}\alpha \left(
s,Y\right) ,\gamma \left( s\right) \right\rangle \leq -\lambda \chi u\left(
s,Y,Y^{\prime }\right) \frac{\left\vert Y-Y^{\prime }\right\vert ^{2}}{%
\varepsilon },
\end{equation*}%
so, by putting $\lambda =C/\chi $ all integrals with respect to $\left\vert
k\right\vert $ and $\left\vert k^{\prime }\right\vert $ vanish. Dropping the
$\lambda $-dependence from the constants, \eqref{testlemma41} and \eqref{Ito}
give%
\begin{eqnarray*}
\frac{1}{C}\frac{\left\vert Y\left( t\right) -Y^{\prime }\left( t\right)
\right\vert ^{2}}{\varepsilon } &\leq &v\left( t,Y\left( t\right) ,Y^{\prime
}\left( t\right) \right) \leq v\left( 0,x,x^{\prime }\right) +\varepsilon
+N\left( t\right) \\
&&+\int_{0}^{t}\frac{\left\vert Y-Y^{\prime }\right\vert ^{2}}{\varepsilon }%
ds+\int_{0}^{t}\frac{\left\vert X-X^{\prime }\right\vert ^{2}}{\varepsilon }%
ds.
\end{eqnarray*}%
Now applying \eqref{testlemma42} to $v\left( 0,x,x^{\prime }\right) $,
multiplying by $\varepsilon $, squaring, taking supremum and expectations on
both sides, we obtain%
\begin{eqnarray*}
E\left[ \sup_{0\leq s\leq t}\left\vert Y\left( s\right) -Y^{\prime }\left(
s\right) \right\vert ^{4}\right] &\leq &C\left( \left\vert x-x^{\prime
}\right\vert ^{4}+\varepsilon ^{4}+\varepsilon ^{2}E\left[ \sup_{0\leq s\leq
t}\left( N\left( s\right) \right) ^{2}\right] \right. \\
&&\left. +\int_{0}^{t}E\left[ \left\vert X-X^{\prime }\right\vert
^{4}+\left\vert Y-Y^{\prime }\right\vert ^{4}\right] ds\right) .
\end{eqnarray*}%
Then proceeding as in \eqref{driftterms}, the Doob-Kolmogorov inequality
gives%
\begin{eqnarray*}
E\left[ \sup_{0\leq s\leq t}\left( N\left( s\right) \right) ^{2}\right]
&\leq &4E\left[ \left( N\left( t\right) \right) ^{2}\right] \\
&\leq &C\int_{0}^{t}\left( \varepsilon ^{2}+E\left[ \frac{\left\vert
Y-Y^{\prime }\right\vert ^{4}}{\varepsilon ^{2}}+\frac{\left\vert
X-X^{\prime }\right\vert ^{4}}{\varepsilon ^{2}}\right] \right) ds,
\end{eqnarray*}%
Letting $\varepsilon $ tend to zero, we obtain%
\begin{equation*}
E\left[ \sup_{0\leq s\leq t}\left\vert Y\left( s\right) -Y^{\prime }\left(
s\right) \right\vert ^{4}\right] \leq C\left( \left\vert x-x^{\prime
}\right\vert ^{4}+\int_{0}^{t}E\left[ \left( \left\vert X-X^{\prime
}\right\vert ^{4}+\left\vert Y-Y^{\prime }\right\vert ^{4}\right) \right]
ds\right) ,
\end{equation*}%
from which the requested inequality follows by a simple application of
Gronwall's inequality.\hfill $\Box $

\vspace{0.2cm}

\noindent \textbf{Proof of Theorem \ref{main}.}
%
Given Lemma \ref{contexist} and Lemma \ref{rsdetheorem}, the proof of
Theorem \ref{main} follows exactly along the lines of the proof of Corollary
5.2 in \cite{DupuisIshii1993}, which in turn follows the same outline as in
\cite{LionsSznitman1984}, Theorem 4.3. Note that the main problem is
verifying the adaptedness property of the solutions to the reflected SDE.
This property follows from an approximation of continuous $\mathcal{F}_t$%
-adapted semimartingales by bounded variation processes, for which one can
show existence of unique bounded variation solutions to the Skorohod
problem, and these bounded variation solutions will be $\mathcal{F}_t$%
-adapted. We omit further details. \hfill $\Box $

\setcounter{equation}{0} \setcounter{theorem}{0}

\section{Fully nonlinear second-order parabolic PDEs\label{PDE}}

In this section, we prove the results on partial differential equations.
First, we recall the definition of viscosity solutions. Let $E\subset
\mathbb{R}^{n+1}$ be arbitrary. If $u:E\rightarrow \mathbb{R}$, then the
parabolic superjet $\mathcal{P}_{E}^{2,+}u\left(s,z\right)$ contains all triplets $%
\left( a,p,X\right) \in \mathbb{R}\times \mathbb{R}^{n}\times \mathbb{S}^{n}$
such that if $\left(s,z\right)\in E$ then
\begin{align*}
u\left(t,x\right)& \leq u\left(s,z\right)+a\left(t-s\right)+\langle p,x-z\rangle +\frac{1}{2}\langle
X\left(x-z\right),x-z\rangle  \\
& +o\left(|t-s|+|x-z|^{2}\right)\quad \text{as }E\ni \left(t,x\right)\rightarrow \left(s,z\right).
\end{align*}%
The parabolic subjet is defined as $\mathcal{P}_{E}^{2,-}u\left(s,z\right)=-\mathcal{P}%
_{E}^{2,+}\left(-u\left(s,z\right)\right)$. The closures $\overline{\mathcal{P}}_{E}^{2,+}u\left(s,z\right)$
and $\overline{\mathcal{P}}_{E}^{2,-}u\left(s,z\right)$ are defined in analogue with
(2.6) and (2.7) in \cite{CrandallIshiiLions1992}.
A function $u\in USC(\widetilde{\Omega })$ is a \textit{viscosity subsolution}
of \eqref{huvudekvationen} in $\Omega ^{\circ }$ if, for all $\left(a,p,A\right)\in
\mathcal{P}_{\widetilde{\Omega }}^{2,+}u\left(t,x\right)$, it holds that
\begin{equation*}
a+F\left(t,x,u\left(t,x\right),p,A\right)\leq 0,\quad \text{for}\;\left(t,x\right)\in \Omega ^{\circ }.
\end{equation*}%
If, in addition, for $\left(t,x\right)\in \partial \Omega $ it holds that
\begin{equation}
\min \{a+F\left(t,x,u\left(t,x\right),p,A\right),\;\langle p,\widetilde{\gamma }\left(t,x\right)\rangle
+f\left(t,x,u\left(t,x\right)\right)\}\leq 0,  \label{eq:BC_viscosity_sub}
\end{equation}%
then $u$ is a viscosity subsolution of \eqref{huvudekvationen}-%
\eqref{randvillkor} in $\widetilde{\Omega }$. Similarly, a function $v\in
LSC(\widetilde{\Omega })$ is a \textit{viscosity supersolution} of %
\eqref{huvudekvationen} in $\Omega ^{\circ }$ if, for all $\left(a,p,A\right)\in
\mathcal{P}_{\widetilde{\Omega }}^{2,-}v\left(t,x\right)$, it holds that
\begin{equation*}
a+F\left(t,x,v\left(t,x\right),p,A\right)\geq 0,\quad \text{for}\;\left(t,x\right)\in \Omega ^{\circ }.
\end{equation*}%
If, in addition, for $\left(t,x\right)\in \partial \Omega $ it holds that
\begin{equation}
\max \{a+F\left(t,x,v\left(t,x\right),p,A\right),\;\langle p,\widetilde{\gamma }\left(t,x\right)\rangle
+f\left(t,x,v\left(t,x\right)\right)\}\geq 0,  \label{eq:BC_viscosity_sup}
\end{equation}%
then $v$ is a viscosity supersolution of \eqref{huvudekvationen}-%
\eqref{randvillkor} in $\widetilde{\Omega }$. A function is a \textit{%
viscosity solution} if it is both a viscosity subsolution and a viscosity
supersolution. We remark that in the definition of viscosity solutions
above, we may replace $\mathcal{P}_{\widetilde{\Omega }}^{2,+}u\left(t,x\right)$ and $%
\mathcal{P}_{\widetilde{\Omega }}^{2,-}v\left(t,x\right)$ by $\overline{\mathcal{P}}_{%
\widetilde{\Omega }}^{2,+}u\left(t,x\right)$ and $\overline{\mathcal{P}}_{\widetilde{%
\Omega }}^{2,-}v\left(t,x\right)$, respectively. In the following, we often skip
writing \textquotedblleft viscosity" before subsolutions, supersolutions and
solutions. Note also that, given any set $E\subset \mathbb{R}^{n+1}$ and $%
t\in \lbrack 0,T]$, we denote, in the following, the time sections of $E$ as
$E_{t}=\{x:\left(t,x\right)\in E\}$.

Next we give two lemmas. The first clarifies that the maximum principle for
semicontinuous functions \cite{CrandallIshii1990}, \cite%
{CrandallIshiiLions1992}, holds true in time-dependent domains.

\begin{lemma}
\label{le:timdep_max} Suppose that $\mathcal{O}^{i}=\mathcal{\widehat{O}}%
^{i}\cap \left( \left(0,T\right)\times \mathbb{R}^{n}\right) $ for $i=1,\dots ,k$ where
$\mathcal{\widehat{O}}^{i}$ are locally compact subsets of $\mathbb{R}^{n+1}$%
. Assume that $u_{i}\in USC(\mathcal{O}^{i})$ and let $\varphi
:\left(t,x_{1},\dots ,x_{k}\right)\rightarrow \varphi \left(t,x_{1},\dots ,x_{k}\right)$ be
defined on an open neighborhood of $\{\left(t,x\right):t\in \left(0,T\right)\;\text{and}\;x_{i}\in
\mathcal{O}_{t}^{i}\;\text{for}\;i=1,\dots ,k\}$ and such that $\varphi $ is
once continuously differentiable in $t$ and twice continuously
differentiable in $\left(x_{1},\dots ,x_{k}\right)$. Suppose that $s\in \left(0,T\right)$ and $%
z_{i}\in \mathcal{O}_{s}^{i}$ and
\begin{equation*}
w\left(t,x_{1},\dots ,x_{k}\right)\equiv u_{1}\left(t,x_{1}\right)+\dots +u_{k}\left(t,x_{k}\right)-\varphi
\left(t,x_{1},\dots ,x_{k}\right)\leq w\left(s,z_{1},\dots ,z_{k}\right),
\end{equation*}%
for $0<t<T$ and $x_{i}\in \mathcal{O}_{t}^{i}$. Assume, moreover, that there
is an $r>0$ such that for every $M>0$ there is a $C$ such that, for $%
i=1,\dots ,k$,
\begin{align}
& b_{i}\leq C,\;\text{whenever }\left( b_{i},q_{i},X_{i}\right) \in \mathcal{%
P}_{\mathcal{O}^{i}}^{2,+}u_{i}\left(t,x\right)\text{ with }\left\Vert X_{i}\right\Vert
\leq M\text{ and}  \notag  \label{eq:besvarlig_assumption} \\
&
|x_{i}-z_{i}|+|t-s|+|u_{i}\left(t,x_{i}\right)-u_{i}\left(s,z_{i}\right)|+|q_{i}-D_{x_{i}}\varphi
\left(s,z_{1},\dots ,z_{k}\right)|\leq r.
\end{align}%
Then, for each $\varepsilon >0$ there exist $\left(b_{i},X_{i}\right)$ such that
\begin{equation*}
\left(b_{i},D_{x_{i}}\varphi \left(s,z_{1},\dots ,z_{k}\right),X_{i}\right)\in \overline{\mathcal{P%
}}_{\mathcal{O}^{i}}^{2,+}u_{i}\left(s,z\right),\quad \text{for}\;i=1,...,k,
\end{equation*}%
\begin{equation*}
-\left( \frac{1}{\varepsilon }+||A||\right) I\leq \left(
\begin{array}{ccc}
X_{1} & \dots & 0 \\
\vdots & \ddots & \vdots \\
0 & \dots & X_{k}%
\end{array}%
\right) \leq A+\varepsilon A^{2},
\end{equation*}%
and
\begin{equation*}
b_{1}+\dots +b_{k}=D_{t}\varphi \left(s,z_{1},\dots ,z_{k}\right),
\end{equation*}%
where $A=\left( D_{x}^{2}\varphi \right) \left(s,z_{1},\dots ,z_{k}\right)$.
\end{lemma}

\noindent \textbf{Proof.} Following ideas from page 1008\textbf{\ }in \cite%
{CrandallIshii1990} we let $K_{i}$ be compact neighborhoods of $\left(s,z\right)$ in $%
\mathcal{O}^{i}$ and define the extended functions $\widetilde{u}_{1},\dots ,%
\widetilde{u}_{k}$, $\widetilde{u}_{i}\in USC\left(\mathbb{R}^{n}\right)$ for $%
i=1,\dots ,k$, by
\begin{equation*}
\widetilde{u}_{i}\left(t,x\right)=\left\{
\begin{array}{rl}
u\left(t,x\right), & \text{if}\quad \left(t,x\right)\in K_{i}, \\
-\infty , & \text{otherwise.}%
\end{array}%
\right.
\end{equation*}%
From the definitions of sub and superjets it follows, for $i=1,\dots ,k$,
that
\begin{equation}
\mathcal{P}_{\mathbb{R}^{n+1}}^{2,+}\widetilde{u}_{i}\left(t,x\right)=\mathcal{P}_{%
\mathcal{O}^{i}}^{2,+}u_{i}\left(t,x\right),  \label{eq:subjet_lika_1}
\end{equation}%
for $\left(t,x\right)$ in the interior of $K_{i}$ relative to $\mathcal{O}^{i}$.
Excluding the trivial case $u_{i}\left(t,x\right)=-\infty $, then the function $%
\widetilde{u}_{i}\left(t,x\right)$ cannot approach $u_{i}\left(s,z\right)$ unless $\left(t,x\right)\in K_{i}$
and it follows that
\begin{equation}
\overline{\mathcal{P}}_{\mathbb{R}^{n+1}}^{2,+}\widetilde{u}_{i}\left(t,x\right)=%
\overline{\mathcal{P}}_{\mathcal{O}^{i}}^{2,+}u_{i}\left(t,x\right).
\label{eq:subjet_lika_2}
\end{equation}%
Setting $\widetilde{w}\left(t,x_{1},\dots ,x_{k}\right)=\widetilde{u}%
_{1}\left(t,x_{1}\right)+\dots +\widetilde{u}_{k}\left(t,x_{k}\right)$ we see that $\left(s,z_{1},\dots
,z_{k}\right)$ is also a maximum of the function $\left(\widetilde{w}-\varphi
\right)\left(t,x_{1},\dots ,x_{k}\right)$. Moreover, we note that the proof of Lemma 8 in
\cite{CrandallIshii1990} still works if (27) in \cite{CrandallIshii1990} is
replaced by assumption \eqref{eq:besvarlig_assumption}. These facts,
together with \eqref{eq:subjet_lika_1} and \eqref{eq:subjet_lika_2}, allows
us to complete the proof of Lemma \ref{le:timdep_max} by using Theorem 7 in
\cite{CrandallIshii1990}. \hfill $\Box $

\vspace{0.2cm}

Before proving the next lemma, let us note that standard arguments imply
that we can assume $\lambda >0$ in \eqref{ass_F_nondecreasing}. Indeed, if $%
\lambda \leq 0$ then for $\bar{\lambda}<\lambda $ the functions $e^{\bar{%
\lambda}t}u\left(t,x\right)$ and $e^{\bar{\lambda}t}v\left(t,x\right)$ are, respectively, sub- and
supersolutions of \eqref{huvudekvationen}-\eqref{randvillkor} with $%
F\left(t,x,r,p,X\right)$ and $f\left(t,x,r\right)$ replaced by
\begin{equation}
-\bar{\lambda}r+e^{\bar{\lambda}t}F\left(t,x,e^{-\bar{\lambda}t}r,e^{-\bar{\lambda%
}t}p,e^{-\bar{\lambda}t}X\right)\quad \text{and}\quad e^{\bar{\lambda}t}f\left(t,x,e^{-%
\bar{\lambda}t}r\right).  \label{assume_lambda_positive}
\end{equation}%
Hence, in the following proof we assume $\lambda >0$ in %
\eqref{ass_F_nondecreasing}. Next we prove the following version of the
comparison principle.

\begin{lemma}
\label{maxrand} Let $\Omega ^{\circ }$ be a time-dependent domain satisfying %
\eqref{timesect}. Assume \eqref{ass_F_cont}-\eqref{ass_F_XY}. Let $u\in
USC(\widetilde{\Omega })$ be a viscosity subsolution and $v\in LSC(%
\widetilde{\Omega })$ a viscosity supersolution of \eqref{huvudekvationen}
in $\Omega ^{\circ }$. Then $\sup_{\widetilde{\Omega}}u-v\leq \sup_{\partial
\Omega \cup \overline{\Omega }_{0}}\left(u-v\right)^{+}$.
\end{lemma}

\noindent \textbf{Proof.} We may assume, by replacing $T>0$ by a smaller
number if necessary, that $u$ and $-v$ are bounded from above on $\widetilde{%
\Omega}$. We can also assume that $\sup_{\widetilde{\Omega}}u-v$ is attained
by using the well known fact that if $u$ is a subsolution of %
\eqref{huvudekvationen}, then so is
\begin{equation*}
u_{\beta }\left(t,x\right)=u\left(t,x\right)-\frac{\beta }{T-t},
\end{equation*}%
for all $\beta >0$. Assume that $\sup_{\widetilde{\Omega }%
}u-v=u\left(s,z\right)-v\left(s,z\right)>u\left(t,x\right)-v\left(t,x\right)$ for some $\left(s,z\right)\in \Omega ^{\circ }$ and
for all $\left(t,x\right)\in \partial \Omega \cup \overline{\Omega }_{0}$. As in
Section 5.B in \cite{CrandallIshiiLions1992}, we use the fact that if $u$ is
a viscosity subsolution, then so is $\bar{u}=u-K$ for every constant $K>0$.
Choose $K>0$ such that $\bar{u}\left(t,x\right)-v\left(t,x\right)\leq 0$ for all $\left(t,x\right)\in
\partial \Omega \cup \overline{\Omega }_{0}$ and such that $\bar{u}%
\left(s,z\right)-v\left(s,z\right):=\delta >0$. Using Lemma \ref{le:timdep_max} in place of
Theorem 8.3 in \cite{CrandallIshiiLions1992} and by observing that
assumptions \eqref{ass_F_cont}-\eqref{ass_F_XY} imply (assuming $\lambda >0$
as is possible by \eqref{assume_lambda_positive}) the corresponding
assumptions in \cite{CrandallIshiiLions1992}, we see that we can proceed as
in the proof of Theorem 8.2 in \cite{CrandallIshiiLions1992} to complete the
proof by deriving a contradiction. \hfill $\Box $

\vspace{0.2cm}

\noindent \textbf{Proof of Theorem \ref{comparison}. }In the following we
may assume, by replacing $T>0$ by a smaller number if necessary, that $u$
and $-v$ in Theorem \ref{comparison} are bounded from above on $\widetilde{%
\Omega }$. We will now produce approximations of $u$ and $v$ which allow us
to deal only with the inequalities involving $F$ and not the boundary
conditions. To construct these approximating functions, we note that Lemma %
\ref{testlemma5} applies with $\gamma $ replaced by $\widetilde{\gamma }$ as
well. Thus, there exists a $\mathcal{C}^{1,2}$ function $\alpha $ defined on
an open neighborhood of $\widetilde{\Omega }$ with the property that $\alpha
\geq 0$ on $\widetilde{\Omega }$ and $\left\langle D_{x}\alpha \left(t,x\right),%
\widetilde{\gamma }\left(t,x\right)\right\rangle \geq 1$ for $x\in \partial \Omega _{t}$%
, $t\in [0,T]$. For $\beta _{1}>0$, $\beta _{2}>0$ and $\beta _{3}>0$ we
define, for $\left(t,x\right)\in \widetilde{\Omega}$,
\begin{align}
u_{\beta _{1},\beta _{2},\beta _{3}}\left(t,x\right)& =u\left(t,x\right)-\beta _{1}\alpha
\left(t,x\right)-\beta _{2}-\frac{\beta _{3}}{T-t},  \notag  \label{approxdef} \\
v_{\beta _{1},\beta _{2}}\left(t,x\right)& =v\left(t,x\right)+\beta _{1}\alpha \left(t,x\right)+\beta _{2}.
\end{align}%
Given $\beta _{3},\beta _{2}>0$ there is $\beta _{1}=\beta _{1}\left(\beta
_{2}\right)\in \left(0,\beta _{2}\right)$ for which $u_{\beta _{1},\beta _{2},\beta _{3}}$
and $v_{\beta _{1},\beta _{2}}$ are sub- and supersolutions of %
\eqref{huvudekvationen}-\eqref{randvillkor}, with $f\left(t,x,r\right)$ replaced by $%
f\left(t,x,r\right)+\beta _{1}$ and $f\left(t,x,r\right)-\beta _{1}$, respectively. Indeed, if $%
\left(a,p,X\right)\in \mathcal{P}_{\widetilde{\Omega }}^{2,+}u_{\beta _{1},\beta
_{2},\beta _{3}}\left(t,x\right)$, then
\begin{equation}
\left( a+\beta _{1}\alpha _{t}\left(t,x\right)+\frac{\beta _{3}}{\left(T-t\right)^{2}},p+\beta
_{1}D\alpha \left(t,x\right),X+\beta _{1}D^{2}\alpha \left(t,x\right)\right) \in \mathcal{P}_{%
\widetilde{\Omega }}^{2,+}u\left(t,x\right).  \label{eq:punkt_i_subjet}
\end{equation}%
Hence, if $u$ satisfies \eqref{randvillkor}, then $\langle p+\beta
_{1}D\alpha \left(t,x\right),\widetilde{\gamma }\left(t,x\right)\rangle +f\left(t,x,u\left(t,x\right)\right)\leq 0$ and
since $\langle D\alpha \left(t,x\right),\widetilde{\gamma }\left(t,x\right)\rangle \geq 1$, $%
u_{\beta _{1},\beta _{2},\beta _{3}}\leq u$ and by %
\eqref{ass_f_nondecreasing} we obtain
\begin{equation}
\langle p,\widetilde{\gamma }\left(t,x\right)\rangle +f\left(t,x,u_{\beta _{1},\beta
_{2},\beta _{3}}\right)+\beta _{1}\leq 0.  \label{eq:RV_uppfylld_approx}
\end{equation}%
Using \eqref{eq:punkt_i_subjet} we also see that if $u$ satisfies %
\eqref{huvudekvationen} then
\begin{equation*}
a+\beta _{1}\alpha _{t}\left(t,x\right)+\frac{\beta _{3}}{\left(T-t\right)^{2}}+F\left(t,x,u,p+\beta
_{1}D\alpha \left(t,x\right),X+\beta _{1}D^{2}\alpha \left(t,x\right)\right)\leq 0.
\end{equation*}%
Using \eqref{ass_F_nondecreasing} and \eqref{ass_F_boundary}, assuming also
that the support of $\alpha $ lies within $U$, we have
\begin{align}
a+\beta _{1}\alpha _{t}\left(t,x\right)+F\left(t,x,u_{\beta _{1},\beta _{2},\beta
_{3}},p,X\right)+\lambda \beta _{2}&  \label{eq:EQ_uppfylld_approx} \\
-m_{2}\left( |\beta _{1}D\alpha \left(t,x\right)|+||\beta _{1}D^{2}\alpha
\left(t,x\right)||\right) & \leq 0.  \notag
\end{align}%
From \eqref{eq:RV_uppfylld_approx} and \eqref{eq:EQ_uppfylld_approx} it
follows that, given $\beta _{2},\beta _{3}>0$, there exist $\beta _{1}\in
\left(0,\beta _{2}\right)$ such that $u_{\beta _{1},\beta _{2},\beta _{3}}$ is a
subsolution of \eqref{huvudekvationen}-\eqref{randvillkor} with $f\left(t,x,u\right)$
replaced by $f\left(t,x,u\right)+\beta _{1}$. The fact that $v_{\beta _{1},\beta _{2}}$
is a supersolution follows by a similar calculation.

To complete the proof of the comparison principle, it is sufficient to prove
that%
\begin{equation*}
\max_{\widetilde{\Omega }}\left(u_{\beta _{1},\beta _{2},\beta _{3}}-v_{\beta
_{1},\beta _{2}}\right)\leq 0,
\end{equation*}%
holds for all $\beta _{2}>0$ and $\beta _{3}>0$. Assume that
\begin{equation*}
\sigma =\max_{\widetilde{\Omega }}\left(u_{\beta _{1},\beta _{2},\beta
_{3}}-v_{\beta _{1},\beta _{2}}\right)>0.
\end{equation*}%
We will derive a contradiction for any $\beta _{3}$ if $\beta _{2}$ (and
hence $\beta _{1}$) is small enough. To simplify notation, we write, in the
following, $u,v$ in place of $u_{\beta _{1},\beta _{2},\beta _{3}},v_{\beta
_{1},\beta _{2}}$. By Lemma \ref{maxrand}, $u\left(0,\cdot \right)\leq v\left(0,\cdot \right)$,
upper semicontinuity of $u-v$ and boundedness from above of $u-v$, we
conclude that for any $\beta _{3}>0$
\begin{equation}
\sigma =\left(u-v\right)\left(s,z\right),\quad \text{for some }z\in \partial \Omega _{s}\text{ and
}s\in \left( 0,T\right) .  \label{sigma}
\end{equation}%
Let $\widetilde{B}\left(\left(s,z\right),\delta \right)=\{\left(t,x\right):\left\vert \left(t,x\right)-\left(s,z\right)\right\vert
\leq \delta \}$ and define
\begin{equation*}
E:=\widetilde{B}\left(\left(s,z\right),\delta \right)\cap \widetilde{\Omega }.
\end{equation*}%
By Remark \ref{spaceremark}, there exists $\theta \in \left(0,1\right)$ such that
\begin{equation}
\left\langle x-y,\widetilde{\gamma }\left( t,x\right) \right\rangle \geq
-\theta \left\vert x-y\right\vert ,\quad \text{for all }\left( t,x\right)
\in E\setminus \Omega ^{\circ }\text{ and }\left( t,y\right) \in E.
\label{make_use_of_cone}
\end{equation}%
By decreasing $\delta $ if necessary, we may assume that %
\eqref{ass_F_boundary} holds in $E$. From now on, we restrict our attention
to events in the set $E$. By Lemma \ref{testlemma4} we obtain, for any $%
\theta \in \left(0,1\right)$, a family $\left\{ w_{\varepsilon }\right\} _{\varepsilon
>0}$ of functions $w_{\varepsilon }\in \mathcal{C}^{1,2}\left( \left[ 0,T%
\right] \times
\mathbb{R}
^{n}\times
\mathbb{R}
^{n},%
\mathbb{R}
\right) $ and positive constants $\chi ,C$ (independent of $\varepsilon $)
such that \eqref{testlemma41}, \eqref{testlemma42}, \eqref{testlemma45}-%
\eqref{testlemma47} as well as
\begin{equation}
\left\langle D_{x}w_{\varepsilon }\left( t,x,y\right) ,\widetilde{\gamma }%
\left( t,x\right) \right\rangle \geq -C\frac{\left\vert x-y\right\vert ^{2}}{%
\varepsilon },\quad \text{if\quad }\left\langle x-y,\widetilde{\gamma }%
\left( t,x\right) \right\rangle \geq -\theta \left\vert x-y\right\vert ,
\label{test3}
\end{equation}%
\begin{equation}
\left\langle D_{y}w_{\varepsilon }\left( t,x,y\right) ,\widetilde{\gamma }%
\left( t,y\right) \right\rangle \geq -C\frac{\left\vert x-y\right\vert ^{2}}{%
\varepsilon },\quad \text{if\quad }\left\langle y-x,\widetilde{\gamma }%
\left( t,y\right) \right\rangle \geq -\theta \left\vert x-y\right\vert ,
\label{test4}
\end{equation}%
hold. Note that \eqref{test3} and \eqref{test4} are direct analogues to %
\eqref{testlemma43} and \eqref{testlemma44} but with $\gamma $ replaced by $%
\widetilde{\gamma }$.

Let $\varepsilon >0$ be given and define

\begin{equation*}
\Phi \left(t,x,y\right)=u\left(t,x\right)-v\left(t,y\right)-\varphi \left(t,x,y\right),
\end{equation*}%
where
\begin{equation*}
\varphi \left(t,x,y\right)=w_{\varepsilon }\left(t,x,y\right)+f\left(s,z,u\left(s,z\right)\right)\langle y-x,\widetilde{%
\gamma }\left(s,z\right)\rangle +\beta _{1}|x-z|^{2}+\left(t-s\right)^{2}.
\end{equation*}%
Let $\left(t_{\varepsilon },x_{\varepsilon },y_{\varepsilon }\right)$ be a maximum
point of $\Phi $. From \eqref{testlemma41} and \eqref{testlemma42} we have
\begin{align}
\sigma & -C\varepsilon \leq \Phi \left(s,z,z\right)\leq \Phi \left(t_{\varepsilon
},x_{\varepsilon },y_{\varepsilon }\right)\leq u\left(t_{\varepsilon },x_{\varepsilon
}\right)-v\left(t_{\varepsilon },y_{\varepsilon }\right)-\chi \frac{\left\vert x_{\varepsilon
}-y_{\varepsilon }\right\vert ^{2}}{\varepsilon }
\label{eq:maxet_gor_att_allt_konvergerar} \\
& -f\left(s,z,u\left(s,z\right)\right)\langle y_{\varepsilon }-x_{\varepsilon },\widetilde{\gamma }%
\left(s,z\right)\rangle -\beta _{1}|x_{\varepsilon }-z|^{2}-\left(t_{\varepsilon }-s\right)^{2}.
\notag
\end{align}%
From this we first see that
\begin{equation*}
|x_{\varepsilon }-y_{\varepsilon }|\rightarrow 0\qquad \text{as}\qquad
\varepsilon \rightarrow 0.
\end{equation*}%
Therefore, using the upper semi-continuity of $u-v$ and %
\eqref{eq:maxet_gor_att_allt_konvergerar} we also obtain
\begin{align}
\frac{|x_{\varepsilon }-y_{\varepsilon }|^{2}}{\varepsilon }& \rightarrow
0,\qquad x_{\varepsilon },y_{\varepsilon }\rightarrow z,\qquad
t_{\varepsilon }\rightarrow s,  \notag  \label{as_ep_to_0} \\
u\left(t_{\varepsilon },x_{\varepsilon }\right)& \rightarrow u\left(s,z\right),\qquad
v\left(t_{\varepsilon },y_{\varepsilon }\right)\rightarrow v\left(s,z\right),
\end{align}%
as $\varepsilon \rightarrow 0$. In the following we assume $\varepsilon $ to
be so small that $\left(t_{\varepsilon },x_{\varepsilon }\right)\in E$

We introduce the notation
\begin{align*}
\bar{p}& =D_{x}\varphi \left(t_{\varepsilon },x_{\varepsilon },y_{\varepsilon
}\right)=D_{x}w_{\varepsilon }\left(t_{\varepsilon },x_{\varepsilon },y_{\varepsilon
}\right)-f\left(s,z,u\left(s,z\right)\right)\widetilde{\gamma }\left(s,z\right)+2\beta _{1}\left( x_{\varepsilon
}-z\right) , \\
\bar{q}& =D_{y}\varphi \left(t_{\varepsilon },x_{\varepsilon },y_{\varepsilon
}\right)=D_{y}w_{\varepsilon }\left(t_{\varepsilon },x_{\varepsilon },y_{\varepsilon
}\right)+f\left(s,z,u\left(s,z\right)\right)\widetilde{\gamma }\left(s,z\right),
\end{align*}%
and observe that
\begin{align}
& \langle \bar{p},\widetilde{\gamma }\left(t_{\varepsilon },x_{\varepsilon
}\right)\rangle +f\left(t_{\varepsilon },x_{\varepsilon },u\left(t_{\varepsilon
},x_{\varepsilon }\right)\right)  \notag \\
=& \langle D_{x}w_{\varepsilon }\left(t_{\varepsilon },x_{\varepsilon
},y_{\varepsilon }\right),\widetilde{\gamma }\left(t_{\varepsilon },x_{\varepsilon
}\right),\rangle +f\left(t_{\varepsilon },x_{\varepsilon },u\left(t_{\varepsilon
},x_{\varepsilon }\right)\right)  \notag \\
& -f\left(s,z,u\left(s,z\right)\right)\langle \widetilde{\gamma }\left(s,z\right),\widetilde{\gamma }%
\left(t_{\varepsilon },x_{\varepsilon }\right)\rangle +2\beta _{1}\langle
x_{\varepsilon }-z,\widetilde{\gamma }\left(t_{\varepsilon },x_{\varepsilon
}\right)\rangle ,  \label{peq}
\end{align}%
and
\begin{eqnarray}
&&-\langle \bar{q},\widetilde{\gamma }\left(t_{\varepsilon },y_{\varepsilon
}\right)\rangle +f\left(t_{\varepsilon },y_{\varepsilon },v\left(t_{\varepsilon
},x_{\varepsilon }\right)\right)  \notag \\
&=&-\langle D_{y}w_{\varepsilon }\left(t_{\varepsilon },x_{\varepsilon
},y_{\varepsilon }\right),\widetilde{\gamma }\left(t_{\varepsilon },y_{\varepsilon
}\right)\rangle +f\left(t_{\varepsilon },y_{\varepsilon },v\left(t_{\varepsilon
},y_{\varepsilon }\right)\right)  \notag \\
&&-f\left(s,z,u\left(s,z\right)\right)\left\langle \widetilde{\gamma }\left(s,z\right) ,\widetilde{%
\gamma }\left(t_{\varepsilon },y_{\varepsilon }\right)\right\rangle.  \label{qeq}
\end{eqnarray}
Using \eqref{smooth_gamma}, \eqref{f_kontinuerlig}, %
\eqref{ass_f_nondecreasing} and \eqref{as_ep_to_0}-\eqref{qeq} we see that
if $\varepsilon $ is small enough, then
\begin{align}
& \langle D_{x}w_{\varepsilon }\left(t_{\varepsilon },x_{\varepsilon
},y_{\varepsilon }\right),\widetilde{\gamma }\left(t_{\varepsilon },x_{\varepsilon
}\right)\rangle \geq -\frac{\beta _{1}}{2}  \notag
\label{boundary_cond_elimination_2} \\
& \implies \langle \bar{p},\widetilde{\gamma }\left(t_{\varepsilon
},x_{\varepsilon }\right)\rangle +f\left(t_{\varepsilon },x_{\varepsilon
},u\left(t_{\varepsilon },x_{\varepsilon }\right)\right)+\beta _{1}>0,  \notag \\
& \langle D_{y}w_{\varepsilon }\left(t_{\varepsilon },x_{\varepsilon
},y_{\varepsilon }\right),\widetilde{\gamma }\left(t_{\varepsilon },y_{\varepsilon
}\right)\rangle \geq -\frac{\beta _{1}}{2}  \notag \\
& \implies -\langle \bar{q},\widetilde{\gamma }\left(t_{\varepsilon
},y_{\varepsilon }\right)\rangle +f\left(t_{\varepsilon },y_{\varepsilon
},v\left(t_{\varepsilon },x_{\varepsilon }\right)\right)-\beta _{1}<0.
\end{align}%
Moreover, from \eqref{make_use_of_cone}-\eqref{test4}, we also have
\begin{align}
& \langle D_{x}w_{\varepsilon }\left(t_{\varepsilon },x_{\varepsilon
},y_{\varepsilon }\right),\widetilde{\gamma }\left(t_{\varepsilon },x_{\varepsilon
}\right)\rangle \geq -C\frac{|x_{\varepsilon }-y_{\varepsilon }|^{2}}{\varepsilon }%
,\quad \text{if }x_{\varepsilon }\in \partial \Omega _{t_{\varepsilon }},
\notag  \label{boundary_cond_elimination_3} \\
& \langle D_{y}w_{\varepsilon }\left(t_{\varepsilon },x_{\varepsilon
},y_{\varepsilon }\right),\widetilde{\gamma }\left(t_{\varepsilon },y_{\varepsilon
}\right)\rangle \geq -C\frac{|x_{\varepsilon }-y_{\varepsilon }|^{2}}{\varepsilon }%
,\quad \text{if }y_{\varepsilon }\in \partial \Omega _{t_{\varepsilon }}.
\end{align}%
Using \eqref{boundary_cond_elimination_2} and %
\eqref{boundary_cond_elimination_3}, it follows by the definition of
viscosity solutions that if $\varepsilon $ is small enough, say $%
0<\varepsilon <\varepsilon _{\beta _{1}}$, then
\begin{equation}
a+F\left(t_{\varepsilon },x_{\varepsilon },u\left(t_{\varepsilon },x_{\varepsilon }\right),%
\bar{p},X\right)\leq 0\leq -b+F\left(t_{\varepsilon },y_{\varepsilon },v\left(t_{\varepsilon
},y_{\varepsilon }\right),-\bar{q},-Y\right),  \label{flyttauppskiten}
\end{equation}%
whenever
\begin{equation*}
\left( a,\bar{p},X\right) \in \overline{\mathcal{P}}_{\widetilde{\Omega }%
}^{2,+}u\left(t_{\varepsilon },x_{\varepsilon }\right)\quad \text{and}\quad \left( -b,-%
\bar{q},-Y\right) \in \overline{\mathcal{P}}_{\widetilde{\Omega }%
}^{2,-}v\left(t_{\varepsilon },y_{\varepsilon }\right).
\end{equation*}%
We next intend to use Lemma \ref{le:timdep_max} to show the existence of
such matrices $X$, $Y$ and numbers $a,b$. Hence, we have to verify condition %
\eqref{eq:besvarlig_assumption}. To do so, we observe that %
\eqref{boundary_cond_elimination_2} holds true with $\bar{p}$ and $\bar{q}$
replaced by any $p$ and $q$ satisfying $|\bar{p}-p|\leq r$ and $|\bar{q}%
-q|\leq r$ if we choose $r=r\left(\varepsilon \right)$ small enough. It follows that
also \eqref{flyttauppskiten} holds with these $p$ and $q$ and we can
conclude
\begin{equation*}
a\leq -F\left(t_{\varepsilon },x_{\varepsilon },u\left(t_{\varepsilon },x_{\varepsilon
}\right),p,X\right)\leq C\quad \text{and}\quad b\leq F\left(t_{\varepsilon },y_{\varepsilon
},v\left(t_{\varepsilon },y_{\varepsilon }\right),-q,-Y\right)\leq C,
\end{equation*}%
for some $C=C\left(\varepsilon \right)$ whenever $\left(a,p,X\right)$ and $\left( b,q,Y\right) $
is as in \eqref{eq:besvarlig_assumption}. Hence, condition %
\eqref{eq:besvarlig_assumption} holds and Lemma \ref{le:timdep_max} gives
the existence of $X,Y\in \mathbb{S}^{n}$ and $a,b\in \mathbb{R}$ such that
\begin{align}
& -\left( \frac{1}{\varepsilon }+||A||\right) I\leq \left(
\begin{array}{cc}
X & 0 \\
0 & Y%
\end{array}%
\right) \leq A+\varepsilon A^{2},  \notag  \label{eq:result_from_CIL92} \\
& \left(a,\bar{p},X\right)\in \overline{\mathcal{P}}_{\widetilde{\Omega }%
}^{2,+}u\left(t_{\varepsilon },x_{\varepsilon }\right),\quad \left(-b,-\bar{q},-Y\right)\in
\overline{\mathcal{P}}_{\widetilde{\Omega }}^{2,-}v\left(t_{\varepsilon
},y_{\varepsilon }\right),  \notag \\
& a+b=D_{t}\varphi \left( t_{\varepsilon },x_{\varepsilon },y_{\varepsilon
}\right) =D_{t}w_{\varepsilon }\left( t_{\varepsilon },x_{\varepsilon
},y_{\varepsilon }\right) +2\left(t_{\varepsilon }-s\right),
\end{align}%
where $A=D_{x,y}^{2}\left( w_{\varepsilon }\left(t_{\varepsilon },x_{\varepsilon
},y_{\varepsilon }\right)+\beta _{1}|x_{\varepsilon }-z|^{2}\right) $. Using %
\eqref{ass_F_nondecreasing}, \eqref{testlemma45} and \eqref{flyttauppskiten}
we obtain, by recalling that we can assume $\lambda >0$ in %
\eqref{ass_F_nondecreasing}, that%
\begin{eqnarray*}
0 &\geq &D_{t}w_{\varepsilon }\left( t_{\varepsilon },x_{\varepsilon
},y_{\varepsilon }\right) +2\left(t_{\varepsilon }-s\right) \\
&&+F\left(t_{\varepsilon },x_{\varepsilon },u\left(t_{\varepsilon },x_{\varepsilon }\right),%
\bar{p},X\right)-F\left(t_{\varepsilon },y_{\varepsilon },v\left(t_{\varepsilon
},y_{\varepsilon }\right),-\bar{q},-Y\right) \\
&\geq &-C\frac{|x_{\varepsilon }-y_{\varepsilon }|^{2}}{\varepsilon }%
+2\left(t_{\varepsilon }-s\right)+\lambda \left(u\left(t_{\varepsilon },x_{\varepsilon
}\right)-v\left(t_{\varepsilon },y_{\varepsilon }\right)\right) \\
&&+F\left(t_{\varepsilon },x_{\varepsilon },u\left(t_{\varepsilon },x_{\varepsilon }\right),%
\bar{p},X\right)-F\left(t_{\varepsilon },y_{\varepsilon },u\left(t_{\varepsilon
},x_{\varepsilon }\right),-\bar{q},-Y\right).
\end{eqnarray*}%
Next, assumption \eqref{ass_F_boundary} gives%
\begin{eqnarray}
0 &\geq &-C\bar{s}+2\left(t_{\varepsilon }-s\right)+\lambda \left(u\left(t_{\varepsilon
},x_{\varepsilon }\right)-v\left(t_{\varepsilon },y_{\varepsilon }\right)\right)  \notag \\
&&+F\left(t_{\varepsilon },x_{\varepsilon },u\left(t_{\varepsilon },x_{\varepsilon }\right),-%
\bar{q},X-C\bar{s}I\right)-F\left(t_{\varepsilon },y_{\varepsilon },u\left(t_{\varepsilon
},x_{\varepsilon }\right),-\bar{q},-Y+C\bar{s}I\right)  \notag \\
&&-m_{2}\left(|\bar{p}+\bar{q}|+C\bar{s}\right)-m_{2}\left(C\bar{s}\right),
\label{sista_med_extrasteg}
\end{eqnarray}%
where we use the notation $\bar{s}=|x_{\varepsilon }-y_{\varepsilon
}|^{2}/\varepsilon $. Note that since the eigenvalues of $\varepsilon A^{2}$
are given by $\varepsilon \lambda ^{2}$, where $\lambda $ is an eigenvalue
to $A$, and since $\lambda $ is bounded, $A+\varepsilon A^{2}\leq CA$.
Hence, by \eqref{testlemma47} we obtain
\begin{equation*}
A+\varepsilon A^{2}\leq \frac{C}{\varepsilon }\left(
\begin{array}{cc}
I & -I \\
-I & I%
\end{array}%
\right) +C\bar{s}I_{2n},
\end{equation*}%
and since $||A||\leq C/\varepsilon $ for some large $C$, we also conclude
that \eqref{eq:result_from_CIL92} implies
\begin{equation*}
-\frac{C}{\varepsilon }I_{2n}\leq \left(
\begin{array}{cc}
X-C\bar{s}I & 0 \\
0 & Y-C\bar{s}I%
\end{array}%
\right) \leq \frac{C}{\varepsilon }\left(
\begin{array}{cc}
I & -I \\
-I & I%
\end{array}%
\right) .
\end{equation*}%
Using the above inequality, assumption \eqref{ass_F_XY}, %
\eqref{sista_med_extrasteg}, the definition of $\bar{q}$ and %
\eqref{testlemma46} we have%
\begin{eqnarray*}
0 &\geq &-C\bar{s}+2\left(t_{\varepsilon }-s\right)+\lambda \left(u\left(t_{\varepsilon
},x_{\varepsilon }\right)-v\left(t_{\varepsilon },y_{\varepsilon }\right)\right) \\
&&-m_{1}\left(C|x_{\varepsilon }-y_{\varepsilon }|+2C\bar{s}\right)-m_{2}\left(|\bar{p}+\bar{%
q}|+C\bar{s}\right)-m_{2}\left(C\bar{s}\right),
\end{eqnarray*}%
when $0<\varepsilon <\varepsilon _{\beta _{1}}$ and $u\left(t_{\varepsilon
},x_{\varepsilon }\right)\geq v\left(t_{\varepsilon },y_{\varepsilon }\right)$. Sending first
$\varepsilon $ and then $\beta _{2}$ to zero (the latter implies $\beta
_{1}\rightarrow 0$) and using \eqref{testlemma46} we obtain a contradiction.
This completes the proof of the comparison principle in Theorem \ref%
{comparison}.\hfill $\Box $

\vspace{0.2cm}

Using the same methodology as in the proof of Theorem \ref{comparison}, we
are now able to prove the comparison principle for mixed boundary conditions
stated in Corollary \ref{maxrand_partial}. This result will be an important
ingredient in the proof of Theorem \ref{existence}.

\vspace{0.2cm}

\noindent \textbf{Proof of Corollary \ref{maxrand_partial}.} If $u$ is a
viscosity subsolution, then so is $u-K$ for all $K>0$. It thus suffices to
prove that if $u\leq v$ on $\left(\partial \Omega \setminus G\right)\cup \overline{%
\Omega }_{0}$, then $u\leq v$ in $\widetilde{\Omega }$. If $G=\partial
\Omega $, then this implication and its proof is identical to Theorem \ref%
{comparison}. If $G\subset \partial \Omega $ is arbitrary, then we know by
assumption that $u\leq v$ on $\partial \Omega \setminus G$ and so the point $%
\left(s,z\right)$ defined in \eqref{sigma} must belong to the set $G$ where the
boundary condition is satisfied. Hence, we can follow the proof of Theorem %
\ref{comparison} and conclude that $u\leq v$ in $\widetilde{\Omega }$. $%
\hfill \Box $

\vspace{0.2cm}

\noindent \textbf{Proof of Theorem \ref{existence}. }We will prove existence
using Perron's method. In particular, we show that the supremum of all
subsolutions to the initial value problem given by %
\eqref{initial_value_problem} is indeed a solution to the same problem. To
ensure that the supremum is taken over a nonempty set, we need to find at
least one subsolution to the problem. We also need to know that the supremum
is finite. This is obtained by producing a supersolution, which, due to the
comparison principle, provides an upper bound for the supremum.

To find the supersolution, let, for some constants $A$ and $B$ to be chosen
later,
\begin{equation*}
\widehat{v}=A\alpha \left(t,x\right)+B,\quad \text{for}\;\left(t,x\right) \in \widetilde{\Omega},
\end{equation*}
where $\alpha \left(t,x\right)$ is the function guaranteed by Lemma \ref{testlemma5}.
By \eqref{f_kontinuerlig}, \eqref{ass_f_nondecreasing} and the boundedness
of $\Omega ^{\circ }$, we can find $A>0$ such that
\begin{equation*}
\langle D\widehat{v}\left(t,x\right),\widetilde{\gamma }\left(t,x\right)\rangle +f\left(t,x,\widehat{v}%
\left(t,x\right)\right)\geq A+f\left(t,x,0\right)\geq 0,
\end{equation*}%
for $\left(t,x\right)\in \partial \Omega $. Moreover, since the support of $\alpha $
lies in $U$, we have, with $\lambda $ and $m_{2}$ defined in %
\eqref{ass_F_nondecreasing} and \eqref{ass_F_boundary},%
\begin{eqnarray*}
&&D_{t}\widehat{v}\left(t,x\right)+F\left(t,x,\widehat{v}\left(t,x\right),D\widehat{v}\left(t,x\right),D^{2}%
\widehat{v}\left(t,x\right)\right) \\
&\geq &-A\sup_{U}\{|D_{t}\alpha \left(t,x\right)|\}+B\lambda +F\left(t,x,0,0,0\right) \\
&&-\sup_{U}m_{2}\left( A\left( |D\alpha \left(t,x\right)|+||D^{2}\alpha \left(t,x\right)||\right)
\right) .
\end{eqnarray*}%
By \eqref{ass_F_cont}, the boundedness of $\Omega ^{\circ }$ and by
recalling that we can assume $\lambda >0$, we see that taking $B$ large
enough, $\widehat{v}$ is a classical supersolution of %
\eqref{initial_value_problem}. Hence, using \eqref{F_fundamental} and
Proposition 7.2 in \cite{CrandallIshiiLions1992}, $\widehat{v}$ is also a
viscosity supersolution. Next, we observe that $\check{u}=-\widehat{v}$ is a
viscosity subsolution $\check{u}$ to the problem given by %
\eqref{initial_value_problem}.

We now apply Perron's method by defining our solution candidate as
\begin{equation*}
\widetilde{w}:=\sup \{w\left(x\right):\text{$w\in USC(\widetilde{\Omega })$ is a
viscosity subsolution of \eqref{initial_value_problem}}\}.
\end{equation*}%
In the following we let $u^{\ast }$ and $u_{\ast }$ denote the upper and
lower semicontinuous envelopes of a function $u$, respectively. By the
comparison principle and by construction we obtain
\begin{equation}
\check{u}_{\ast }\leq \widetilde{w}_{\ast }\leq \widetilde{w}^{\ast }\leq
\widehat{v}^{\ast }\quad \text{on}\;\widetilde{\Omega }.
\label{eq:by_construction}
\end{equation}%
Let us assume for the moment that $\widetilde{w}^{\ast }$ satisfies the
initial condition of being a subsolution and that $\widetilde{w}_{\ast }$
satisfies the initial condition of being a supersolution, that is
\begin{equation}
\widetilde{w}^{\ast }\left(0,x\right)\leq g\left(x\right)\leq \widetilde{w}_{\ast }\left(0,x\right),\quad
\text{for all }x\in \overline{\Omega }_{0}.  \label{eq:initial_assumed}
\end{equation}%
We can then proceed as in \cite{CrandallIshiiLions1992} (see also \cite%
{Barles1993} and \cite{Ishii1987}) to show that $\widetilde{w}^{\ast }$ is a
viscosity subsolution and $\widetilde{w}_{\ast }$ is a viscosity
supersolution of the initial value problem in \eqref{initial_value_problem}.
Using the comparison principle again, we then have $\widetilde{w}_{\ast
}\geq \widetilde{w}^{\ast }$ and so by \eqref{eq:by_construction} $%
\widetilde{w}_{\ast }=\widetilde{w}^{\ast }$ is the requested viscosity
solution. To complete the proof of Theorem \ref{existence}, it hence
suffices to prove \eqref{eq:initial_assumed}. This will be achieved by
constructing families of explicit viscosity sub- and supersolutions.

We first show that the subsolution candidate $\widetilde{w}^{\ast }$
satisfies the initial conditions for all $x\in \Omega _{0}$. To this end, we
define, for arbitrary $z\in \Omega _{0}$ and $\varepsilon >0$, the barrier
function
\begin{equation*}
V_{z,\varepsilon }\left(t,x\right)=g\left(z\right)+\varepsilon +B|x-z|^{2}+Ct,\quad \text{for}%
\;\left(t,x\right)\in \left[ 0,T\right] \times \mathbb{R}^{n},
\end{equation*}%
where $B$ and $C$ are constants, which may depend on $z$ and $\varepsilon $,
to be chosen later. We first observe that, by continuity of $g$ and
boundedness of $\Omega _{0}$, we can, for any $\varepsilon >0$, choose $B$
so large that $V_{z,\varepsilon }\left(0,x\right)\geq g\left(x\right)$, for all $x\in \overline{%
\Omega }_{0}$. Moreover, since $\widetilde{w}$ is bounded on $\overline{%
\Omega }$, we conclude, by increasing $B$ and $C$ if necessary, that we also
have
\begin{equation*}
V_{z,\varepsilon }\left(t,x\right)\geq \widetilde{w}\left(t,x\right),\quad \text{for}\;\left(t,x\right)\in
\partial \Omega \cup \overline{\Omega }_{0}.
\end{equation*}%
A computation shows that, for $z$, $\varepsilon $, $B$ given, we can choose
the constant $C$ so large that $V_{z,\epsilon }$ is a classical
supersolution of \eqref{huvudekvationen} in $[0,\infty )\times \mathbb{R}%
^{n} $. Hence, by \eqref{F_fundamental}, $V_{z,\epsilon }$ is also a
continuous viscosity supersolution of \eqref{huvudekvationen} in $\Omega
^{\circ }$. By the maximum principle in Lemma \ref{maxrand} applied to $%
V_{z,\varepsilon }$ and each component in the definition of $\widetilde{w}$,
we obtain
\begin{equation}
V_{z,\varepsilon }\left(t,x\right)\geq \widetilde{w}\left(t,x\right),\quad \text{for}\;\left(t,x\right)\in
\widetilde{\Omega }.  \label{eq:first_barrier_downpush_inside}
\end{equation}%
It follows that $\widetilde{w}^{\ast }\leq V_{z,\varepsilon }^{\ast
}=V_{z,\varepsilon }$ in this set and hence the initial condition in $\Omega
_{0}$ follows since for any $x\in \Omega _{0}$
\begin{equation}
\widetilde{w}^{\ast }\left(0,x\right)\leq \inf_{\varepsilon ,z}V_{z,\varepsilon
}\left(0,x\right)=g\left(x\right).  \label{eq:first_barrier_downpush_inside_final}
\end{equation}%
To prove that the supersolution candidate $\widetilde{w}_{\ast }$ satisfies
the initial condition in $\Omega _{0}$, we proceed similarly by studying a
family of subsolutions of the form
\begin{equation*}
U_{z,\epsilon }\left(t,x\right)=g\left(z\right)-B|x-z|^{2}-\varepsilon -Ct.
\end{equation*}

We next prove that $\widetilde{w}^{\ast }$ satisfies the boundary conditions
for each $x\in \partial \Omega _{0}$. In this case the barriers above will
not work as we cannot ensure that they exceed $\widetilde{w}^{\ast }$ on $%
\partial \Omega $. Instead, we will construct barriers that are sub- and
supersolutions only locally, near the boundary, during a short time
interval. These local barriers are useful due to the maximum principle for
mixed boundary conditions proved in Corollary \ref{maxrand_partial}. To
construct the local barriers, fix $\widehat{z}\in \partial \Omega _{0}$ and
let $z\left(t\right)$ be the H\"{o}lder continuous function
\begin{equation*}
z\left(t\right)=\widehat{z}-K\widetilde{\gamma }\left(0,\widehat{z}\right)t^{\widehat{\alpha }},
\end{equation*}%
where $\widehat{\alpha }$ is the H\"{o}lder exponent from \eqref{tempholder}
and $K$ is a constant depending on the H\"{o}lder constant and the shape of
the exterior cones in \eqref{boundarylip}. It follows that $z\left(t\right)$ stays
inside of $\Omega $ for a short time and that $z\left(0\right)=\widehat{z}$. Consider,
for $\varepsilon >0$, the barrier function
\begin{equation*}
\widetilde{V}_{\varepsilon ,\widehat{z}}\left(t,x\right)=g\left(\widehat{z}\right)+A\left( \alpha
\left(t,x\right)-\alpha \left(0,\widehat{z}\right)\right) +e^{\left(\widehat{C}/\chi \right)\alpha
\left(t,x\right)}w_{\varepsilon }\left(t,x,z\left(t\right)\right)+B+Ct^{\widehat{\alpha }},
\end{equation*}%
whenever $\left(t,x\right)\in \lbrack 0,T]\times \mathbb{R}^{n}$, where $\widehat{C}$
and $\chi $ are the constants from Lemma \ref{testlemma4} and $A,B$ and $C$
are constants to be chosen later, possibly depending on $\widehat{z}$ and $%
\varepsilon $. We first show that for any choice of $A$, we can find $B$
such that
\begin{equation}
g\left(x\right)\leq \widetilde{V}_{\varepsilon ,\widehat{z}}\left(0,x\right),\quad \text{for all }%
x\in \overline{\Omega }_{0}\text{\quad and\quad }\inf_{\varepsilon }%
\widetilde{V}_{\varepsilon ,\widehat{z}}\left(0,\widehat{z}\right)=g\left(\widehat{z}\right).
\label{eq:sup_barrier_above_g}
\end{equation}%
Indeed, to prove the left inequality in \eqref{eq:sup_barrier_above_g},
observe that by \eqref{testlemma41} we have $\chi \left\vert x-\widehat{z}%
\right\vert ^{2}/\varepsilon \leq w_{\varepsilon }\left( 0,x,\widehat{z}%
\right) $. Moreover, by the continuity of $g\left( \cdot \right) -A\alpha
\left( 0,\cdot \right) $ in $\overline{\Omega }_{0}$, we can find $B$,
depending on $\varepsilon $ and $A$, so that
\begin{equation*}
g\left(x\right)-g\left(\widehat{z}\right)-A\left( \alpha \left(0,x\right)-\alpha \left(0,\widehat{z}\right)\right) \leq
B+\chi \frac{\left\vert x-\widehat{z}\right\vert ^{2}}{\varepsilon }.
\end{equation*}%
This proves the left inequality in \eqref{eq:sup_barrier_above_g}. Finally,
it is no restriction to assume that $B\rightarrow 0$ as $\varepsilon
\rightarrow 0$, and this implies the right inequality in %
\eqref{eq:sup_barrier_above_g}.

We next show that $\widetilde{V}_{\varepsilon ,\widehat{z}}$ satisfies the
boundary condition in a small neighborhood of $\widehat{z}$ in ${\partial
\Omega }$. To do so, let $E_{\widehat{z}}=\left(0,\kappa \right)\times B\left(\widehat{z}%
,\rho \right)$ for some $\kappa ,\rho >0$ to be chosen. We intend to find $\kappa $%
, $\rho $, $A$ and $C$ such that
\begin{equation}
\langle D_{x}\widetilde{V}_{\varepsilon ,\widehat{z}}\left(t,x\right),\widetilde{\gamma
}\left(t,x\right)\rangle +f\left(t,x,\widetilde{V}_{\varepsilon ,\widehat{z}}\left(t,x\right)\right)\geq
0,\quad \text{for}\;\left(t,x\right)\in E_{\widehat{z}}\cap \partial \Omega .
\label{eq:sup_barrier_RV}
\end{equation}%
First, observe that $\alpha $ is differentiable in time on $\overline{\Omega
}$. Therefore, by taking $C$ large enough and by using %
\eqref{eq:sup_barrier_above_g} we ensure that
\begin{equation*}
\widetilde{V}_{\varepsilon ,\widehat{z}}\left(t,x\right)\geq g\left(\widehat{z}\right),\quad \text{%
for}\;\left(t,x\right)\in \overline{\Omega }.
\end{equation*}%
In general, the choice of $C$ will depend on $A$, but it is evident from the
next inequality that this will not give rise to circular reasoning. By %
\eqref{ass_f_nondecreasing} and the boundedness of $\overline{\Omega }$, we
can choose $A$ so that
\begin{equation*}
f\left(t,x,\widetilde{V}_{\varepsilon ,\widehat{z}}\left(t,x\right)\right)\geq f\left(t,x,g\left(\widehat{z}%
\right)\right)\geq -A,\quad \text{for}\;\left(t,x\right)\in \overline{\Omega }.
\end{equation*}%
Thus, the boundary condition in \eqref{eq:sup_barrier_RV} will follow if we
can prove
\begin{equation}
\langle D_{x}\widetilde{V}_{\varepsilon ,\widehat{z}}\left(t,x\right),\widetilde{\gamma
}\left(t,x\right)\rangle \geq A,\quad \text{for}\;\left(t,x\right)\in E_{\widehat{z}}\cap \partial
\Omega .  \label{eq:RV_andra_halvan}
\end{equation}%
To this end, choose $\rho $ and $\kappa $ so small that
\begin{equation}
\left\langle x-z\left(t\right),\widetilde{\gamma }\left( t,x\right) \right\rangle \geq
-\theta \left\vert x-z\left(t\right)\right\vert \quad \text{whenever}\;x\in B\left(
\widehat{z},\rho \right) \cap \partial \Omega _{t},\;t\in \left[ 0,\kappa %
\right] .
\end{equation}%
Inequality \eqref{test3} then holds with $y=z\left( t\right) $ for all $%
\left(t,x\right)\in E_{\widehat{z}}\cap \partial \Omega $. Together with the properties
of $\alpha $, this gives
\begin{eqnarray*}
&&\langle D_{x}\widetilde{V}_{\varepsilon ,\widehat{z}}\left(t,x\right),\widetilde{%
\gamma }\left(t,x\right)\rangle \\
&=&A\langle D_{x}\alpha \left(t,x\right),\widetilde{\gamma }\left(t,x\right)\rangle +e^{\left(\widehat{C%
}/\chi \right)\alpha \left(t,x\right)} \\
&&\cdot \left\langle D_{x}w_{\varepsilon }\left( t,x,z\left(t\right)\right)
+w_{\varepsilon }\left(t,x,z\left(t\right)\right)\frac{\widehat{C}}{\chi }D_{x}\alpha \left(t,x\right),%
\widetilde{\gamma }\left(t,x\right)\right\rangle \\
&\geq &A-\widehat{C}\frac{|x-z\left(t\right)|^{2}}{\varepsilon }+\chi \frac{|x-z\left(t\right)|^{2}%
}{\varepsilon }\frac{\widehat{C}}{\chi }=A,\quad \text{for}\;\left(t,x\right)\in E_{%
\widehat{z}}\cap \partial \Omega .
\end{eqnarray*}%
This proves \eqref{eq:RV_andra_halvan} and hence the boundary condition %
\eqref{eq:sup_barrier_RV} follows.

We now show that for $C$ large enough, $\widetilde{V}_{\varepsilon ,\widehat{%
z}}$ is a supersolution to \eqref{huvudekvationen}, that is
\begin{equation}
D_{t}\widetilde{V}_{\varepsilon ,\widehat{z}}\left(t,x\right)+F\left(t,x,\widetilde{V}%
_{\varepsilon ,\widehat{z}}\left(t,x\right),D_{x}\widetilde{V}_{\varepsilon ,\widehat{z}%
}\left(t,x\right),D_{x}^{2}\widetilde{V}_{\varepsilon ,\widehat{z}}\left(t,x\right)\right)\geq 0,\quad
\text{for }\left(t,x\right)\in \Omega ^{0}.
\label{eq:ekvationen_f�r_den_sista_barri�ren}
\end{equation}%
With $D_{s}$ and $D_{\eta }$ denoting differentiation with respect to the
first and third arguments of $w_{\varepsilon }$, respectively, we have%
\begin{eqnarray}
D_{t}\widetilde{V}_{\varepsilon ,\widehat{z}}\left(t,x\right) &=&AD_{t}\alpha \left(t,x\right)+e^{\left(%
\widehat{C}/\chi \right)\alpha \left(t,x\right)}\frac{\widehat{C}}{\chi }D_{t}\alpha \left(
t,x\right) w_{\varepsilon }\left( t,x,z\left(t\right)\right) +e^{\left(\widehat{C}/\chi
\right)\alpha \left(t,x\right)}  \notag \\
&&\cdot \left( D_{s}w_{\varepsilon }\left( t,x,z\left(t\right)\right) -2K\widehat{%
\alpha }\left\langle D_{\eta }w_{\varepsilon }\left(t,x,z\left(t\right)\right),\widetilde{\gamma }%
\left(0,\widehat{z}\right)\right\rangle t^{\widehat{\alpha }-1}\right)  \notag \\
&&+C\widehat{\alpha }t^{\widehat{\alpha }-1}.  \label{eq:timederiv}
\end{eqnarray}%
Moreover, by \eqref{ass_F_nondecreasing} with $\lambda =0$ and by %
\eqref{ass_F_boundary} we have%
\begin{eqnarray}
&&F\left(t,x,\widetilde{V}_{\varepsilon ,\widehat{z}}\left(t,x\right),D_{x}\widetilde{V}%
_{\varepsilon ,\widehat{z}}\left(t,x\right),D_{x}^{2}\widetilde{V}_{\varepsilon ,%
\widehat{z}}\left(t,x\right)\right)  \notag \\
&\geq &F\left( t,x,g\left(\widehat{z}\right),0,0\right) -\sup_{\Omega }m_{2}\left(
|D_{x}\widetilde{V}_{\varepsilon ,\widehat{z}}\left(t,x\right)|+||D_{x}^{2}\widetilde{V}%
_{\varepsilon ,\widehat{z}}\left(t,x\right)||\right) .  \label{spacederiv}
\end{eqnarray}%
By \eqref{testlemma45}-\eqref{testlemma47}, \eqref{eq:timederiv} and %
\eqref{spacederiv}, we can find $C$ so that %
\eqref{eq:ekvationen_f�r_den_sista_barri�ren} is satisfied. Hence, using %
\eqref{F_fundamental} and Proposition 7.2. in \cite{CrandallIshiiLions1992},
$\widetilde{V}_{\varepsilon ,z}$ is a viscosity supersolution in $\Omega $
which satisfies the boundary condition \eqref{randvillkor} on $E_{\widehat{z}%
}\cap \partial \Omega $ in the viscosity sense.

We now perform the localized comparison. From the construction of $%
\widetilde{w}$, it is clear that $\widetilde{w}\left(0,x\right)\leq g\left( x\right) $,
for all $x\in \overline{\Omega }_{0}$. Combined with the left inequality in %
\eqref{eq:sup_barrier_above_g}, this yields%
\begin{equation}
\widetilde{V}_{\varepsilon ,\widehat{z}}\left(0,x\right)\geq \widetilde{w}\left(0,x\right),\quad
\text{for }x\in \overline{\Omega }_{0}\text{. }  \label{eq:bottom_comparison}
\end{equation}%
Moreover, for some constant $K$ depending on $g$, $\alpha $, $\widehat{z}$, $%
A$, $\kappa $ and $\rho $, we have%
\begin{equation*}
\widetilde{V}_{\varepsilon ,\widehat{z}}\left(t,x\right)\geq -K+\chi \frac{\left\vert
x-z\left( t\right) \right\vert ^{2}}{\varepsilon }+B,\quad \text{for}%
\;\left( t,x\right) \in \left( \partial E_{\widehat{z}}\setminus \partial
\Omega \right) \cap \left( \left[ 0,\kappa \right) \times \mathbb{R}%
^{n}\right) .
\end{equation*}%
Since $\widetilde{w}$ is bounded, we can conclude, by increasing $B$ if
necessary, that
\begin{equation}
\widetilde{V}_{\varepsilon ,\widehat{z}}\left(t,x\right)\geq \widetilde{w}\left(t,x\right),\quad
\text{for}\;\left(t,x\right)\in \left(\partial E_{\widehat{z}}\setminus \partial \Omega
\right)\cap \left( \lbrack 0,\kappa )\times
\mathbb{R}
^{n}\right) .  \label{supbarrierdominaterest}
\end{equation}%
Now, let $\kappa $ be so small that for some $\widetilde{\varepsilon }>0$,
it holds that
\begin{equation}
\left\vert x-z\left( t\right) \right\vert >\widetilde{\varepsilon }>0\quad
\text{whenever}\;\left( t,x\right) \in \left( \partial E_{\widehat{z}%
}\setminus \partial \Omega \right) \cap \left( \left[ 0,\kappa \right)
\times \mathbb{R}^{n}\right) .  \label{eq:extra_sista_z(t)_egenskap}
\end{equation}%
This choice is possible by the definition of $z\left(t\right)$ and by the properties of
the domain. Inequality \eqref{eq:extra_sista_z(t)_egenskap} implies that it
is no restriction to assume that $B\rightarrow 0$ as $\varepsilon
\rightarrow 0$, which is necessary. By means of \eqref{eq:sup_barrier_RV}, %
\eqref{eq:bottom_comparison} and \eqref{supbarrierdominaterest}, we can use
Corollary \ref{maxrand_partial} to make comparison in $E_{\widehat{z}}\cap
\overline{\Omega }$ of the supersolution $\widetilde{V}_{\varepsilon ,%
\widehat{z}}$ with each subsolution in the definition of $\widetilde{w}$.
Hence
\begin{equation*}
\widetilde{V}_{\varepsilon ,\widehat{z}}\left(t,x\right)\geq \widetilde{w}\left(t,x\right),\quad
\text{for}\;\left(t,x\right)\in \overline{E}_{\widehat{z}}\cap \overline{\Omega },
\end{equation*}%
and, as a consequence, $\widetilde{V}_{\varepsilon ,\widehat{z}}=\widetilde{V%
}_{\varepsilon ,\widehat{z}}^{\ast }\geq \widetilde{w}^{\ast }$ in $%
\overline{E}_{\widehat{z}}\cap \overline{\Omega }$. Thus, for any $x\in
\partial \Omega _{0}$,
\begin{equation*}
\widetilde{w}^{\ast }\left(0,x\right)\leq \inf_{\varepsilon ,\widehat{z}}\widetilde{V}%
_{\varepsilon ,\widehat{z}}\left(0,x\right)=g\left(x\right).
\end{equation*}

To prove that $\widetilde{w}_{\ast }$ satisfies the initial condition on $%
\partial \Omega _{0}$, we proceed similarly by constructing a family of
subsolutions of the form
\begin{equation*}
\tilde{U}_{\varepsilon ,\widehat{z}}\left(t,x\right)=g\left(\widehat{z}\right)-A\left( \alpha
\left(t,x\right)-\alpha \left(0,\widehat{z}\right)\right) -e^{\left(\widehat{C}/\chi \right)\alpha
\left(t,x\right)}w_{\varepsilon }\left(t,x,z\left(t\right)\right)-B-Ct^{\widehat{\alpha }}.
\end{equation*}%
This completes the proof of Theorem \ref{existence}. $\hfill \Box $\vspace{%
0.2cm}


\bibliographystyle{plain}
\bibliography{acompat,thomas2018}

\end{document}